\documentclass[a4paper,11pt,reqno]{amsart}
\usepackage[top=2.5cm,bottom=2.5cm,left=2.5cm,right=2.5cm]{geometry}
\usepackage[numbers,sort&compress]{natbib}
\usepackage[colorlinks]{hyperref}

\usepackage{amsthm,amsmath,amssymb,dsfont}
\usepackage{mathrsfs,amsfonts,functan,extarrows,mathtools,stmaryrd}

\usepackage{xcolor}
\usepackage{bm,scalerel}
\usepackage{enumitem}
\newtheorem{theorem}{Theorem}[section]

\numberwithin{equation}{section}
\allowdisplaybreaks
\usepackage[capitalize,nameinlink]{cleveref}
\usepackage{aliascnt}
  \newaliascnt{definition}{theorem}
  \aliascntresetthe{definition}
  \Crefname{definition}{Definition}{Definitions}
  \newaliascnt{lemma}{theorem}
  \newtheorem{lemma}[lemma]{Lemma}
  \aliascntresetthe{lemma}
  \Crefname{lemma}{Lemma}{Lemmas}
  \newaliascnt{proposition}{theorem}
  \newtheorem{proposition}[proposition]{Proposition} %
  \aliascntresetthe{proposition}
  \Crefname{proposition}{Proposition}{Propositions}
  \newaliascnt{remark}{theorem}
  \newtheorem{remark}[remark]{Remark} %
  \aliascntresetthe{remark}
  \Crefname{remark}{Remark}{Remarks} %
  \crefname{equation}{eq.}{eqs.}
\def\ls{\lesssim}

\def\p{\partial}
\def\d{\mathop{}\!\mathrm{d}}
\def\no{\nonumber}
\def\eps{\varepsilon}
\def\div{\mathrm{div}}

\def\RR{\mathbb{R}}
\def\S{\mathbb{S}}
\def\EE{\mathbb{E}}

\def\u{{\bm u}}

\def\D{\mathcal{D}}
\def\E{\mathcal{E}}
\def\H{\mathcal{H}}
\def\L{\mathcal{L}}
\def\O{\mathcal{O}}

\def\Q{\mathcal{Q}}

\def\ker{\mathrm{Ker}}

\def\IN{\mathrm{in}}
\def\M{{\!\scriptscriptstyle M}}
\def\R{{\!\scriptscriptstyle R}}
\def\MV{{\!\scriptscriptstyle {M^{{\scaleto{-1}{3pt}}}}}}

\newcommand{\lx}{{L^2_x}}
\newcommand{\lv}{{L^2_v}}
\newcommand{\lxv}{{L^2_{x,v}}}
\newcommand{\hx}[1]{{H^ {#1} _x}}
\newcommand{\hxv}[1]{{\H^ {#1} _x L^2_{v}}}

\newcommand{\abs}[1]{\left|#1\right|}
\newcommand{\nm}[1]{\left\|#1\right\|}

\newcommand{\skp}[2]{\left\langle #1,\, #2 \right\rangle}
\newcommand{\skpt}[2]{\langle #1,\, #2 \rangle}
\newcommand{\agl}[1]{\left\langle#1\right\rangle}
\newcommand{\agll}[1]{\left\langle\!\!\left\langle #1 \right\rangle\!\!\right\rangle}
\newcommand{\sskp}[2]{\left\langle\!\!\left\langle #1,\, #2 \right\rangle\!\!\right\rangle}

\begin{document}

\title[Macroscopic Limit for Kinetic Cucker-Smale Model] {From Kinetic Flocking Model of Cucker-Smale Type to Self-Organized Hydrodynamic model}

\author[Ning Jiang]
    {Ning Jiang}
\address[Ning Jiang]
    {\newline School of Mathematics and Statistics, Wuhan University, Wuhan, 430072, P. R. China}
\email{njiang@whu.edu.cn}

\author[Yi-Long Luo]
    {Yi-Long Luo${}^*$}
\address[Yi-Long Luo]
    {\newline School of Mathematics, South China University of Technology, Guangzhou, 510641, P. R. China}
\email{luoylmath@scut.edu.cn}

\author[Teng-Fei Zhang]
    {Teng-Fei Zhang}
\address[Teng-Fei Zhang]
    {\newline School of Mathematics and Physics, China University of Geosciences, Wuhan, 430074, P. R. China}
\email{zhangtf@cug.edu.cn}

\thanks{${}^*$ Corresponding author}

\maketitle
\begin{abstract}
We investigate the hydrodynamic limit problem for a kinetic flocking model. We develop a GCI-based Hilbert expansion method, and establish rigorously the asymptotic regime from the kinetic Cucker-Smale model with a confining potential in a mesoscopic scale to the macroscopic limit system for self-propelled individuals, which is derived formally by Aceves-S\'anchez, Bostan, Carrillo and Degond (2019).

In the traditional kinetic equation with collisions, for example, Boltzmann type equations, the key properties that connect the kinetic and fluid regimes are: the linearized collision operator (linearized collision operator around the equilibrium), denoted by $\mathcal{L}$, is symmetric, and has a nontrivial null space (its elements are called collision invariants) which include all the fluid information, i.e. the dimension of Ker($\mathcal{L}$) is equal to the number of fluid variables. Furthermore, the moments of the collision invariants with the kinetic equations give the macroscopic equations.

The new feature and difficulty of the corresponding problem considered in this paper is: the linearized operator $\mathcal{L}$ is not symmetric, i.e. $\mathcal{L}\neq \mathcal{L}^*$, where $\mathcal{L}^*$ is the dual of $\mathcal{L}$. Moreover, the collision invariants lies in Ker($\mathcal{L}^*$), which is called generalized collision invariants (GCI). This is fundamentally different with classical Boltzmann type equations. This is a common feature of many collective motions of self-propelled particles with alignment in living systems, or many active particle system. Another difficulty (also common for active system) is involved by the normalization of the direction vector, which is highly nonlinear.

In this paper, using Cuker-Smale model as an example, we develop systematically a GCI-based expansion method, and micro-macro decomposition on the dual space, to justify the limits to the macroscopic system, a non-Euler type hyperbolic system. We believe our method has widely application in the collective motions and active particle systems.


\end{abstract}

% \keywords{Hydrodynamic limits; Hilbert expansion; Cucker-Smale model;}

% \subjclass{35Q20, 35Q61, 76P05, 76W05}

% \maketitle

\setcounter{tocdepth}{3}
%\tableofcontents

%%----------------------------------------------------------------------
%%----------------------------------------------------------------------
\section{Introduction} % (fold)
\label{sec:introduction}

% section introduction (end)

\subsection{Problem description} % (fold)
\label{sub:problem_description}

% subsection problem_description (end)

We consider the following kinetic model of Cucker-Smale type for self-propelled individuals, as employed in \cite{ABCD-19mbe,BC-20m3as}, in a rescaled formulation:
  \begin{align}\label{eq:CS}
    \p_t f^\eps + v \cdot \nabla_x f^\eps = \frac{1}{\eps} \mathcal{Q}(f^\eps),
  \end{align}
with $f^\eps = f^\eps(t,x,v)$ denoting the one-particle density distribution function in the position-velocity phase space $(x,v) \in \RR^d \times \RR^d$ ($d \ge 2$) at time $t\ge 0$. The kinetic equation is supplemented with the following initial conditions:
  \begin{align}
    f(t,x,v)|_{t=0} = f^{\IN}(x,v) \,.
  \end{align}
The ``collisional operator'' are expressed as
  \begin{align}\label{eq:collision-op}
    \mathcal{Q}(f) = \div_v \left\{ \sigma \nabla_v f + f \left[ (v- \Omega_{f}) + \eta \nabla_v V \right] \right\}.
  \end{align}
Here the term $\div_v (f\nabla_v V)$ stands for the self-propulsion and friction mechanism leading to
the cruise speed of the particles in the absence of alignment, with a confining potential function acting only in the microscopic velocity variables $V(v)$. The diffusion coefficient $\sigma>0$ denotes the noise intensity in the velocity space, and the constant $\eta>0$ stands for an accommodation coefficient. For the sake of simplicity, both of $\sigma$ and $\eta$ are set to unit constant in the rest of this paper. In particular, we consider here the case $V= V(|v|) = \beta \tfrac{|v|^4}{4} - \alpha \tfrac{|v|^2}{2}$ with constants $\alpha,\, \beta>0$ throughout this paper. Furthermore, the term $\div_v\left\{f(v- \Omega_{f})\right\}$ stands for the relaxation toward the normalized mean velocity, which we explain in the following. 

We define respectively the macroscopic density function and the first-order moment function with respect to velocity variables, as follows:
  \begin{align}
    \rho_f (t,x) = \int_{\RR^d} f(t,x,v) \d v, \qquad
    j_f (t,x) = \int_{\RR^d} v f(t,x,v) \d v.
  \end{align}
The notation $\Omega_{f}$ is the orientation of the mean velocity with respect to $f$, as follows:
  \begin{align}
    \Omega_{f} =
  \begin{cases}
    \frac{\int_{\RR^d} v f(v) \d v}{\abs{\int_{\RR^d} v f(v) \d v}}, & \text{ if } \int_{\RR^d} v f(v) \d v \in \RR^d \setminus {0},  \\[3pt]
    0, & \text{ if } \int_{\RR^d} v f(v) \d v = 0.
  \end{cases}
  \end{align}
Note that, the kinetic model \cref{eq:CS} is a modified model, in a localized and normalized sense, comparing to the classical Cucker-Smale model \cite{CS-07ieee,CS-07jjm} and other localized versions studied before, see \cite{MT-11jsp,CFTV-10b,CFR-10sima,HT-08krm,HL-09cms,BC-20m3as}. In usual localized Cucker-Smale models, the relaxation mechanism arising from the alignment between particles is towards the mean velocity $\tfrac{\int_{\RR^d} v f(v) \d v}{\abs{\int_{\RR^d} f(v) \d v}} = \tfrac{j_f}{\abs{\rho_f}}$, appearing as $\div_v [ f ( v- \tfrac{j_f}{\abs{\rho_f}} ) ] $ or with an additional multiplier $ \rho_f \cdot \div_v [ f ( v- \tfrac{j_f}{\abs{\rho_f}} ) ] $. Whereas, the relaxation considered in \cref{eq:CS} is towards the (local and) normalized velocity orientation $\Omega_{f} = \tfrac{j_f}{\abs{j_f}}$ with unit length.

We aim in this paper to consider the large time and space scale regime, i.e., the asymptotic behavior of \cref{eq:CS} when the small parameter $\eps$ goes to zero. The parameter $\eps$ stands for the relaxation time scale, or equivalently, the asymptotic cruise speed (the balance between friction and self-propulsion), see \cite{BC-17m3as,ABCD-19mbe}. Here, the relevant scaling is a hydrodynamic scaling, which means that $\eps$ also equals the ratio of the microscopic time scale (the mean time between interactions) to the
macroscopic observation time. Corresponding to the formal derivation in \cite{ABCD-19mbe}, the limit system governs the dynamics of the density $\rho$ and mean orientation functions $\Omega$, which are functions of space and time variables, namely, $(\rho,\,\Omega) = (\rho (t,x),\,\Omega (t,x))$. This macroscopic limit system, as referred to the \emph{Self-Organized Hydrodynamics} (SOH) system, reads as follows:
  \begin{align}\label{sys:limit-0-order}
  \begin{cases}
    \p_t \rho + \div_x (c_0 \rho \Omega) = 0, \\[3pt]
    \p_t \Omega + c_1 (\Omega \cdot \nabla_x ) \Omega + \sigma (\mathrm{Id} - \Omega \otimes \Omega) \nabla_x \ln \rho = 0,
  \end{cases}
  \end{align}
where the constants $c_0,\, c_1$ will be defined specifically below, see \cref{prop:limit-0-order}. The corresponding initial conditions are
  \begin{align}
    \rho(t,x)|_{t=0} = \rho^{\IN}(x) = \int_{\RR^d} {f^{\IN}(x,v)} \d v, \quad
    \Omega(t,x)|_{t=0} = \Omega^{\IN}(x) = \tfrac{\int_{\RR^d} {vf^{\IN}(x,v) \d v}}{\abs{\int_{\RR^d} {vf^{\IN}(x,v) \d v}}}\,.
  \end{align}

The present paper focuses on the rigorous justification of the asymptotic regime from the kinetic Cucker-Smale type model \cref{eq:CS} with a confining potential towards the above self-organized hydrodynamic system \cref{sys:limit-0-order}. This is a continuation of the series of previous works \cite{JXZ-16sima,JLZ-20arma}, concerning the hydrodynamic limit of kinetic self-organized models from a viewpoint of mathematical analysis. We emphasize that previous works only consider very special cases, while this current is the first result concerning the fully general case.

%--------------------------------------------------------------------
\subsection{Brief review} % (fold)
\label{sub:brief_review}

% subsection brief_review (end)

Collective motions of self-propelled particles are ubiquitous in many disciplines, such as, in physical, biological or chemical systems. Well-known examples contains birds flocking, fish schooling, animal collective behaviors, human crowds and social dynamics, and so on. Such collective behaviors generates complexity interactions among individuals, presenting many nonlinear and nonlocal phenomena, and hence are viewed as complex system consisting of active (living) particles. Describing the complex emergence collective behaviors exhibited in multiscale levels brings new challenges from the viewpoints of phenomenological interpretation, modelling and numerical simulations, and rigorous analysis, and has already attracted more and more attentions from different areas for last decades. For more introductions, we may refer the readers to, for example, \cite{BDT-17b,BBGO-17b} or review papers \cite{MT-11jsp,KCBFL-13physD,BF-17m3as,CHL-17chap,ABF8-19m3as,Tad-21nAMS}.

Since from the pioneering works \cite{VCBCS-95prl} on the so-called Vicsek model, and \cite{CS-07ieee,CS-07jjm} on the Cucker-Smale model, alignment dynamics, in which the particles (agents) tend to align with their neighbors, has been particularly studied with a lot of literatures. This phenomenon is also described by the term \emph{flocking} \cite{ABF8-19m3as}.
% Reynolds proposed in \cite{Rey87} a general communication interaction driven by alignment, repulsion, and attractions.
In \cite{DM-08m3as}, Degond-Motsch proposed the self-organized kinetic model of Vicsek type; and a kinetic version of Cucker-Smale model was proposed by Ha-Tadmor \cite{HT-08krm}.
Inspired by these works, variants of interesting self-propelled (self-organized) models have been developed, to model more complex interactions such as repulsion, attractions, nematic alignment, suspensions, body-attitude coordination, et al., see \cite{DFLMN-14swz,DFL-15arma,DA-20m3as,DMVY-19jmfm,DFA-17m3as,DFAT-18mms} and the references therein.

Note that, the kinetic model \cref{eq:CS} is a Cucker-Smale type model, taking account for local alignment, random Brownian noise in the velocity space, and the friction and self-propulsion mechanisms. Comparing with previous treatments \cite{MT-11jsp,CFTV-10b,CFR-10sima,HT-08krm,HL-09cms}, the model \cref{eq:CS} appears to be a localized and normalized one. As indicated in \cite{ABCD-19mbe}, it does not present phase transition, which is also an important topic in the fields of active particle systems far from equilibrium. Notice Bostan-Carrillo \cite{BC-20m3as} investigated the phase transition when considering a similar model as \cref{eq:CS}, but with a different relaxation towards the mean velocity $\tfrac{j_f}{\abs{\rho_f}}$ rather than $\Omega_f$.

In this paper, we focus on the hydrodynamic limit problem, by using an appropriate asymptotic expansion scheme, from the mesoscopic kinetic Cucker-Smale model with local alignment along with the normalized mean orientation, \cref{eq:CS}, towards the self-organized hydrodynamic model \cref{sys:limit-0-order} at a macroscopic level. We aim to provide a rigorous justification on the asymptotic regime as the parameter $\eps$ tends to zero. The expansion approach, referred as the Hilbert expansion approach, was initialized from the classic work of Caflisch \cite{Caf-80cpam}, concerning the compressible Euler limit form the Boltzmann equation. This approach aims to construct a class of special solutions to the scaled kinetic equation around the solutions to the limiting equation (note, the well-posedness of the solutions to the limiting equations is assumed to be known or already proven). The Hilbert expansion approach has been proved to be effective in justifying the various hydrodynamic limit issues for kinetic and kinetic-fluids coupled equations,
to name a few, the limit from the Boltzmann equation to the incompressible Navier-Stokes and the incompressible/compressible Euler equations in \cite{MEL-89cpam,Guo-06cpam,JLT-22axv},
the limit from the Vlasov-Maxwell-Boltzmann (VMB) system to the incompressible Navier-Stokes-Fourier-Maxwell (InNSFM) system \cite{JLZ-21axv},
and the limit from the Doi-Onsager equation to the Ericksen-Leslie equation in liquid crystal theory \cite{WZZ-15cpam}, and so on.
The key ingredients for the Hilbert expansion lies in obtaining the uniform \emph{a-priori} estimates for the remainder equation, besides of proving, usually in a framework of classical solutions, the well-posedness of the limiting equation.

In our previous works, corresponding to a formal coarse-graining process in Degond-Mostch \cite{DM-08m3as} and related models \cite{DFLMN-14swz,DFL-15arma}, and micro-macro coupling models \cite{DMVY-19jmfm}, the Hilbert expansion approach is also employed to rigorously justify the hydrodynamic limit from self-organized kinetic model of Vicsek type to self-organized hydrodynamic model in \cite{JXZ-16sima}, and in \cite{JLZ-20arma} the convergence regime from a kinetic-fluid model coupling of Vicsek-Navier-Stokes model towards the self-organized hydrodynamics and Navier-Stokes equations.
We emphasize that, in these works, the notion of \emph{Generalized Collision Invariants} (GCI) employed in \cite{DM-08m3as} plays a key role in derivation and analysis of the limit regime. The GCI hypothesis, forcing the class of solutions to be constructed always present a first-moment along with the known limit orientation, restricts our the expansion ansatz to a linear case with respect to the limit orientation. Note this remains consistent with the limit regime, due to the basic fact the limit equilibrium lies in the kernel of Fokker-Planck operator, with the exactly same orientation. From a viewpoint of analysis, under the GCI hypothesis, the linear collision operator is reduced to the Fokker-Planck(-type) operator.

By contrast, we consider in the present paper the linearized operator with respect to the original nonlinear collision operator through its first-order variation, (not just a linear operator as before). The linearized operator, denoted later by $\L_{\Omega_0}$, contains the GCI structure in a natural way, see \Cref{sub:derivation_of_macroscopic_system}. The operator $\L_{\Omega_0}$ presents new complex features, in particular, it is not self-adjoint any more. This is exactly the reason we have to consider the dual operator $\L_{\Omega_0}^*$ and employ a duality solvability with respect to the orthogonal set of its kernel, we refer the readers to \cite{WZZ-15cpam} for similar ideas. The key fact causes several difficulties in the processes of both derivation and analysis, our treatment is hence different from that on standard kinetic theory like the Boltzmann equation. The formal derivation process is required to consider the interactions of kernel and its kernel orthogonal parts, presenting more complex nonlinearity effects. In this sense, we could refer the duality-based expansion method as a nonlinear version, relative to that of previous works \cite{JXZ-16sima,JLZ-20arma}.
As for the rigorous analysis, the Poincar\'e inequality applied to the remainder function is with a non vanishing mean value, which brings some singular terms of order $\tfrac{1}{\eps}$. To overcome this difficulty, we employ a new coercivity estimate (of type II, see \Cref{prop:coercivity}), by finding an additional dissipation effect from the structure of the confining potentials $V(|v|)$. In order to obtain a closed estimate, we need to choose carefully the appropriate function space that the energy functionals belongs to. For more explanations, see \Cref{sub:strategy_of_the_proof} below.

Note that, the Hilbert expansion scheme, based on the duality solvability in the present paper, contains the GCI-based expansion method studied in \cite{JXZ-16sima,JLZ-20arma} as a special case on one hand; and on the other hand, it could be extended to higher-orders expansion, as stated in \Cref{sub:induction_scheme_for_higher_order_expansions}. In this sense, this paper aims to provide a general framework, from a viewpoint of rigorous mathematics, to justify the passage from mesoscopic/kinetic scale to macroscopic scale for collective motions of self-propelled particles with alignment.

As a last comment, there are two basic methodology to derive and justify the fluid limits from kinetic equations, one is the Hilbert expansion method mentioned above (also referred as a \emph{top-down} strategy), and the other one the moment method (correspondingly, referred as a \emph{bottom-up} strategy). The most successful achievement by the moment method is the so-called Bardos-Golse-Levermore (BGL) program. The BGL program justifies the Leray solutions of the incompressible Navier-Stokes equations from the DiPerna-Lions renormalized solutions to the Boltzmann equation, which is initialized from \cite{BGL-91jsp,BGL-93cpam} and with the first complete convergence result given by Golse and Saint-Raymond \cite{GSR-04inv}.
As a comparison with the top-down result \cite{JLZ-21axv} concerning the limit from VMB to InNSFM, we refer to Ars\'enio and Saint-Raymond \cite{ASR-19b}, and Jiang-Luo \cite{JL-22apde}, both of which belong to the bottom-up type, in the contexts of renormalized solutions and classical solutions, respectively. We also refer to \cite{LW-18jfa} as a comparison with \cite{WZZ-15cpam}, relating to the liquid crystal theory.
In particular, Degond-Frouvelle-Liu \cite{DFL-22krm} studied the convergence from the Doi-Navier-Stokes coupled model for rod-like polymer suspensions to the Ericksen-Leslie model in liquid crystal theory, in a bottom-up strategy associating to the GCI concept. Usually, when the alignment force is considered, the conservation of momentum could not be preserved. Therefore, it is a very challenging topic on the hydrodynamic limit of the kinetic Vicsek and Cucker-Smale type models, and more general complex systems of active particles, requiring more attentions and studies.

It also should be mentioned, as for the Cucker-Smale and related modified models, there are many other research topics in analysis, such as well-posedness issues, emergent behavior analysis, justification of the mean field limits from a microscopic particle model to the mesoscopic kinetic model, sensitivity analysis with random inputs, and so on, see \cite{CHL-17chap,ABF8-19m3as,BDMa-22sima,FK-19apde} and references therein for instance.

%--------------------------------------------------------------------
% \subsection{Notations} % (fold)
% \label{sub:notations}

% % subsection notations (end)
\subsection{Main results} % (fold)
\label{sub:overview_of_main_results}

% subsection overview_of_main_results (end)

Before we state the main results, we firstly introduce some notations.

\smallskip\noindent\underline{\textbf{Notations}}.
Let $\alpha=(\alpha_1,\, \alpha_2,\, \alpha_3) \in \mathbb{N}^3 $ be a multi-index with its length defined as $\textstyle |\alpha| = \sum_{i=1}^3 \alpha_i $. We define here the multi-derivative operator $\nabla^\alpha_x = \p_{x_1}^{\alpha_1} \p_{x_2}^{\alpha_2} \p_{x_3}^{\alpha_3}$, and we also denote $\nabla_x^k$ for $|\alpha| = k$ for simplicity. In addition, the notation $A \ls B$ means that there exists some positive constant $C>0$ such that $A \le CB$.

We will work in Sobolev spaces with respect to $x$, and in weighted Sobolev spaces with respect to the microscopic velocity variables $v$.
For simplicity, we denote by $\nm{\cdot}_{L^2_x}$ the usual $L^2$-norm over spatial variables $x$, and by $\nm{\cdot}_\hx{s}$ the higher-order spatial derivatives $H^s$-norm. Besides, we denote by $\nm{\cdot}_{L^2_{x,v}}$ the mixed $L^2$-norm with respect to both variables $x$ and $v$.

The notation $\langle \cdot,\, \cdot \rangle$ will stand for the inner product in a specific space marked by a subscript, for instance, by a $L^2_v$ or $L^2_{x}$. Sometimes we will omit the subscript of space marks for simplicity, when no confusion is possible. The angle bracket $\langle \cdot \rangle$ denotes the integral over the velocity space, i.e., $\langle f \rangle = \int_{\mathbb{R}^d} f \d v$. Furthermore, the double angle brackets $\langle\langle \cdot,\, \cdot \rangle\rangle$ stands for the inner product in $L^2_{x,v}$.

On the other hand, we introduce the weighted Sobolev spaces over velocity variables $v$. The weighted $L^2$ space with a weight $M(v)$, denoted by $L^2_v(M)$, is defined as the set of measurable function $f$ satisfying $\int_{\RR^d} |f|^2 M \d v < \infty$. Similarly, the higher-order weighted Sobolev spaces $H^s_v(M)$ consists of those measurable function such that $\int_{\RR^d} |\nabla_v^k f|^2 M \d v < \infty$ for all integers $k \in [0,\,s] $.

By defining the following weighted inner product in $L^2_v$:
  \begin{align*}
    \skp{f}{g}_{L^2_v(M)} \triangleq \left \langle f,\, g M \right\rangle
      = \int_{\mathbb{R}^d} fg \,M \d v,
  \end{align*}
for any pairs $f(v),\, g(v) \in L^2_v(M)$. The associating weighted norm can be defined, and we denote it by the notation $\nm{\cdot}_{L^2_v(M)}$.
Note, the weight we will use in sequels usually appears as $M_{\Omega_0}$ (or $M_{\Omega_0}^{-1}$), which are sometimes simplified as a single subscript $M$ (or $M^{-1}$, respectively), when no possible confusion arises.

We will use the mixed weighted Sobolev spaces with the norm $\nm{\cdot}_{\hxv{s}(M^{-1})}$ for $s\ge 0$:
  \begin{align}
    \nm{f}^2_{\hxv{s}(M^{-1})}
    = \sum_{0\le k \le s} \iint_{\RR^{2d}} \tfrac{\abs{\nabla_x^k f}^2}{ M_{\Omega_0} } \d v \d x
    = \sum_{0\le k \le s} \nm{ \tfrac{\nabla_x^k f}{M_{\Omega_0}} }_\M^2 \,.
  \end{align}
We mention that the weighted Sobolev norm $\nm{\cdot}_{\hxv{s}(M^{-1})}$ is not the same as the usual mixed Sobolev norm, since the weight $M_{\Omega_0}$ is a function over $x$.

Before we state the main theorem of this paper, we introduce some notations, for which we will explain in details later. 
\begin{itemize}
  \item The the von Mised-Fisher (VMF) distributions $f_0(t,x,v) = \rho_0(t,x) M_{\Omega_0(t,x)}(v)$, where $(\rho_0, \Omega_0) \in C([0,T_0];\, H^{12}_x(\mathbb{R}^3))$ is a local solution to the Cauchy problem of the hydrodynamic system \cref{sys:limit-0-order} with given initial data $(\rho_0^{\IN}, \Omega_0^{\IN}) \in H^{12}_x(\mathbb{R}^3)$ and $\rho_0^{\IN}>0$;
      
  \item $f_1^{\shortparallel}(t,x,v) = f_1^{\shortparallel}(\rho_1,\u_1,\Omega_0) = \rho_1(t,x) M_{\Omega_0} + c_0^{-1} v_{\bot} \cdot {\bm u}_1(t,x) M_{\Omega_0}$, where $(\rho_1, \u_1) \in C([0,T_1];$ $H^{9}_x(\mathbb{R}^3))$ is a local solution to the Cauchy problem of the first order macroscopic system \cref{sys:limit-1st-order-rho}-\eqref{sys:limit-1st-order-j} with given initial data $(\rho_1^{\IN}, \Omega_1^{\IN}) \in H^{9}_x(\mathbb{R}^3)$;

  \item $f_i^{\perp}(t,x,v) \in \ker^{\bot}(\L_{\Omega_0})$ ($i=1,2$) are respectively the unique solution to \cref{eq:1st-order} and \cref{eq:2nd-order}, satisfying $f_1^{\perp} \in H^1([0,T];\, \H^{9}_x L^2_v(M^{-1}))$ and $f_2^{\perp} \in H^1([0,T];\, \H^{6}_x L^2_v(M^{-1}))$.
\end{itemize}

% In addition, the notation $A \ls B$ means that there exists some positive constant $C>0$ such that $A \le CB$.
% The notation $A \sim B$ stands for the equivalence between both sides up to a constant, i.e., some $C>0$ exists such that $C^{-1}B \le A \le CB$. Sometimes we will also omit the integral domain symbol for simplicity.

%% ------------------------------------------------------------------------
% \subsection{Main results} % (fold)
% \label{sub:main_results}

% % subsection main_results (end)

%% ------------------------------------------------------------------------ %%
% \iffalse
% \subsubsection{Main result} % (fold)
% \label{ssub:main_result}

% % subsubsection main_result (end)
\medskip\noindent\underline{\textbf{Main result}}.
We now state our main result of the hydrodynamic limit, for simplicity we set the spatial dimension $d=3$.
\begin{theorem}[Hydrodynamic limit] \label{thm:limit}
  % Let $s\ge 2$. $s=3$
Assume that the coefficients $\alpha,\beta$ in the confining potential $V(|v|)= \beta \tfrac{|v|^4}{4} - \alpha \tfrac{|v|^2}{2}$ satisfy $\alpha > 4^3 (c_0^{-1} \lambda_0^{-1} \beta^{\frac{1}{2}} (d+2)^\frac{3}{2} +1)$, with $\lambda_0$ being the Poincar\'e constant in Fokker-Planck operator.

% {eq:decomps-f1}

  Denoted by $f_1 (t,x,v)= f_1^{\shortparallel} (\rho_1,\u_1,\Omega_0) + f_1^{\perp}(t,x,v)$, and $f_2 = f_2^{\perp}(t,x,v)$. Let the initial data being $f_1 (t,x,v)|_{t=0} $ $= f_1^{\shortparallel} (\rho_1^{\IN},\u_1^{\IN},\Omega_0^{\IN}) + f_1^{\perp}(t,x,v)|_{t=0}$, and $f_2 (t,x,v)|_{t=0} = f_2^{\perp}(t,x,v)|_{t=0}$, both of which are denoted by $f_1^{\IN}(x,v)$ and $f_2^{\IN}(x,v)$, respectively.

  Suppose the initial data obeys the following formulation:
    \begin{align}\label{eq:ini-data}
      f^{\eps,\IN}(x,v) = f_0^{\IN}(x,v) + \eps f_1^{\IN}(x,v) + \eps^2 f_2^{\IN}(x,v) + \eps^2 f_\R^{\eps,\IN}(x,v),
    \end{align}
  with a bounded initial data $\|f_\R^{\eps,\IN}(x,v)\|_{\H^3_xL^2_v(M^{-1})} \le C$.

  Then there exists an $\eps_0 >0$ such that, for all $\eps \in (0,\eps_0)$, the Cauchy problem of Cucker-Smale system \eqref{eq:CS} admits a unique solution $f^\eps \in L^\infty ([0,T_1];\ \H^3_xL^2_v(M^{-1}))$, being of the form
  \begin{align}
    f^\eps(t,x,v) = f_0(t,x,v) + \eps f_1(t,x,v) + \eps^2 f_2(t,x,v) + \eps^2 f_\R^{\eps}(t,x,v),
  \end{align}
  with $\sup_{t\in[0,T_1]} \|(f_\R^\eps,\, \eps^{\frac{1}{2}} \nabla_x f_\R^\eps,\, \eps \nabla_x^2 f_\R^\eps,\, \eps^{\frac{3}{2}} \nabla_x^3 f_\R^\eps)(t)\|_{L^2_{x,v}(M^{-1})} \le C$, where the bound $C$ is independent of $\eps$.% and $t\in[0,T]$.
\end{theorem}

% \begin{remark}\label{rem:for-thm}

Note, the main result concerns only a truncated expansion formula up to order 2, more general higher-order expansion formula could also be considered in the same manner, cf. \Cref{sub:induction_scheme_for_higher_order_expansions}.

We also point out here, the initial data of macroscopic parameters $(\rho_0, \Omega_0)$ and $(\rho_1, \u_1)$, and of the remainder $f_\R^{\eps}$, are the only quantities to be \emph{given}. The functions $f_i^{\perp}$'s ($i=1,2$) are determined through \cref{eq:1st-order} and \cref{eq:2nd-order}, and hence their initial data are \emph{taken} as the value of that in $t=0$. Thus, we refer the initial expansion formula \cref{eq:ini-data} as a well-prepared formula.
% \end{remark}

The crucial part is to prove a uniform-in-$\eps$ estimate on the remainder equation \cref{eq:Rem-order}. We mention that, the remainder $f_\R^\eps$ is of order $\O(\eps^2)$, indicating that the convergence $f^\eps \to f_0 = \rho_0 \Omega_0 $ stated in our main result \Cref{thm:limit} holds valid in $H^3_x$-spaces with respect to spatial variables. If only consider the remainder itself, we can have a uniform-in-$\eps$ estimate in any $H^s_x$-spaces with $s \ge 3$.

\begin{proposition}[Uniform-in-$\eps$ estimate]\label{prop:uniform-esm}
  Define the energy and dissipation functionals as follows, for $s\ge 3$:
  \begin{align}\label{eq:functionals-ED}
    \E(t) = \sum_{\mathclap{0\le k\le s}} \eps^k \nm{\tfrac{\nabla_x^k f_\R^\eps}{M_0}}_{\lxv(\M)}^2, \
    \D(t) = \sum_{\mathclap{0\le k\le s}} \eps^{k-1} \left( \nm{\nabla_v (\tfrac{\nabla_x^k f_\R^\eps}{M_0})}_{\lxv(\M)}^2 + \frac{\alpha}{c_0} \nm{P_{\Omega_0^\bot} {\nabla_x^k j_\R^\eps}}_{L^2_x}^2 \right).
  \end{align}
  Then there exist some constants $\eps_0>0$ such that, for any $\eps \in (0,\ \eps_0)$ and $t \in [0,\ T_1]$, we have
  \begin{align}\label{esm:apriori-uniform}
    \tfrac{\d}{\d t} \E(t) + \D(t) \le C \left( 1 + \E(t) + \eps^\frac{1}{2} \E^2(t) + \eps^\frac{5}{2} \E^3(t) \right) \,,
  \end{align}
  where $C$ depends only on the initial data $\|(\rho_0^{\IN}, \Omega_0^{\IN})\|_{H^{s+9}_x}$.
\end{proposition}
% \fi

%% ------------------------------------------------------------------------ %%
\subsection{Strategy of the proof} % (fold)
\label{sub:strategy_of_the_proof}

% subsection strategy_of_the_proof (end)

Our analysis on the asymptotic behavior of the kinetic Cucker-Smale model \cref{eq:CS} relies on the Hilbert expansion scheme. As pointed out in \cite{ABCD-19mbe}, the key observation lies in the linear collision operator is not an self-adjoint operator, i.e., $\L_{\Omega_0} \neq \L_{\Omega_0}^*$. By \Cref{prop:operator-L*} below, the collision invariants do not lie in the kernel of the linear operator, but in the kernel of its dual operator $\ker(\L_{\Omega_0}^*)$. This is a fundamental difference from previous researches on classical hydrodynamic limit issues, where the set of collision invariants is shown to be consistence with the kernel of linear Boltzmann operator. We sketch here our treatments from viewpoints of both derivation and analysis.

From a derivation viewpoint, we consider in this paper a macro-micro decomposition in a dual sense, and moreover, the duality solvability, see \Cref{lemm:solvability} in the next section. This property provides a way for us to construct an expansion ansatz, and to derive, at a formal level, the macroscopic limit equation associated to the leading order function, the next-order macroscopic equation, and the remainder equation, as stated in the main convergence result \Cref{thm:limit}.
Note that, the high nonlinearity exhibited in the normalized mean orientation brings an additional nonlinear term in deriving the next-order asymptotic macroscopic model. This is also different from classical results for kinetic equations by a standard Hilbert expansion approach, in which the next-order asymptotic model appears as a linear equation with a similar structure as the limit system.

On the other hand, the rigorous analysis mainly relies on obtaining the uniformly \emph{a-priori} estimate for the remainder equation, with respect to the singular parameter $\eps$. However, the two parts coming from macro-micro decomposition of the remainder function are not directly orthogonal in the sense of the usual inner product. Indeed, the macro-micro decomposition is orthogonal in a dual sense, with respect to the kernel and kernel orthogonal of the dual linear operator $\ker(\L_{\Omega_0}^*)$. This indicates that, we have to consider the estimates for the remainder function itself instead of the two parts from its macro-micro decomposition.

In order to obtain the uniform-in-$\eps$ estimates, a lower bound estimate for the linear operator will play a key role. As for the next-order functions, precisely, for the kernel orthogonal part $f_i^\perp$'s with $i=1,2$, the Poincar\'e inequality suffices to get the estimate, by noticing the mean value vanishes in the kernel of Fokker-Planck operator. However, a mean value would present in dealing with the Poincar\'e inequality for the remainder function $f_\R^\eps$ itself. For that, we employ a structure condition for the coefficients $\alpha$ and $\beta$ in the confining potential function, see the assumption (4) in \Cref{thm:limit}. This contributes to providing an additional dissipation, and hence enables us to obtain a delicate coercivity estimate for the linear operator with respect to $f_\R^\eps$ itself. We referred the two types of coercivity estimate as, respectively, type I (without any condition), and type II (with conditions on $\alpha,\beta$). The details are postponed to \Cref{sub:two_types_of_coercivity_estimate} below. Besides, the structure of the confining potential also yields a clear control on moments integral over the microscopic variables $v$, in the spirit of the proof of the Poincar\'e inequality \cref{eq:Poinc-ineq-Phi} in \Cref{lemm:Poincare-ineq}. This leads us to closing the estimate in an appropriate weighted Sobolev space, without any additional estimates on moment integrals, like in \cite{LLZ-07cpam,JLZ-18sima} to control the growing of $v \in \RR^d$.

We are still required to estimate the higher-order derivatives on spatial variables $x\in \RR^3$. A commutator arises in the interaction between $\nabla_x^s$ and $\L_{\Omega_0}$, which brings some energy terms with a singularity, see the factor $\nm{\nabla_x^{s_1} f_\R^\eps}_{\MV}$ in \cref{eq:coercive-HiOd} for instance. As explained above, the Poincar\'e inequalities in \Cref{lemm:Poincare-ineq} does not work here
because of the non-vanishing mean value. To overcome this difficulty, we employ a multiplier $\eps$ on every orders of derivatives estimate, and adjust the powers of $\eps$ to match different derivative orders, see, for instance, \cref{eq:coercive-HiOd-nu,eq:coercive-HiOd-eps} and the statements around them. Therefore, by constructing modified energy and dissipation functionals with multiplier $\eps$, see their definitions in \cref{eq:functionals-ED}, we can obtain finally a closed \emph{a-priori} estimate. Since the energy functional depends on $\eps$, it is reasonable that the convergence for the expansion ansatz holds valid in a space, with its index also depending on $\eps$, as stated in the main result \Cref{thm:limit}. As a comparison, in our previous paper \cite{JLZ-20arma}, we consider the hydrodynamic limit for a self-organized kinetic-fluid coupling model in a simpler ``linear'' case by the GCI-based expansion formula, where the Poincar\'e inequalities without a mean value holds, and thus, the convergence result holds valid in any Sobolev paces with index $s \ge 2$.

Note in the process, the high nonlinearity of mean orientation $\Omega_{f} = \tfrac{j_f}{\abs{j_f}}$ requires more careful treatment as well. The expansion ansatz separates the numerator and denominator factors by different orders, and we have to deal with the Taylor formula with remainders, see \Cref{lemm:Taylor-remainders} below. Another key ingredient is to keep the denominator $\abs{j_f}$ away from zero, which could be obtained by deriving the lower bound of $\nm{j_\eps}_{L^\infty_x}$. We prove a persistence property for the lower bound of moment function, which depends on its initial data $j_\eps^{\IN}$, see \Cref{prop:nonzero-j} in \Cref{ssub:persistence_of_moment}. It should be mentioned, the persistence estimate is consistence with the uniformly \emph{a-priori} estimate above for the remainder function, by setting a proper iteration scheme in constructing the local solution for the remainder equation, see \Cref{sub:local_existence_of_remainder_equation}.

%--------------------------------------------------------------------
\subsection{Organization of this paper} % (fold)
\label{sub:organization}

% subsection Organization of this paper (end)

The rest of this paper is organized as follows: the next section is devoted to the formal derivation, by using an order analysis from the Hilbert expansion ansatz, on the limit self-organized hydrodynamic system \cref{sys:limit-0-order}, and moreover, on the first-order macroscopic system \cref{sys:limit-1st-order-rho,sys:limit-1st-order-j}, and general higher-order asymptotic system by an induction scheme. The duality solvability plays a crucial role in the process. By truncating the expansion formula up to second order with a remainder, see \cref{eq:expansion-truncate} below, we obtain the remainder equation \cref{eq:Rem-order}.

\Cref{sec:proof_of_the_main_result} contains the main proof for the convergence result \Cref{thm:limit}. Firstly, we prove two types of coercivity estimates in \Cref{sub:two_types_of_coercivity_estimate}, in which the type II estimate (with a structure condition) contributes to providing an additional dissipation and to closing the \emph{a-priori} estimate. Secondly, some useful lemmas, concerning the Poincar\'e inequalities and estimates on remainders in Taylor's formula, are collected in \Cref{sub:some_useful_lemmas}. Thirdly, we complete in \Cref{sub:completion_of_main_result} the proof of main convergence result \Cref{thm:limit}, by obtaining the uniformly in $\eps$ estimate for the remainder equation, through \Cref{sub:uniform-esm}. Besides, the persistence of moment function $j_{f^\eps}$ is shown in \Cref{prop:nonzero-j}, providing a positive lower bound away from zero.

For the completeness of the proof, a local existence result for remainder equation remains to be proven. This is provided in \Cref{sub:local_existence_of_remainder_equation}. In addition, the estimates on the lower-orders functions, $f_1$ and $f_2$, are also required, which are given in the last \Cref{sec:estimates_on_lower_order_functions}.

% In \ref{sec:lemmas}, some useful lemmas are presented, involving the structure of the Fokker-Planck operator and the weighted Poincar\'e inequalities. \S \ref{sec:a_priori_estimates} is devoted to closing the \emph{a-priori} estimate, see Proposition \ref{prop:a-priori}, for that we need to derive closed estimates for fluctuations, including purely spatial higher-order derivatives on fluctuations themselves and fluctuation moments, mixed spatial-configurational derivatives, and furthermore, additional dissipation on number densities at the macroscopic level. In the last section \S \ref{sec:completion_of_proof}, we complete the proof of global existence result (Theorem \ref{thm:main}), by combining Proposition \ref{prop:a-priori}, the basic energy-dissipation law, and a standard bootstrap principle of continuity argument together.

%--------------------------------------------------------------------
\section{Formal Analysis from the Hilbert Expansion Method} % (fold)
\label{sec:formal_analysis_from_the_hilbert_expansion_method}

% section formal_analysis_from_the_hilbert_expansion_method (end)

\subsection{The asymptotic expansion ansatz and the leading order} % (fold)
\label{sub:the_asymptotic_expansion_ansatz_and_the_leading_order}

% subsection the_asymptotic_expansion_ansatz_and_the_leading_order (end)

We aim here to justify the asymptotic limit by performing a Hilbert expansion:
  \begin{align}\label{eq:expansion-ansatz}
    f^\eps = f_0 + \eps f_1 + \eps^2 f_2 + \cdots.
  \end{align}
As a consequence, the collision operator $\mathcal{Q}(f)$ can be expanded as:
  \begin{multline}
    \mathcal{Q}(f^\eps)
    = \div_v \left\{ \nabla_v (f_0 + \eps f_1 + \eps^2 f_2 + \cdots) \right. \\
      \left.
      + (f_0 + \eps f_1 + \eps^2 f_2 + \cdots) \cdot \left[ (v- \Omega_0 - \eps \widehat{\Omega}_1 - \eps^2 \widehat{\Omega}_2 - \cdots) + \nabla_v V(v) \right] \right\}.
  \end{multline}

Inserting this expansion formula into the kinetic Cucker-Smale model, \cref{eq:CS}, we can write that
  \begin{align}\label{eq:expansion-orders}
    (\p_t + v \cdot \nabla_x) & (f_0 + \eps f_1 + \eps^2 f_2 + \cdots) \no \\
    & = \frac{1}{\eps} \div_v \left\{ \nabla_v (f_0 + \eps f_1 + \eps^2 f_2 + \cdots) \right. \\ \no
    & \qquad \qquad \quad \left.
      + (f_0 + \eps f_1 + \eps^2 f_2 + \cdots) \cdot \left[ (v- \Omega_0 - \eps \widehat{\Omega}_1 - \eps^2 \widehat{\Omega}_2 - \cdots) + \nabla_v V(v) \right] \right\}.
  \end{align}
Here we have used the Taylor expansion formula for mean orientation:
  \begin{align}\label{eq:taylor-Omega}
    \Omega_{f^\eps} = \frac{j_{f^\eps}}{\abs{j_{f^\eps}}} = \frac{\int_{\RR^d} v f^\eps(v) \d v}{\abs{\int_{\RR^d} v f^\eps(v) \d v}} % \\[5pt] \no &
    = \Omega_0 + \eps \widehat{\Omega}_1 + \eps^2 \widehat{\Omega}_2 + \cdots,
  \end{align}
where we can compute
  \begin{align}\label{eq:expansion-Omega}
  \begin{cases}
    \Omega_0 = \frac{j_{f_0}}{\abs{j_{f_0}}} = \frac{j_0}{\abs{j_0}} , \\[5pt]
    \widehat{\Omega}_1^i = \left( \delta_{ik} - \Omega_0^i \otimes \Omega_0^k \right) \frac{j_1^k}{\abs{j_0}}
                       = P_{\Omega_0^\bot}^{ik} \frac{j_1^k}{\abs{j_0}}
                       = (c_0\rho_0)^{-1} P_{\Omega_0^\bot}^{ik} {j_1^k}, \\[7pt]
    \begin{aligned}
      \widehat{\Omega}_2^i & = P_{\Omega_0^\bot}^{ik} \frac{j_2^k}{\abs{j_0}} - \frac{j_1 \cdot j_1}{\abs{j_0}^2} \frac{j_0^i}{\abs{j_0}} - 2 \frac{j_0 \cdot j_1}{\abs{j_0}^2} \frac{j_1^i}{\abs{j_0}} + 3 \frac{j_0^i \otimes j_0^k}{\abs{j_0}^2} \frac{j_0 \cdot j_1}{\abs{j_0}^2} \frac{j_1^k}{\abs{j_0}} \\[7pt]
    & = (c_0\rho_0)^{-1} P_{\Omega_0^\bot}^{ik} {j_2^k}
        - (c_0\rho_0)^{-2} (P_{\Omega_0^\bot}^{kl} : j_1^k \otimes j_1^l) \Omega_0^i
        - 2 (c_0\rho_0)^{-2} ({j_1^l} \Omega_0^l) \cdot (P_{\Omega_0^\bot}^{ik} {j_1^k}) \,.
    \end{aligned}
  \end{cases}
  \end{align}
In order to ensure the original information of the mean orientation will not be lost in the above Taylor expansion, we have to consider remainder of $n$-th order truncation scheme in the following form:
  \begin{align}\label{eq:taylor-Omega-n}
    \eps^n \widehat{\Omega}_\R^\eps
    = \Omega_{f^\eps} - (\Omega_0 + \eps \widehat{\Omega}_1 + \eps^2 \widehat{\Omega}_2 + \cdots + \eps^n \widehat{\Omega}_n).
  \end{align}

Therefore, the balance of orders $\mathcal{O}(\tfrac{1}{\eps})$ in \cref{eq:expansion-orders} gives us that,
  \begin{align}
    \div_v \left\{ \nabla_v f_0 + f_0 \left[ (v- \Omega_0) + \nabla_v V(v) \right] \right\} = 0,
  \end{align}
or say,
  \begin{align}\label{eq:0-order}
    \mathcal{Q} (f_0) = 0,
  \end{align}
which means, from the definition of equilibrium for $\mathcal{Q}$, that the zeroth order (leading order) expansion of $f^\eps$ can be expressed through the von Mised-Fisher (VMF) distributions, i.e.,
  \begin{align}
    f_0 (t,x,v) = \rho_0 (t,x) M_{\Omega_0(t,x)}(v),
  \end{align}
with parameters $\rho_0 =\rho_{f_0}$ and $\Omega_0 = \Omega_{f_0}$ being the density and mean orientation, respectively. Here we have used the definition of a probability function $M_{\Omega}(v)$, that
  \begin{align}
    M_{\Omega}(v) = \tfrac{1}{Z_{\Omega}} e^{- \frac{\Phi_{\Omega}(v)}{\sigma} }, \text{\quad with } \Phi_{\Omega}(v) = \tfrac{\abs{v- \Omega}^2}{2} + V(|v|).
  \end{align}
Here the normalized constant $Z_{\Omega} = \int_{\RR^d} e^{- \frac{\Phi_{\Omega}(v)}{\sigma} } \d v$, and $\Omega \in \S^{d-1}$.

% %%--------------------------------------------------------------------
% \section{Formal Analysis from the Hilbert Expansion Method} % (fold)
% \label{sec:formal_analysis_from_the_hilbert_expansion_method}

% % section formal_analysis_from_the_hilbert_expansion_method (end)

%%--------------------------------------------------------------------
\subsection{Solvability and the derivation of macroscopic system} % (fold)
\label{sub:derivation_of_macroscopic_system}

% subsection derivation_of_macroscopic_system (end)
We collect here some facts on the structures of linear operator $\L_{\Omega_0}$ and its dual operator $\L_{\Omega_0}^*$, see \cite{ABCD-19mbe} for more details. Here $\L_{\Omega_0}$ is the linearization of $\mathcal{Q}$ around the equilibrium von Mised-Fisher distributions $f_0 = \rho_0 M_{\Omega_0}$.  Note the formulation may be slightly different, due to the sake for our analysis later. Besides, we point out that, the coefficients $c_0$ and $c_1$ used in this paper, are denoted in \cite{ABCD-19mbe} by $c_1$ and $c_2$, respectively.
In fact, the operator $\L_{\Omega_0}$ can be defined by the first variation of $\mathcal{Q}$ with respect to $f$, i.e. the critical point of $\mathcal{Q}$ at $f_0$:
    \begin{align}
      \mathcal{L}_{\Omega_0} g & = \frac{\d}{\d s}\Big|_{s=0} \mathcal{Q}(f_0 + s g) \\ \no
       & = \div_v \left\{ \nabla_v g + g \left[ (v- \Omega_0) + \nabla_v V(v) \right] - f_0 P_{\Omega_0^\bot} \frac{j_g}{\abs{j_0}} \right\} \\ \no
       % & = \div_v \left\{ \nabla_v g + g \nabla_v \Phi_{\Omega_0}(v) - f_0 P_{\Omega_0^\bot} \frac{j_g}{\abs{j_0}} \right\} \\ \no
       & = \div_v \left\{ \nabla_v g + g \nabla_v \Phi_{\Omega_0} - \rho_0 M_{\Omega_0} (c_0\rho_0)^{-1} \cdot P_{\Omega_0^\bot} {j_g} \right\} \\ \label{eq:linearized-op-sep}
       & = \div_v \left\{ M_{\Omega_0} \nabla_v \left( \frac{g}{M_{\Omega_0}} \right) - c_0^{-1} P_{\Omega_0^\bot} {j_g}\, M_{\Omega_0} \right\}
       \\ \no
       & = \div_v \left\{ M_{\Omega_0} \nabla_v \left( \frac{g}{M_{\Omega_0}} \right) - c_0^{-1} \nabla_v \left( v \cdot P_{\Omega_0^\bot} {j_g} \right) \, M_{\Omega_0} \right\}
       \\\label{eq:linearized-op}
       & = \div_v \left\{ M_{\Omega_0} \nabla_v \left( \frac{g}{M_{\Omega_0}} - c_0^{-1} P_{\Omega_0^\bot} v \cdot P_{\Omega_0^\bot} {j_g} \right) \right\},
    \end{align}
where $c_0 = \int_{\RR^d} (v \cdot \Omega_0) M_{\Omega_0} (v) \d v $. It implies $j_0 = \int_{\RR^d} v f_0 (v) \d v = c_0 \rho_0 \Omega_0$, and $\abs{j_0} = c_0 \rho_0 $.

Therefore, the kernel of the linearized collision operator $\mathcal{L}_{\Omega_0}$ is spanned by $\{1, P_{\Omega_0^\bot} v \}$, and is described as:
  \begin{align}\label{Ker-L-Omega0}
    \ker (\L_{\Omega_0})
      & = \left\{ g = \rho_g M_{\Omega_0} + c_0^{-1} \left( P_{\Omega_0^\bot} v \cdot P_{\Omega_0^\bot} {j_g} \right) M_{\Omega_0} \,\Big|\, \rho_g \in \RR, \text{ and } j_g \in \RR^d \right\} \no \\
      & = \left\{ g = \rho_g M_{\Omega_0} + c_0^{-1} \left( v_{\bot} \cdot {\tilde{j}_g} \right) M_{\Omega_0} \,\Big|\, \rho_g \in \RR, \text{ and } \tilde{j}_g \in {(\RR \Omega_0)^\bot} \right\}.
  \end{align}
Here we denoted $v_{\bot} = P_{\Omega_0^\bot} v$ for short. Note that $ \tilde{j}_g = P_{\Omega_0^\bot} j_g = \int_{\RR^d} P_{\Omega_0^\bot} v g M_{\Omega_0} (v) \d v = \agl{v_{\bot} g}_{M}\,.$

It is worth mentioning, the formulation of $\mathcal{L}_{\Omega_0}$ contains the GCI constraint as a natural part, by observing the factor $P_{\Omega_0^\bot} {j_g}$ forces any function with a direction parallel to the known $\Omega_0$ to be vanishing.

\begin{proposition}[\cite{ABCD-19mbe}] \label{prop:operator-L*}
  The formal adjoint of linearized operator $\mathcal{L}_{\Omega_0}$, denoted by $\mathcal{L}_{\Omega_0}^*$, is given by
    \begin{align}\label{eq:kernel-L*}
      \mathcal{L}_{\Omega_0}^* \psi
       & = \frac{\div_v \left( M_{\Omega_0} \nabla_v \psi \right)}{M_{\Omega_0}} + P_{\Omega_0^\bot} v \cdot W[\psi] \\ \no
       & = M_{\Omega_0}^{-1} \div_v \left( M_{\Omega_0} \nabla_v \psi \right) + c_0^{-1} v_{\bot} \cdot \agl{\nabla_v \psi}_{M},
    \end{align}
where we have used the notations of $\agl{g}_{M} = \int_{\RR^d} g M_{\Omega_0} \d v$, and
  \begin{align*}
    W[\psi]
    = \frac{\int_{\RR^d} M_{\Omega_0} \nabla_v \psi \d v}{\int_{\RR^d} (v \cdot \Omega_0) M_{\Omega_0} (v) \d v} = c_0^{-1} \agl{\nabla_v \psi}_{M}.
  \end{align*}

We also have the following facts:
\begin{enumerate}[label = (\arabic*)]
  \item The kernel of the adjoint operator $\mathcal{L}_{\Omega_0}^*$ is spanned by 1 and $\chi \left( \tfrac{v \cdot \Omega_0}{|v|},\, |v| \right) \tfrac{P_{\Omega_0^\bot} v}{\abs{P_{\Omega_0^\bot} v}} $ $= \chi(v) \cdot \tfrac{v_{\bot}}{\abs{v_{\bot}}} $.

  \item The statement of $\psi (v)$ being a collision invariant, is equivalent to say, $\psi (v) \in \ker (\L_{\Omega_0}^*)$ holds for some $W \in \ker \left( \mathcal{M}_{\Omega_0} - c_0 \mathrm{Id} \right)$, with a pressure tensor defined in \cite{ABCD-19mbe}: $\mathcal{M}_{\Omega_0} = \int_{\RR^d} (v- \Omega_0) \otimes P_{\Omega_0^\bot} (v- \Omega_0) M_{\Omega_0} \d v.$

  \item The linear map $W[\psi] : \ker (\L_{\Omega_0}^*) \to \ker \left( \mathcal{M}_{\Omega_0} - c_0 \mathrm{Id} \right) $ induces an isomorphism between the vector spaces $\ker (\L_{\Omega_0}^*)/\ker W$ and $\ker \left( \mathcal{M}_{\Omega_0} - c_0 \mathrm{Id} \right)$, where $\ker W$ is the set of the constant functions. As a consequence, we have $$\mathrm{dim} \left\{ \ker (\L_{\Omega_0}^*)/\ker W \right\} = \mathrm{dim}\ \ker \left( \mathcal{M}_{\Omega_0} - c_0 \mathrm{Id} \right) = d-1 .$$

\end{enumerate}

\end{proposition}

We can write the order of $\mathcal{O}(1)$ in \cref{eq:expansion-orders} as:
  \begin{align}
    \div_v \left\{ \nabla_v f_1 + f_1 \left[ (v- \Omega_0) + \nabla_v V(v) \right] - f_0 \cdot \widehat{\Omega}_1 \right\}
    = (\p_t + v \cdot \nabla_x) f_0,
  \end{align}
or equivalently,
  \begin{align}\label{eq:1st-order}
    \L_{\Omega_0} f_1 = (\p_t + v \cdot \nabla_x) f_0.
  \end{align}

% \dnf{A solvability lemma is needed here, in a way of functional analysis. }

\begin{lemma}[Duality solvability lemma] \label{lemm:solvability}
  The restricted operator $\L_{\Omega_0}$ onto $\ker^{\bot}(\L_{\Omega_0})$, i.e., $\L_{\Omega_0} \big|_{\ker^{\bot}(\L_{\Omega_0})}:\, \ker^{\bot}(\L_{\Omega_0}) \to \ker^{\bot}(\L_{\Omega_0}^*)$, is a Fredholm operator, and has a pseudo-inverse. Namely, \cref{eq:1st-order} can be solved if and only if the right-hand side belongs to $\ker^{\bot}(\L_{\Omega_0}^*)$.

\end{lemma}

\proof
It should be emphasized here that, since the linear operator $\L_{\Omega_0}$ is not a self-adjoint operator, $\L_{\Omega_0} \neq \L_{\Omega_0}^*$. As stated in \Cref{prop:operator-L*}, the collision invariant should be considered in the kernel of the dual operator $\L_{\Omega_0}^*$. Thus, we can transfer the solvability of \cref{eq:1st-order} to the following sense, for any $\phi $ belonging to the domain of $\L_{\Omega_0}^*$, it holds
  \begin{align}
    \skp{f_1}{\L_{\Omega_0}^* \phi}_{L^2_v} = \skp{\L_{\Omega_0} f_1}{\phi}_{L^2_v}.
  \end{align}
This indicates us to consider the orthogonal decomposition of the domain of $\L_{\Omega_0}^*$: $\textrm{Dom} ({\L_{\Omega_0}^*}) = \ker(\L_{\Omega_0}^*) + \ker^\bot(\L_{\Omega_0}^*)$. By virtue of proposition for linear bounded operator, the latter part can be rewritten as $\ker^\bot(\L_{\Omega_0}^*) = \textrm{Im} \L_{\Omega_0}$, since the image of $\L_{\Omega_0}$ is closed. Thus, we have already obtained, the operator $\L_{\Omega_0}$ restricted onto $\ker^{\bot}(\L_{\Omega_0})$, denoting by $\L_{\Omega_0} \big|_{\ker^{\bot}(\L_{\Omega_0})} $, is a bijection from $\ker^{\bot}(\L_{\Omega_0})$ to $\ker^{\bot}(\L_{\Omega_0}^*)$.

Notice for the operator $\L_{\Omega_0}$, the dimension of its kernel space and the codimension of its image space are finite, due to \cref{Ker-L-Omega0}, and the fact $ \textrm{Dom}({\L_{\Omega_0}}) / \textrm{Im} \L_{\Omega_0} = \textrm{Dom}({\L_{\Omega_0}^*}) / \ker^\bot(\L_{\Omega_0}^*) = \ker(\L_{\Omega_0}^*)$. So the operator $\L_{\Omega_0}$ is a Fredholm operator, and it has a pseudo-inverse. In turn, this provides the solvability of \cref{eq:1st-order}, that the kernel orthogonal part of $f_1$ can be solved if and only if the right-hand side belongs to $\textrm{Im} \L_{\Omega_0}$, namely, $\ker^{\bot}(\L_{\Omega_0}^*)$.

Now it remains to prove $\L_{\Omega_0}$ is a bounded operator. Perform the following weighted inner product:
  \begin{align}
    \skp{\mathcal{L}_{\Omega_0}^* \psi}{\phi}_{L^2_v(M)}
       & = \skp{ \div_v \left( M_{\Omega_0} \nabla_v \psi \right) }{\phi}_{L^2_v} + c_0^{-1} \skp{v_{\bot}}{\phi}_{L^2_v(M)} \cdot \agl{\nabla_v \psi}_{M} \\ \no
       & = - \skp{ \nabla_v \psi }{\nabla_v \phi}_{L^2_v(M)} + c_0^{-1} \skp{v_{\bot}}{\phi}_{L^2_v(M)} \cdot \agl{\nabla_v \psi}_{M} \\ \no
       & \le \nm{\nabla_v \psi}_{\M} \nm{\nabla_v \phi}_{\M} + c_0^{-1} \agl{\abs{v_{\bot}}^2 M_{\Omega_0}}^{\frac{1}{2}} \nm{\phi}_{\M} \cdot \nm{\nabla_v \psi}_{\M} \agl{M_{\Omega_0}}^{\frac{1}{2}}
       \\ \no
       & \le C \nm{\nabla_v \psi}_{L^2_v(M)} \nm{\phi}_{H^1_v(M)}.
  \end{align}
This yields the boundness of the operators $\L_{\Omega_0}^*$ in a weighted Sobolev space $H^1_v(M)$, and the boundness of $\L_{\Omega_0}$ as well, by noticing the fact $\skp{\L_{\Omega_0}^* \psi}{\phi}_{L^2_v(M)} =\skp{\L_{\Omega_0} (\phi M_{\Omega_0})}{\psi}_{L^2_v}$.
\endproof

Based on the solvability in \Cref{lemm:solvability}, we perform an inner product on \cref{eq:1st-order} with the functions belonging to $\ker(\L_{\Omega_0}^*)$. On one hand, this scheme provides the solvability of the first-order function $f_1$ (modulo $\ker(\L_{\Omega_0})$), and on the other hand, this scheme also enables us to deduce the balance laws for the macroscopic parameters $\rho_0(t,x)$ and $\Omega_0(t,x)$ arising from the leading order function $f_0$.

Indeed, this macroscopic limit system satisfied by $\rho_0(t,x)$ and $\Omega_0(t,x)$ has been derived, appears to be a self-organized hydrodynamic (SOH) system, see more details in \cite[Theorem 1.1, pp. 7886]{ABCD-19mbe}:
\begin{proposition}[\cite{ABCD-19mbe}] \label{prop:limit-0-order}

The limit self-organized hydrodynamic system, satisfied by the concentration $\rho_0$ and the orientation $\Omega_0$ reads as follows:
  \begin{align}%\label{sys:limit-0-order}
  \begin{cases}
    \p_t \rho + \div_x (c_0 \rho \Omega) = 0, \\[3pt]
    \p_t \Omega + c_1 (\Omega \cdot \nabla_x ) \Omega + \sigma (\mathrm{Id} - \Omega \otimes \Omega) \nabla_x \ln \rho = 0,
  \end{cases}
  \end{align}
where $c_0 = \int_{\RR^d} (v \cdot \Omega_0) M_{\Omega_0} (v) \d v $ as before, and $c_1= \tfrac{\tilde c_1}{\tilde c_0} $ with the following notations:
  \begin{align}\label{eq:notation-tc-01}
    \tilde{c}_0 & = \tfrac{1}{d-1} \int_{\RR^d} { \left[ |v|^2 - (v \cdot \Omega)^2 \right]^{\frac{1}{2}} \chi \left ( \tfrac{v \cdot \Omega}{|v|}, |v|\right ) M_{\Omega_0} (v) \d v},
  \\[5pt] \no
    \tilde{c}_1 & = \tfrac{1}{d-1} \int_{\RR^d} { (v \cdot \Omega_0) \left[ |v|^2 - (v \cdot \Omega)^2 \right]^{\frac{1}{2}} \chi \left ( \tfrac{v \cdot \Omega}{|v|}, |v|\right ) M_{\Omega_0} (v) \d v}.
  \end{align}

\end{proposition}

The above limit SOH system preserves the same hyperbolic structure as that obtained from the Vicsek dynamics (only with different coefficients), so the well-posedness issue is naturally obtained from the studies before. In \cite{DLMP-13maa}, the authors gave a first (local) well-posedness result for the SOH system in $H^m$ space ($m \ge 5/2$), with initial hypotheses on the positivity of initial density, and on a structural condition $\sin \theta^{\IN} >0$ to avoid possible singularity, arising from the spherical coordinates, i.e. $\Omega (\theta,\, \varphi)= ( \sin \theta \cos \varphi ,\, \sin \theta \sin \varphi ,\, \cos \theta )^\top$. As noticed in our previous work \cite{JLZ-20arma}, this singularity just comes from the choice of the coordinate, and hence can be avoided by adopting \emph{stereographic projection transform}. Let $\Omega (\phi,\,\psi)= \left( \tfrac{2 \phi}{W} ,\, \tfrac{2 \psi}{W} ,\, \tfrac{\phi^2 + \psi^2 - 1}{W} \right)^\top$, with $W = 1 + \phi^2 + \psi^2$. The well-posedness result holds in $H^m$ space for $m \ge 3$, without any other hypothesis except for the positivity of initial density. The same strategy is also applied here, so we omit to prove the well-posedness of the limit system.

Corresponding to the solvability stated above, we decompose the function $f_1$ into two parts: a kernel part $\ker(\L_{\Omega_0})$ and its orthogonal part $\ker^{\bot}(\L_{\Omega_0})$, denoted by $f_1^{\shortparallel} = \mathcal{P}_{\L} f_1$ % $f_1^{\text{\tiny $\sslash$}}$
and $f_1^{\perp} = (\mathrm{Id} - \mathcal{P}_{\L}) f_1 = \mathcal{P}_{\L}^{\perp} f_1$, respectively. Here the notations $\mathcal{P}_{\L}$ and $\mathcal{P}_{\L}^{\perp}$ denote the projection operators onto the two subspaces $\ker(\L_{\Omega_0})$ and $\ker^{\bot}(\L_{\Omega_0})$, respectively.
In other words, we have the following expression:
  \begin{align*}
    f_1 = f_1^{\shortparallel} + f_1^{\perp} \in \ker(\L_{\Omega_0}) + \ker^{\bot}(\L_{\Omega_0}).
  \end{align*}
The solvability \Cref{lemm:solvability} applied to \cref{eq:1st-order} ensures that, the kernel orthogonal part $f_1^{\perp}$ has been uniquely determined, in the sense of modulo $\ker(\L_{\Omega_0})$.

Recalling the structure of $\ker(\L_{\Omega_0})$ in \cref{Ker-L-Omega0}, we can write the following decomposition:
  \begin{align}\label{eq:decomps-f1}
    f_1 = f_1^{\shortparallel} + f_1^{\perp}
    = \rho_1 M_{\Omega_0} + c_0^{-1} (v_{\bot} \cdot {\bm u}_1) M_{\Omega_0} + f_1^{\perp},
  \end{align}
for some parameters $\rho_1(t,x) \in \RR$, and ${\bm u}_1(t,x) = P_{\Omega_0^\bot} j_1^{\shortparallel}(t,x) \in {(\RR \Omega_0)^\bot}$ which will be determined in the next round.

We mention that, the part $\rho_1(t,x) \in \ker(\L_{\mathrm{FP}})$ is orthogonal to the other two part which belongs to the kernel orthogonal part of Fokker-Planck operator $\L_{\mathrm{FP}}$, while the decomposition of $f_1^{\shortparallel}$ and $f_1^{\perp}$ is orthogonal in the sense of $\L_{\mathrm{FP}}^{-1}$. Indeed, note that \cref{Ker-L-Omega0} and \Cref{prop:operator-L*} ensure an isomorphism between $\ker(\L_{\Omega_0}^*)$ and $\ker(\L_{\Omega_0})$, moreover, we have the relation $\ker(\L_{\Omega_0}^*) = \L_{\mathrm{FP}}^{-1} \ker(\L_{\Omega_0})$, by virtue of \cite[Lemma 3.1, pp. 7893]{ABCD-19mbe} (or \cref{eq:kernel-L*} above).
In fact, we have $\chi(v) \tfrac{v_{\bot}}{\abs{v_{\bot}}} = \L_{\mathrm{FP}}^{-1} (v_{\bot})$, satisfying the weighted inner product $\skp{\L_{\mathrm{FP}} \left( \chi(v) \tfrac{v_{\bot}}{\abs{v_{\bot}}} \right) }{ \chi(v) \tfrac{v_{\bot}}{\abs{v_{\bot}}} }_M = \skp{v_{\bot}}{\L_{\mathrm{FP}}^{-1} (v_{\bot})}_M $.

\subsection{Derivation of the first order macroscopic system} % (fold)
\label{sub:derivation_of_the_first_order_macroscopic_system}

% subsection derivation_of_the_first_order_macroscopic_system (end)

We turn to consider the second-orders in the expansion formula \cref{eq:expansion-orders}, the balance of order $\mathcal{O}(\eps)$ can be written as:
  \begin{align}\label{eq:2nd-order-ori}
    \div_v \left\{ \nabla_v f_2 + f_2 \left[ (v- \Omega_0) + \nabla_v V(v) \right] - f_0 \cdot \widehat{\Omega}_2 \right\} - \div_v ( f_1 \cdot \widehat{\Omega}_1 )
    = (\p_t + v \cdot \nabla_x) f_1.
  \end{align}

Corresponding to the decomposition of $f_1 = f_1^{\shortparallel} + f_1^{\perp}$ in \cref{eq:decomps-f1}, we write
  \begin{align}\label{eq:Omega-1-re}
    \widehat{\Omega}_1 = (c_0\rho_0)^{-1} P_{\Omega_0^\bot} {j_1}
    & = (c_0\rho_0)^{-1} P_{\Omega_0^\bot} \int_{\RR^d} v (f_1^{\shortparallel} + f_1^{\perp}) \d v  \no \\
    & = (c_0\rho_0)^{-1} P_{\Omega_0^\bot} \left\{ j_1^{\shortparallel} + \int_{\RR^d} v f_1^{\perp} \d v \right\}  \no \\
    & = (c_0\rho_0)^{-1} \u_1 + (c_0\rho_0)^{-1} P_{\Omega_0^\bot} j_1^{\perp}  \no \\
    & = \widehat{\Omega}_1^{\shortparallel} + \widehat{\Omega}_1^{\perp} \,,
  \end{align}
where we have used the fact $\u_1 = P_{\Omega_0^\bot} j_1^{\shortparallel} \in {(\RR \Omega_0)^\bot}$ since \cref{eq:decomps-f1} ensures $j_1^{\shortparallel} = \agl{v f_1^{\shortparallel}}_{M} = c_0 \rho_1 \Omega_0 + {\bm u}_1$. Note that, $\widehat{\Omega}_1^{\perp} = (c_0\rho_0)^{-1} P_{\Omega_0^\bot} j_1^{\perp} = (c_0\rho_0)^{-1} \agl{v_{\bot} f_1^{\perp}}$ could be viewed as a known source term later, due to the solvability of $f_1^{\perp}$.

It follows from \cref{eq:expansion-Omega} that
    \begin{align}\label{eq:Omega-2-re}
      \widehat{\Omega}_2^i & = P_{\Omega_0^\bot}^{ik} \frac{j_2^k}{\abs{j_0}} - \frac{j_1 \cdot j_1}{\abs{j_0}^2} \frac{j_0^i}{\abs{j_0}} - 2 \frac{j_0 \cdot j_1}{\abs{j_0}^2} \frac{j_1^i}{\abs{j_0}} + 3 \frac{j_0^i \otimes j_0^k}{\abs{j_0}^2} \frac{j_0 \cdot j_1}{\abs{j_0}^2} \frac{j_1^k}{\abs{j_0}} \no \\[5pt]
    & = (c_0\rho_0)^{-1} P_{\Omega_0^\bot}^{ik} {j_2^k}
        - (c_0\rho_0)^{-2} (P_{\Omega_0^\bot}^{kl} : j_1^k \otimes j_1^l) \Omega_0^i
        - 2 (c_0\rho_0)^{-2} ({j_1^l} \Omega_0^l) \cdot P_{\Omega_0^\bot}^{ik} {j_1^k} \no \\[5pt]
    & = (c_0\rho_0)^{-1} P_{\Omega_0^\bot}^{ik} {j_2^k}
        - (c_0\rho_0)^{-2} \abs{P_{\Omega_0^\bot} j_1}^2 \Omega_0^i
        - 2 (c_0\rho_0)^{-2} (\Omega_0 j_1) P_{\Omega_0^\bot}^{ik} {j_1^k}
    \no \\[5pt]
    & = (c_0\rho_0)^{-1} P_{\Omega_0^\bot}^{ik} {j_2^k} + \widehat{\gamma}_2 (j_1) \,.
    \end{align}
Here we used $\widehat{\gamma}_2 (j_1)$ to denote some degree-2 homogeneous function depending on $j_1$. Combined with the above expression \cref{eq:Omega-1-re} of $\widehat{\Omega}_1$, this homogeneous function reduces to
    \begin{align}\label{eq:gamma-2-j1}
      -\widehat{\gamma}_2 (j_1)
    & = (c_0\rho_0)^{-2} \abs{P_{\Omega_0^\bot} j_1}^2 \Omega_0^i
          + 2 (c_0\rho_0)^{-2} (\Omega_0 j_1) P_{\Omega_0^\bot}^{ik} {j_1^k} \no \\[5pt]
    & = \abs{\widehat{\Omega}_1^{\shortparallel} + \widehat{\Omega}_1^{\perp}}^2 \Omega_0^i
        + 2 (c_0\rho_0)^{-1} [\Omega_0 (j_1^{\shortparallel} + j_1^{\perp})] (\widehat{\Omega}_1^{\shortparallel} + \widehat{\Omega}_1^{\perp})^i
    \no \\[5pt]
    & = (c_0\rho_0)^{-2} (|\u_1|^2 \Omega_0^i + 2 c_0 \rho_1 \u_1^i) \no \\[5pt] & \qquad
        + 2 \left[ (c_0\rho_0)^{-1} (\widehat{\Omega}_1^{\perp} \cdot \u_1) \Omega_0^i
             + (c_0\rho_0)^{-2} (\Omega_0 j_1^{\perp}) \u_1^i
             + \rho_0^{-1} \rho_1 (\widehat{\Omega}_1^{\perp})^i \right] \no \\[5pt]
      & \qquad
        + \left[ \abs{\widehat{\Omega}_1^{\perp}}^2 \Omega_0^i + 2 (c_0\rho_0)^{-1} (\Omega_0 j_1^{\perp}) (\widehat{\Omega}_1^{\perp})^i \right]
    \no \\[5pt]
    & = - \left( \gamma_{2,n} (j_1^\shortparallel) + \gamma_{2,l} (j_1^\shortparallel) + \gamma_{2,s} (j_1^\perp) \right) \,,
    \end{align}
where we denoted $\gamma_{2,n} (j_1^\shortparallel) = -(c_0\rho_0)^{-2} (|\u_1|^2 \Omega_0^i + 2 c_0 \rho_1 \u_1^i)$ by a quadratic \emph{nonlinear} function with respect to undetermined parameters $\rho_1$ and $\u_1$. The notation $\gamma_{2,l} (j_1^\shortparallel)$ indicates some \emph{linear} function with respect to undetermined parameter $\rho_1$ or $\u_1$ and known source part $j_1^\perp$, while $\gamma_{2,s} (j_1^\perp)$ stands for some known \emph{source} function with respect to determined function $j_1^\perp$. As a remark, the subscripts $n$, $l$ and $s$ indicates ``nonlinear'', ``linear'', and ``source'', respectively.

In order to analyze the solvability of \cref{eq:2nd-order-ori} involving the second-order expanded function $f_2$ (or say more strictly, $f_2^{\perp}$), we rewrite \cref{eq:2nd-order-ori} as
  \begin{align}\label{eq:2nd-order}
    \L_{\Omega_0} f_2 - \div_v ( f_1 \cdot \widehat{\Omega}_1 ) - \div_v ( f_0 \cdot \widehat{\gamma}_2 (j_1)) = (\p_t + v \cdot \nabla_x) f_1,
  \end{align}
more precisely,
  \begin{align}\label{eq:2nd-order-decomp}
    \L_{\Omega_0} f_2
    & = (\p_t + v \cdot \nabla_x) (f_1^{\shortparallel} + f_1^{\perp})
        + \div_v \left\{ (f_1^{\shortparallel} + f_1^{\perp}) \cdot
                      ( \widehat{\Omega}_1^{\shortparallel} + \widehat{\Omega}_1^{\perp} )
                \right\}
        \no \\[5pt] \no
    & \quad + \div_v \left\{ f_0 ( \gamma_{2,n} (j_1^\shortparallel) + \gamma_{2,l} (j_1^\shortparallel) + \gamma_{2,s} (j_1^\perp) ) \right\}
    \no \\[5pt] \no
    & = (\p_t + v \cdot \nabla_x) f_1^{\shortparallel}
        + \div_v ( f_1^{\shortparallel} \cdot \widehat{\Omega}_1^{\shortparallel}
                   + f_0 \cdot \gamma_{2,n} (j_1^\shortparallel) )
        + \underbrace{\div_v ( f_1^{\perp} \cdot \widehat{\Omega}_1^{\shortparallel}
                       + f_1^{\shortparallel} \cdot \widehat{\Omega}_1^{\perp}
                       + f_0 \cdot \gamma_{2,l} (j_1^{\shortparallel}) )}_{\text{linear terms}}
          \\[5pt] \no & \qquad
        + \underbrace{(\p_t + v \cdot \nabla_x) f_1^{\perp} + \div_v ( f_1^{\perp} \cdot \widehat{\Omega}_1^{\perp} + f_0 \cdot \gamma_{2,s} (j_1^{\perp}) )}_{\text{source terms}}
    \\[5pt]
    & = (\p_t + v \cdot \nabla_x) f_1^{\shortparallel}
        + \div_v ( f_1^{\shortparallel} \cdot \widehat{\Omega}_1^{\shortparallel} + f_0 \cdot \gamma_{2,n} (j_1^\shortparallel) )
        + U_{\text{linear}}(f_1^{\shortparallel}) + U_{\text{source}}(f_1^{\perp}) \,.
  \end{align}

By the solvability requirement, we can derive the governing system for macroscopic parameters $\rho_1(t,x)$ and ${\bm u}_1(t,x) \in {(\RR \Omega_0)^\bot}$ by performing the inner product with the collision invariant functions $\psi (v) \in \ker (\L_{\Omega_0}^*)$, more precisely, with $1$ and $\chi(v) \cdot \tfrac{v_{\bot}}{\abs{v_{\bot}}} = \chi \left( \tfrac{v \cdot \Omega_0}{|v|},\, |v| \right) \tfrac{P_{\Omega_0^\bot} v}{|P_{\Omega_0^\bot} v|}$, respectively.
\begin{proposition}\label{prop:limit-1st-order}
  The macroscopic system, satisfied by the first order concentration and momentum functions $(\rho_1,\,\u_1)$, reads as follows:
  \begin{align}\label{sys:limit-1st-order-rho}
    \p_t \rho_1 + \div_x (c_0 \Omega_0 \rho_1 + \u_1) = S_{\rho}(U_{\text{source}}(f_1^{\perp})) \,,
  \end{align}
and
  \begin{align}\label{sys:limit-1st-order-j}
    \p_t \u_1^i \, + &\, c_1 (\Omega_0 \cdot \nabla_x) \u_1^i
    \\[5pt]
      & + [ (\u_1 \cdot \nabla_x) \ln \rho_0 \cdot \Omega_0^i + c_2 (\u_1 \cdot \nabla_x) \Omega_0^i
      + (c_2 - c_1) P_{\Omega_0^\bot}^{il} \nabla_{x^l} \Omega_0^k \u_1^k
      + c_2 \div_{x} \Omega_0 \cdot \u_1^i]
    \no \\[5pt]
      & + c_0 \rho_1 P_{\Omega_0^\bot}^{il} \nabla_{x^l} \ln \left( \frac{\rho_1}{\rho_0} \right)
      + c_0 (\tilde{c}_0 \rho_0)^{-1} \rho_1 \u_1^i
      % - (1 + e_0) (c_0 \rho_0)^{-1} (|\u_1|^2 \Omega_0^i + c_0\rho_1 \u_1^i)
    \no \\[5pt] \no
    = & -(\tilde{c}_0 \rho_0)^{-1} \left[ c_{f_1^{\perp}}^{il} \u_1^l + 2(\Omega_0 j_1^{\perp}) \u_1^i \right]
      % \no \\[5pt] & \hspace*{7.5cm}
      - c_0 (\tilde{c}_0 \rho_0)^{-1} \rho_1 (P_{\Omega_0^\bot} j_1^{\perp})^i
      + S_{\u}(U_{\text{source}}(f_1^{\perp})) \,,
    \end{align}
with some linear terms on the right-hand side, corresponding to the terms in $U_{\text{linear}}(f_1^{\shortparallel})$, and some known source terms $S(U_{\text{source}}(f_1^{\perp}))$. They can be computed with a  dependence on the determined part $f_1^{\perp}$.

\end{proposition}

Note that, a quadratic nonlinearity arises in the last term of the left-hand side, indicated by $\rho_1 \u_1$. % and $|\u_1|^2$.

Here we have used the notation $c_2 = \tfrac{\tilde{c}_2}{\tilde{c}_0}$ with (see \cite{BGL-93cpam,DA-20m3as} for similar definitions):
  \begin{align}\label{eq:notation-4-tensor}
    \tilde{c}_2 = & \tfrac{1}{(d-1)(d+1)} \int_{\RR^d} { {\abs{v_{\bot}}}^3 \chi \left ( \tfrac{v \cdot \Omega}{|v|}, |v|\right ) M_{\Omega_0} (v) \d v} \,,
  \\ \no
    X_{\Omega_0^\bot}^{ikml} = & \int_{\RR^d} { \tfrac{v_{\bot}^i v_{\bot}^k v_{\bot}^m v_{\bot}^l}{\abs{v_{\bot}}} \chi \left ( \tfrac{v \cdot \Omega}{|v|}, |v|\right ) M_{\Omega_0} (v) \d v}
    = \tilde{c}_2 (P_{\Omega_0^\bot}^{im} P_{\Omega_0^\bot}^{kl} + P_{\Omega_0^\bot}^{il} P_{\Omega_0^\bot}^{km} + P_{\Omega_0^\bot}^{ik} P_{\Omega_0^\bot}^{ml}) \,.
  \end{align}
% and $e_0 = \tfrac{\tilde{e}_0}{\tilde{c}_0}$ with
%   \begin{align}\label{eq:notation-e0}
%     \tilde{e}_0 = \tfrac{1}{d-1} \int_{\RR^d} { \tfrac{V'(|v|)}{|v|} \left[ |v|^2 - (v \cdot \Omega)^2 \right]^{\frac{1}{2}} \chi \left ( \tfrac{v \cdot \Omega}{|v|}, |v|\right ) M_{\Omega_0} (v) \d v} \,.
%   \end{align}

%%----------------------------------------------------
% We write the detailed calculations as follows.

The detailed derivation for \Cref{prop:limit-1st-order} is postponed to \cref{sub:derivation_of_limit-1st-order} later.

%%--------------------------------------------------------------------
\subsection{Induction scheme for higher-order expansions} % (fold)
\label{sub:induction_scheme_for_higher_order_expansions}

% subsection induction_scheme_for_higher_order_expansions (end)

We are now already to write the following expression of $f_2$,
 \begin{align}\label{eq:decomps-f2}
    f_2 = f_2^{\shortparallel} + f_2^{\perp}
    = \rho_2 M_{\Omega_0} + c_0^{-1} (v_{\bot} \cdot {\bm u}_2) M_{\Omega_0} + f_2^{\perp},
  \end{align}
with the completely determined part $f_2^{\perp}$, and some undetermined parameters $\rho_2(t,x) \in \RR$, and ${\bm u}_2(t,x) \in {(\RR \Omega_0)^\bot}$. For that, we need to consider the next order equation, namely, the balance relation of order $\mathcal{O}(\eps^2)$ in \cref{eq:expansion-orders}:
  \begin{align}\label{eq:3rd-order-ori}
    \div_v \left\{ \nabla_v f_3 + f_3 \left[ (v- \Omega_0) + \nabla_v V(v) \right] - f_0 \cdot \widehat{\Omega}_3 \right\} - \div_v ( f_1 \cdot \widehat{\Omega}_2 + f_2 \cdot \widehat{\Omega}_1 )
    = (\p_t + v \cdot \nabla_x) f_2.
  \end{align}

Performing a similar but more complicated calculation as that in \cref{eq:Omega-2-re} yields that,
    \begin{align}\label{eq:Omega-3}
      \widehat{\Omega}_3
      = (c_0\rho_0)^{-1} P_{\Omega_0^\bot} j_3 + \widehat{\gamma}_3 (j_1,j_2) \,.
    \end{align}
The notation $\widehat{\gamma}_3 (j_1,j_2)$ denotes a degree-3 homogeneous function depending on $j_1$ and $j_2$, containing a cubic nonlinear dependence part $\gamma_{3,\text{nonlinear}} (j_2^\shortparallel)$ with coefficients depending on known part $j_2^\perp$ and (whole) $j_1$.

As a result, we can write
  \begin{align}\label{eq:3rd-order}
    \L_{\Omega_0} f_3
    & = (\p_t + v \cdot \nabla_x) f_2
      + \div_v \left[ f_0 \cdot \widehat{\gamma}_3 (j_1,j_2) + f_1 \cdot \widehat{\Omega}_2 + f_2 \cdot \widehat{\Omega}_1 \right]
    \no \\[5pt]
    & = (\p_t + v \cdot \nabla_x) f_2^{\shortparallel}
        + \div_v \left[ f_0 \cdot \gamma_{3,\text{nonlinear}} (j_2^\shortparallel) \right]
        + U_{\text{linear}}(f_2^{\shortparallel}) + U_{\text{source}}(f_2^{\perp},f_1) \,.
  \end{align}
Note that the whole information of first-order function $f_1$ should be viewed as completely known source in this round. Now the kernel-orthogonal part $f_3^\perp$ is able to be solved from the above equation \cref{eq:3rd-order}. At the same time, the solvability of \cref{eq:3rd-order} will provide the macroscopic system obeyed by $(\rho_2(t,x),\, {\bm u}_2(t,x) )$, which should satisfy a similar formulation as \cref{sys:limit-1st-order-rho}-\eqref{sys:limit-1st-order-j}.

Therefore, we can write $f_3 = f_3^{\shortparallel} + f_3^{\perp} = \rho_3 M_{\Omega_0} + c_0^{-1} (v_{\bot} \cdot {\bm u}_3) M_{\Omega_0} + f_3^{\perp}$, and obtain the next order balance relation in \cref{eq:expansion-orders}, i.e.,
  \begin{align}\label{eq:4th-order}
    \L_{\Omega_0} f_4
    & = (\p_t + v \cdot \nabla_x) f_3
      + \div_v \left[ f_0 \cdot \widehat{\gamma}_4 (j_1,j_2,j_3) + f_1 \cdot \widehat{\Omega}_3 + f_2 \cdot \widehat{\Omega}_2 + f_3 \cdot \widehat{\Omega}_1 \right]
    \no \\[5pt]
    & = (\p_t + v \cdot \nabla_x) f_3^{\shortparallel}
        + \div_v \left[ f_0 \cdot \gamma_{4,\text{nonlinear}} (j_3^\shortparallel) \right]
        + U_{\text{linear}}(f_3^{\shortparallel}) + U_{\text{source}}(f_3^{\perp},f_2,f_1) \,.
  \end{align}
The above argument will enable us to obtain the kernel-orthogonal part $f_4^\perp$, and to obtain the macroscopic system obeyed by $(\rho_3(t,x),\, {\bm u}_3(t,x) )$ for determining the kernel part $f_3^{\shortparallel}$.

We can continue the same solvability argument to obtain information for higher-order expansion. Concerning a general expansion ansatz,
  \begin{align}\label{eq:expansion-ansatz-general}
    f^\eps = f_0 + \eps f_1 + \eps^2 f_2 + \cdots + \eps^n f_n + \eps^{n+1} f_{n+1} + \cdots.
  \end{align}
We write the $n$-th order function, that
  \begin{align}
    f_n = f_n^{\shortparallel} + f_n^{\perp}
        = \rho_n M_{\Omega_0} + c_0^{-1} (v_{\bot} \cdot {\bm u}_n) M_{\Omega_0} + f_n^{\perp},
  \end{align}
in which the kernel-orthogonal part $f_n^{\perp}$ could be determined from the previous round of order balance, and hence it plays a role of known source. To determine the macroscopic parameters $(\rho_n(t,x),\, {\bm u}_n(t,x) )$ in kernel part, we consider the next order balance relation in \cref{eq:expansion-orders}, i.e.,
  \begin{align}\label{eq:n-order}
    \L_{\Omega_0} f_{n+1}
    & = (\p_t + v \cdot \nabla_x) f_n
      + \div_v \left[ f_0 \cdot \widehat{\gamma}_{n+1} (j_1,j_2,\cdots,j_n) + f_1 \cdot \widehat{\Omega}_n + f_2 \cdot \widehat{\Omega}_{n-1} + \cdots + f_n \cdot \widehat{\Omega}_1 \right]
    \no \\[5pt]
    & = (\p_t + v \cdot \nabla_x) f_n^{\shortparallel}
        + \div_v \left[ f_0 \cdot \gamma_{n+1,\text{nonlinear}} (j_n^\shortparallel) \right]
        + U_{\text{linear}}(f_n^{\shortparallel}) + U_{\text{source}}(f_n^{\perp},f_{n-1},\cdots,f_1) \,.
  \end{align}
Besides that the kernel-orthogonal part $f_{n+1}^\perp$ can be solved from \cref{eq:n-order}, the solvability of \cref{eq:3rd-order} also provide a way to obtain the macroscopic system for $(\rho_n(t,x),\, {\bm u}_n(t,x) )$, as expressed in following abstract formulation:
  \begin{align}\label{sys:limit-n-order}
  \begin{cases}
    \p_t \rho_n + \div_x (c_0 \Omega_0 \rho_n + \u_n) = S_{\rho}(U_{\text{source}}(f_n^{\perp})) \,, \\[3pt]
    \p_t \u_n + c_1 (\Omega_0 \cdot \nabla_x ) \u_n + \sigma c_0 P_{\Omega_0^\bot} \nabla_x \rho_n
      + \mathcal{N}_{n+1} (\rho_n,\u_n) = S_{\u}(U_{\text{linear}}(f_n^{\shortparallel}),\, U_{\text{source}}(f_n^{\perp},f_{n-1},\cdots,f_1)) \,.
  \end{cases}
  \end{align}
This enjoys a similar structure as \cref{sys:limit-1st-order-rho}-\eqref{sys:limit-1st-order-j}. Here the linear and source terms on the right-hand side could be computed, depending on the determined part $f_1^{\perp}$. We emphasize that, there is a nonlinear term $\mathcal{N}_{n+1} (\rho_n,\u_n)$ with nonlinarity of order $n+1$, associating to the degree-$(n+1)$ homogeneous function $\gamma_{n+1,\text{nonlinear}} (j_n^\shortparallel)$ in \cref{eq:n-order}. Indeed, this arises from the Taylor expansion formula \cref{eq:taylor-Omega} with nonlinear coefficients for mean orientation $\Omega_{f^\eps}$.

%%--------------------------------------------------------------------
\subsection{Truncation and the remainder equation} % (fold)
\label{sub:truncation_and_the_remainder_equation}

% subsection truncation_and_the_remainder_equation (end)

For sake of analyzing this expansion ansatz, we consider the following truncated expansion formula with a remainder:
% we need to truncate it, and to analyze the remainder system. More precisely, we consider the following expansion formula:
  \begin{align}\label{eq:expansion-truncate}
    f^\eps = f_0 + \eps f_1 + \eps^2 f_2 + \eps^2 f_\R^\eps \,.
  \end{align}
More precisely, we will justify that
  \begin{align*}%\label{eq:expansion-truncate}
    f^\eps = f_0 + \eps (f_1^{\shortparallel} + f_1^{\perp}) + \eps^2 f_2^{\perp} + \eps^2 f_\R^\eps \,,
  \end{align*}
where the whole information of $f_1$ is completely determined, while the kernel part $f_2^{\shortparallel}$ is assumed to be vanished such that $f_2 = f_2^{\perp}$. In particular, by a Taylor's formula with remainder of order 3, we can infer that
  \begin{align}\label{eq:Omega-ansatz}
    \Omega_{f^\eps}
    & = \Omega_0 + \eps (c_0\rho_0)^{-1} P_{\Omega_0^\bot} j_1
        + \eps^2 \left[ (c_0\rho_0)^{-1} P_{\Omega_0^\bot} j_2 + \widehat{\gamma}_2 (j_1) \right]
        + \eps^2 (c_0\rho_0)^{-1} P_{\Omega_0^\bot} j_\R^\eps
        + \eps^3 \widehat{\Omega}_{3,R}^\eta
     \no\\
    & = \Omega_0 + \eps \widehat{\Omega}_1 + \eps^2 \widehat{\Omega}_2 + \eps^2 \widehat{\Omega}_\R^\eps.
  \end{align}
As mentioned in \cref{eq:taylor-Omega-n} before, the remainder term $\widehat{\Omega}_\R^\eps$ appears just as a notation, defined through the difference $\eps^2 \widehat{\Omega}_\R^\eps = \Omega_{f^\eps} - ( \Omega_0 + \eps \widehat{\Omega}_1 + \eps^2 \widehat{\Omega}_2 ) = \eps^2 (c_0\rho_0)^{-1} P_{\Omega_0^\bot} j_\R^\eps + \eps^3 \widehat{\Omega}_{3,R}(\eta)$, where the term $\widehat{\Omega}_{3,R}^\eta$ stands for a Taylor's remainder term of order 3, with the parameter $\eta \in (0,\, \eps)$.

Recalling the order balance relations \cref{eq:0-order,eq:1st-order,eq:2nd-order}, we can get from \cref{eq:expansion-orders} the remainder equations, as follows:
  \begin{multline*}
    \div_v \left\{ \nabla_v f_\R^\eps + f_\R^\eps \left[ (v- \Omega_0) + \nabla_v V(v) \right] - f_0 \cdot \widehat{\Omega}_\R^\eps \right\}  \\
      - \div_v \left(\eps f_1 \cdot (\widehat{\Omega}_2 + \widehat{\Omega}_\R^\eps) \right)
      - \div_v \left\{ (f_2 + f_\R^\eps) \cdot \left[ \eps \widehat{\Omega}_1 + \eps^2 (\widehat{\Omega}_2 + \widehat{\Omega}_\R^\eps) \right] \right\} \\
    = \eps (\p_t + v \cdot \nabla_x) (f_2 + f_\R^\eps).
  \end{multline*}
Note the first term on the left-hand side is not the same as the linear operator.

Noticing $\eps \widehat{\Omega}_{3,R} = \widehat{\Omega}_\R^\eps - (c_0\rho_0)^{-1} P_{\Omega_0^\bot} j_\R^\eps$, we write
  \begin{multline}\label{eq:Rem-order}
    (\p_t + v \cdot \nabla_x) f_\R^\eps
    = \frac{1}{\eps} \L_{\Omega_0} f_\R^\eps - \div_v \left( f_0 \cdot \widehat{\Omega}_{3,R} \right)
      \\
      - \div_v \left[ (f_1 + \eps f_2) \widehat{\Omega}_\R^\eps \right]
    % \\\no
      - \div_v \left[ f_\R^\eps \cdot \tfrac{\left( \Omega_{f^\eps} - \Omega_0 \right)}{\eps} \right]
        - \S(f_1,f_2) \,,
  \end{multline}
where we have employed a notation, $\S(f_1, f_2)$, to denote a known source term with dependence on $f_1$ and $f_2$ (and certainly, $f_0$). It is defined as:
  \begin{align*}
    \S(f_1,f_2) = \div_v \left[ f_1 \widehat{\Omega}_2 + f_2 (\widehat{\Omega}_1 + \eps \widehat{\Omega}_2) \right] + (\p_t + v \cdot \nabla_x) f_2 \,.
  \end{align*}
We also have used the expansion formula \cref{eq:Omega-ansatz}: %, with a remainder of order 1:
  \begin{align*}
    \widehat{\Omega}_1 + \eps (\widehat{\Omega}_2 + \widehat{\Omega}_\R^\eps) = \tfrac{\Omega_{f^\eps} - \Omega_0 }{\eps}.
  \end{align*}
Actually, this is a remainder of order 1 in Taylor's formula. Note that the only nonlinear dependence effect on $f_\R^\eps$ arises from the product $\eps f_\R^\eps \widehat{\Omega}_\R^\eps$ with one more order of $\eps$, which helps us get a better estimate than that on the original kinetic Cucker-Smale system \cref{eq:CS}.

\subsection{Detailed derivation of \texorpdfstring{\Cref{prop:limit-1st-order}}{Proposition 2.5}} % (fold)
\label{sub:derivation_of_limit-1st-order}

% subsection derivation_of_limit-1st-order (end)

% \smallskip\noindent\emph{Proof of \Cref{prop:limit-1st-order}}.
We provide the detailed derivation of the first-order macroscopic sytem \cref{sys:limit-1st-order-rho,sys:limit-1st-order-j} as follows.

%%--------------------------------------------------------------------
\subsubsection{Derivation of \texorpdfstring{\cref{sys:limit-1st-order-rho}}{eq. (2.29)}} % (fold)
\label{ssub:derivation_of_sys:limit-1st-order-rho}

% subsubsection derivation_of_sys:limit-1st-order-rho (end)
% \smallskip\noindent\textbf{Derivation of \cref{sys:limit-1st-order-rho}}.

Perform the inner product on \cref{eq:2nd-order-decomp} with $1 \in \ker (\L_{\Omega_0}^*)$, we have
  \begin{align}
    \int_{\RR^d} (\p_t + v \cdot \nabla_x) (f_1^{\shortparallel} + f_1^{\perp}) \d v
    = \int_{\RR^d} f_2\, \L_{\Omega_0}^* (1) \d v + \int_{\RR^d} \div_v \left( f_1 \, \widehat{\Omega}_1 + f_0 \cdot \widehat{\gamma}_2 (j_1) \right) \d v
    = 0 \,.
  \end{align}

The expression $f_1^{\shortparallel} = \mathcal{P}_{\L} f_1 = \rho_1 M_{\Omega_0} + c_0^{-1} (v_{\bot} \cdot \u_1) M_{\Omega_0}$, and the simple relation $\Omega_0 \cdot \p_t \Omega_0 =0$ yield that
  \begin{align}\label{eq:f1-para-dt}
    \p_t f_1^{\shortparallel}
    & = \p_t \left[ \rho_1 M_{\Omega_0} + c_0^{-1} \left( v_{\bot} \cdot \u_1 \right) M_{\Omega_0} \right]
    \no\\
    & = \left[ \p_t \rho_1 M_{\Omega_0} + \rho_1 M_{\Omega_0} (v- \Omega_0) \cdot \p_t \Omega_0 \right] \no\\
    & \quad
      + \left[c_0^{-1} ( v_{\bot} \cdot \p_t \u_1 + \p_t P_{\Omega_0^\bot} v \cdot \u_1 ) M_{\Omega_0}
            + c_0^{-1} \left( v_{\bot} \cdot \u_1 \right) M_{\Omega_0} (v- \Omega_0) \cdot \p_t \Omega_0 \right]
    \no\\
    & = \p_t \rho_1 M_{\Omega_0} + \rho_1 M_{\Omega_0} (v_{\bot} \cdot \p_t \Omega_0) + c_0^{-1} (v_{\bot} \cdot \p_t \u_1) M_{\Omega_0} \no\\
    & \quad
      - c_0^{-1} (\p_t \Omega_0 \cdot \u_1) (v \cdot \Omega_0) M_{\Omega_0}
      + c_0^{-1} \left( v_{\bot} \cdot \u_1 \right) M_{\Omega_0} (v_{\bot} \cdot \p_t \Omega_0) \,,
  \end{align}
then we can infer, from the facts $\int_{\RR^d} M_{\Omega_0} \d v=1$, and $\int_{\RR^d} v_{\bot} M_{\Omega_0} \d v= 0$, that the second and third terms will vanish in integrating, and hence
  \begin{align}\label{calc:rho-1-dt}
    \int_{\RR^d} \p_t f_1^{\shortparallel} \d v = \p_t \rho_1,
  \end{align}
where we have used
  \begin{align*}
    \int_{\RR^d} - (\p_t \Omega_0 \cdot \u_1) (v \cdot \Omega_0) M_{\Omega_0} \d v = - c_0 (\p_t \Omega_0 \cdot \u_1),
  \end{align*}
and,
  \begin{align*}
    \int_{\RR^d} ( v_{\bot}^k \cdot \u_1^k ) (v_{\bot}^l \cdot \p_t \Omega_0^l) M_{\Omega_0} \d v
   = \u_1^k\, \p_t \Omega_0^l \int_{\RR^d} v_{\bot}^k v_{\bot}^l M_{\Omega_0} \d v
   = c_0 \,\delta_{kl}\, \u_1^k \p_t \Omega_0^l
   = c_0 (\p_t \Omega_0 \cdot \u_1) \,,
  \end{align*}
due to the definition $c_0 = \int_{\RR^d} (v \cdot \Omega_0) M_{\Omega_0} (v) \d v $, and
% $$\int_{\RR^d} v_{\bot}^k v_{\bot}^l M_{\Omega_0} \d v = \delta_{kl}\, \tfrac{1}{d-1} \int_{\RR^d} \left[ |v|^2-(v \cdot \Omega_0)^2 \right]  M_{\Omega_0} \d v = c_0\, \delta_{kl} \,.$$
$$\int_{\RR^d} v \otimes v_{\bot} M_{\Omega_0} \d v = \tfrac{1}{d-1} \int_{\RR^d} \left[ |v|^2-(v \cdot \Omega_0)^2 \right]  M_{\Omega_0} \d v P_{\Omega_0^\bot}= c_0 P_{\Omega_0^\bot}\,.$$

We can also compute that
  \begin{align}\label{calc:rho-1-divx}
    \int_{\RR^d} v \cdot \nabla_x f_1^{\shortparallel} \d v
    & = \nabla_x \cdot \int_{\RR^d} v \left[ \rho_1 M_{\Omega_0} + c_0^{-1} \left( v_{\bot} \cdot \u_1 \right) M_{\Omega_0} \right] \d v \no\\
    & = \nabla_x \cdot \left( \rho_1 \int_{\RR^d} v M_{\Omega_0} \d v + c_0^{-1} \u_1 \int_{\RR^d} v \otimes v_{\bot} M_{\Omega_0} \d v \right) \no\\
    & = \div_x \left( c_0 \Omega_0 \rho_1 + \u_1 \right) \,.
  \end{align}

Therefore, we obtain \cref{sys:limit-1st-order-rho} by combining the above calculations \cref{calc:rho-1-dt} and \cref{calc:rho-1-divx}, with some source terms depending on the known function $f_1^{\perp}$.

%%--------------------------------------------------------------------
\subsubsection{Derivation of \texorpdfstring{\cref{sys:limit-1st-order-j}}{eq. (2.30)}} % (fold)
\label{ssub:derivation_of_sys:limit-1st-order-j}

% subsubsection derivation_of_sys:limit-1st-order-j (end)
% \smallskip\noindent\textbf{Derivation of \cref{sys:limit-1st-order-j}}.

Next we turn to consider the derivation of \cref{sys:limit-1st-order-j}. Perform the inner product on \cref{eq:2nd-order-decomp} with $\chi(v) \cdot \tfrac{v_{\bot}}{\abs{v_{\bot}}} = \chi \left( \tfrac{v \cdot \Omega_0}{|v|},\, |v| \right) \tfrac{P_{\Omega_0^\bot} v}{\abs{P_{\Omega_0^\bot} v}} \in \ker (\L_{\Omega_0}^*)$, we have
  \begin{multline}
    \int_{\RR^d} (\p_t + v \cdot \nabla_x) (f_1^{\shortparallel} + f_1^{\perp}) \chi(v) \cdot \tfrac{v_{\bot}}{\abs{v_{\bot}}} \d v \\
    = \int_{\RR^d} f_2\, \L_{\Omega_0}^* (\chi(v) \tfrac{v_{\bot}}{\abs{v_{\bot}}}) \d v
      + \int_{\RR^d} \div_v \left( f_1 \, \widehat{\Omega}_1 + f_0 \cdot \widehat{\gamma}_2 (j_1) \right) \chi(v) \cdot \tfrac{v_{\bot}}{\abs{v_{\bot}}} \d v \,,
  \end{multline}
or, by the expression in the last line of \cref{eq:2nd-order-decomp},
  \begin{align}\label{eq:skproduct-f1-para}
    0 = \skp{f_2} { \L_{\Omega_0}^* ( \chi(v) \tfrac{v_{\bot}}{\abs{v_{\bot}}} )  }
    = & \skp{ (\p_t + v \cdot \nabla_x) f_1^{\shortparallel} }
          { \chi(v) \cdot \tfrac{v_{\bot}}{\abs{v_{\bot}}} } \\\no
      & + \skp{ \div_v ( f_1^{\shortparallel} \cdot \widehat{\Omega}_1^{\shortparallel}
                        + f_0 \cdot \gamma_{2,n} (j_1^\shortparallel))}
            { \chi(v) \cdot \tfrac{v_{\bot}}{\abs{v_{\bot}}} } \\\no
      & + \skp{ U_{\text{linear}}(f_1^{\shortparallel}) + U_{\text{source}}(f_1^{\perp}) }
            { \chi(v) \cdot \tfrac{v_{\bot}}{\abs{v_{\bot}}} } \,.
  \end{align}

We focus on the first three terms on the right-hand side of \cref{eq:skproduct-f1-para}. Due to the expression \cref{eq:f1-para-dt}, we write
  \begin{align}
    & \skp{\p_t f_1^{\shortparallel}}{ \chi(v) \tfrac{v_{\bot}}{\abs{v_{\bot}}} } \\ \no
    & = \int_{\RR^d} \p_t \rho_1 M_{\Omega_0} { \chi(v) \tfrac{v_{\bot}}{\abs{v_{\bot}}} } \d v
       + \int_{\RR^d} \rho_1 M_{\Omega_0} (v_{\bot} \cdot \p_t \Omega_0) { \chi(v) \tfrac{v_{\bot}}{\abs{v_{\bot}}} } \d v
       + \int_{\RR^d} c_0^{-1} (v_{\bot} \cdot \p_t \u_1) M_{\Omega_0} { \chi(v) \tfrac{v_{\bot}}{\abs{v_{\bot}}} } \d v \\ \no
    & \quad
      - \int_{\RR^d} c_0^{-1} (\p_t \Omega_0 \cdot \u_1) (v \cdot \Omega_0) M_{\Omega_0} { \chi(v) \tfrac{v_{\bot}}{\abs{v_{\bot}}} } \d v
      + \int_{\RR^d} c_0^{-1} \left( v_{\bot} \cdot \u_1 \right) M_{\Omega_0} (v_{\bot} \cdot \p_t \Omega_0) { \chi(v) \tfrac{v_{\bot}}{\abs{v_{\bot}}} } \d v \,.
  \end{align}
Notice the first, the fourth and the last terms vanish here because of the odd-order factors $v_{\bot} = P_{\Omega_0^\bot} v$ in integrand, and in particular, of that
  \begin{align*}
    \int_{\RR^d} v_{\bot}^k v_{\bot}^l \tfrac{v_{\bot}^i}{\abs{v_{\bot}}} \chi ( \tfrac{v \cdot \Omega_0}{|v|},\, |v| ) M_{\Omega_0} \d v =0, \quad \text{ for any } i,j,k \in {1,2,\cdots,d-1}.
  \end{align*}
Due to the notation \cref{eq:notation-tc-01} and the relation
  \begin{align*}
    \int_{\RR^d} P_{\Omega_0^\bot} v \otimes \tfrac{P_{\Omega_0^\bot} v}{|P_{\Omega_0^\bot} v|} \chi (\tfrac{v \cdot \Omega}{|v|},\, |v| ) M_{\Omega_0} \d v
    = P_{\Omega_0^\bot} \tfrac{1}{d-1} \int_{\RR^d} [|v|^2-(v \cdot \Omega_0)^2]^{\frac{1}{2}} \chi (\tfrac{v \cdot \Omega}{|v|},\, |v| ) M_{\Omega_0} \d v
    = \tilde{c}_0 P_{\Omega_0^\bot} \,,
  \end{align*}
the left two terms are computed as:
  \begin{align}
    & \int_{\RR^d} \rho_1 M_{\Omega_0} (v_{\bot}^k \cdot \p_t \Omega_0^k) { \chi(v) \tfrac{v_{\bot}^i}{\abs{v_{\bot}}} } \d v
    + \int_{\RR^d} c_0^{-1} (v_{\bot}^k \cdot \p_t \u_1^k) M_{\Omega_0} { \chi(v) \tfrac{v_{\bot}^i}{\abs{v_{\bot}}} } \d v \\ \no
    & = \rho_1 \p_t \Omega_0^k \int_{\RR^d} v_{\bot}^k \cdot \tfrac{v_{\bot}^i}{\abs{v_{\bot}}} \chi ( \tfrac{v \cdot \Omega_0}{|v|},\, |v| ) M_{\Omega_0} \d v
    + c_0^{-1} \p_t \u_1^k \int_{\RR^d} v_{\bot}^k \cdot \tfrac{v_{\bot}^i}{\abs{v_{\bot}}} \chi ( \tfrac{v \cdot \Omega_0}{|v|},\, |v| ) M_{\Omega_0} \d v \\ \no
    & = \tilde{c}_0 \rho_1 \p_t \Omega_0^i + \tilde{c}_0 c_0^{-1} P_{\Omega_0^\bot}^{ik} \p_t \u_1^k\,.
  \end{align}
Here we have used again $\p_t \Omega_0 \in {(\RR \Omega_0)^\bot}$. So we get
  \begin{align}\label{calc:limit-j1-dt}
    \skp{\p_t f_1^{\shortparallel}}{ \chi(v) \tfrac{v_{\bot}^i}{\abs{v_{\bot}}} }
    = \tilde{c}_0 c_0^{-1} P_{\Omega_0^\bot}^{ik} \p_t \u_1^k + \tilde{c}_0 \rho_1 \p_t \Omega_0^i \,.
  \end{align}

We next compute that
  \begin{align}\label{eq:f1-para-dx}
    v \cdot \nabla_x f_1^{\shortparallel}
    & = v \cdot \nabla_x \left[ \rho_1 M_{\Omega_0} + c_0^{-1} \left( v_{\bot} \cdot \u_1 \right) M_{\Omega_0} \right] \no\\
    & = \left[ v \cdot \nabla_x \rho_1 M_{\Omega_0} + \rho_1 M_{\Omega_0} v \otimes (v- \Omega_0) : \nabla_x \Omega_0 \right] \no\\
    & \quad
      + \left[c_0^{-1} ( v \otimes v_{\bot} : \nabla_x \u_1 + v \cdot \nabla_x P_{\Omega_0^\bot} v \cdot \u_1 ) M_{\Omega_0}
            + c_0^{-1} (v_{\bot} \cdot \u_1) M_{\Omega_0} v \otimes (v- \Omega_0) : \nabla_x \Omega_0 \right]
    \no\\
    & = v \cdot \nabla_x \rho_1 M_{\Omega_0} + \rho_1 M_{\Omega_0} (v \otimes v_{\bot}: \nabla_x \Omega_0)
        + c_0^{-1} (v \otimes v_{\bot} : \nabla_x \u_1) M_{\Omega_0} \no\\
    & \quad
      - c_0^{-1} (v \otimes \u_1: \nabla_x \Omega_0) (v \cdot \Omega_0) M_{\Omega_0}
      + c_0^{-1} (v_{\bot} \cdot \u_1) M_{\Omega_0} (v \otimes v_{\bot}: \nabla_x \Omega_0) \,,
  \end{align}
and correspondingly, we split $\skpt{ v \cdot \nabla_x f_1^{\shortparallel} }{ \chi(v) \tfrac{v_{\bot}}{\abs{v_{\bot}}} }$ into five terms:
  \begin{align}
    \skp{ v \cdot \nabla_x \rho_1 M_{\Omega_0} }{ \chi(v) \tfrac{v_{\bot}}{\abs{v_{\bot}}} }
    & = \nabla_x \rho_1 \int_{\RR^d} v \otimes \tfrac{v_{\bot}}{\abs{v_{\bot}}} \chi ( \tfrac{v \cdot \Omega_0}{|v|},\, |v| ) M_{\Omega_0} \d v
    = \tilde{c}_0 P_{\Omega_0^\bot} \nabla_x \rho_1 \,,
  \no\\
    \skp{ \rho_1 M_{\Omega_0} (v \otimes v_{\bot}: \nabla_x \Omega_0) }{ \chi(v) \tfrac{v_{\bot}}{\abs{v_{\bot}}} }
    & = \rho_1 \int_{\RR^d} (v \cdot \Omega_0) v_{\bot} \otimes \tfrac{v_{\bot}}{\abs{v_{\bot}}} \chi ( \tfrac{v \cdot \Omega_0}{|v|},\, |v| ) M_{\Omega_0} \d v (\Omega_0 \cdot \nabla_x) \Omega_0 \no\\
    & = \tilde{c}_1 \rho_1 (\Omega_0 \cdot \nabla_x) \Omega_0 \,,
  \no\\
    \skp{ (v \otimes v_{\bot} : \nabla_x \u_1) M_{\Omega_0} }{ \chi(v) \tfrac{v_{\bot}}{\abs{v_{\bot}}} }
    & = \int_{\RR^d} (v \cdot \Omega_0) v_{\bot} \otimes \tfrac{v_{\bot}}{\abs{v_{\bot}}} \chi ( \tfrac{v \cdot \Omega_0}{|v|},\, |v| ) M_{\Omega_0} \d v (\Omega_0 \cdot \nabla_x) \u_1 \no\\
    & = \tilde{c}_1 P_{\Omega_0^\bot} [(\Omega_0 \cdot \nabla_x) \u_1] \,,
  \no\\ \no
    \skp{ (v \otimes \u_1: \nabla_x \Omega_0) (v \cdot \Omega_0) M_{\Omega_0} }{ \chi(v) \tfrac{v_{\bot}}{\abs{v_{\bot}}} }
    & = \tilde{c}_1 P_{\Omega_0^\bot}^{il} \nabla_{x^l} \Omega_0^k \cdot \u_1^k \,,
  \end{align}
  % \begin{align}
  %   & \skp{ v \cdot \nabla_x \rho_1 M_{\Omega_0} }{ \chi(v) \tfrac{v_{\bot}}{\abs{v_{\bot}}} }
  %   = \nabla_x \rho_1 \int_{\RR^d} v \otimes \tfrac{v_{\bot}}{\abs{v_{\bot}}} \chi ( \tfrac{v \cdot \Omega_0}{|v|},\, |v| ) M_{\Omega_0} \d v
  %   = \tilde{c}_0 P_{\Omega_0^\bot} \nabla_x \rho_1 \,,
  % \no\\
  %   & \skp{ \rho_1 M_{\Omega_0} (v \otimes v_{\bot}: \nabla_x \Omega_0) }{ \chi(v) \tfrac{v_{\bot}}{\abs{v_{\bot}}} }
  %   = \rho_1 \int_{\RR^d} (v \cdot \Omega_0) v_{\bot} \otimes \tfrac{v_{\bot}}{\abs{v_{\bot}}} \chi ( \tfrac{v \cdot \Omega_0}{|v|},\, |v| ) M_{\Omega_0} \d v (\Omega_0 \cdot \nabla_x) \Omega_0
  %   = \tilde{c}_1 \rho_1 (\Omega_0 \cdot \nabla_x) \Omega_0 \,,
  % \no\\
  %   & \skp{ (v \otimes v_{\bot} : \nabla_x \u_1) M_{\Omega_0} }{ \chi(v) \tfrac{v_{\bot}}{\abs{v_{\bot}}} }
  %   = \int_{\RR^d} (v \cdot \Omega_0) v_{\bot} \otimes \tfrac{v_{\bot}}{\abs{v_{\bot}}} \chi ( \tfrac{v \cdot \Omega_0}{|v|},\, |v| ) M_{\Omega_0} \d v (\Omega_0 \cdot \nabla_x) \u_1
  %   = \tilde{c}_1 P_{\Omega_0^\bot} [(\Omega_0 \cdot \nabla_x) \u_1] \,,
  % \no\\ \no
  %   & \skp{ (v \otimes \u_1: \nabla_x \Omega_0) (v \cdot \Omega_0) M_{\Omega_0} }{ \chi(v) \tfrac{v_{\bot}}{\abs{v_{\bot}}} }
  %   = \tilde{c}_1 P_{\Omega_0^\bot}^{il} \nabla_{x^l} \Omega_0^k \cdot \u_1^k \,,
  % \end{align}
and
  \begin{align}
    \skp{ (v_{\bot} \cdot \u_1) M_{\Omega_0} (v \otimes v_{\bot}: \nabla_x \Omega_0) }{ \chi(v) \tfrac{v_{\bot}}{\abs{v_{\bot}}} } % \\\no
    & = X_{\Omega_0^\bot}^{ikml} \nabla_{x^l} \Omega_0^m \u_1^k  \\ \no
    & = \tilde{c}_2 (P_{\Omega_0^\bot}^{im} P_{\Omega_0^\bot}^{kl} + P_{\Omega_0^\bot}^{il} P_{\Omega_0^\bot}^{km} + P_{\Omega_0^\bot}^{ik} P_{\Omega_0^\bot}^{ml}) \nabla_{x^l} \Omega_0^m \u_1^k \\ \no
    & = \tilde{c}_2 ( \nabla_{x^l} \Omega_0^i \u_1^l + P_{\Omega_0^\bot}^{il} \nabla_{x^l} \Omega_0^k \u_1^k + \div_{x} \Omega_0 \cdot \u_1^i ) \,.
  \end{align}

Therefore, we are readily to write:
  \begin{align}\label{calc:limit-j1-dx}
    \skp{ v \cdot \nabla_x f_1^{\shortparallel} }{ \chi(v) \tfrac{v_{\bot}^i}{\abs{v_{\bot}}} } % \\\no
    % = & \tilde{c}_0 P_{\Omega_0^\bot}^{il} \nabla_{x^l} \rho_1 + \tilde{c}_1 \rho_1 (\Omega_0 \cdot \nabla_x) \Omega_0^i % \\ \no &
    %   + c_0^{-1} \tilde{c}_1 [P_{\Omega_0^\bot}^{ik} (\Omega_0 \cdot \nabla_x) \u_1^k - P_{\Omega_0^\bot}^{il} \nabla_{x^l} \Omega_0^k \cdot \u_1^k] \no \\
    % & + c_0^{-1} \tilde{c}_2 ( \nabla_{x^l} \Omega_0^i \u_1^l + P_{\Omega_0^\bot}^{il} \nabla_{x^l} \Omega_0^k \u_1^k + \div_{x} \Omega_0 \cdot \u_1^i )
    % \no \\
    = & \tilde{c}_1 \rho_1 (\Omega_0 \cdot \nabla_x) \Omega_0^i +\tilde{c}_0 P_{\Omega_0^\bot}^{il} \nabla_{x^l} \rho_1
      + c_0^{-1} \tilde{c}_1 P_{\Omega_0^\bot}^{ik} (\Omega_0 \cdot \nabla_x) \u_1^k
      \\\no
    & + c_0^{-1} (\tilde{c}_2 - \tilde{c}_1) P_{\Omega_0^\bot}^{il} \nabla_{x^l} \Omega_0^k \cdot \u_1^k
      + c_0^{-1} \tilde{c}_2 [ (\u_1 \cdot \nabla_x) \Omega_0^i + \div_{x} \Omega_0 \cdot \u_1^i ] \,.
  \end{align}

Concerning the third term $\skpt{ \div_v ( f_1^{\shortparallel} \cdot \widehat{\Omega}_1^{\shortparallel} )}{ \chi(v) \tfrac{v_{\bot}}{\abs{v_{\bot}}} }$ in \cref{eq:skproduct-f1-para}, we can write
  \begin{align}
    \div_v (f_1^{\shortparallel} \cdot \widehat{\Omega}_1^{\shortparallel})
    & = \div_{v^l} \left\{ c_0^{-1} \rho_0^{-1} \u_1^l [\rho_1 M_{\Omega_0} + c_0^{-1} (v_{\bot}^k \cdot \u_1^k) M_{\Omega_0}] \right\} \no \\
    & = c_0^{-1} \rho_0^{-1} \rho_1 \u_1^l \nabla_{v^l} M_{\Omega_0}
        + c_0^{-2} \rho_0^{-1} \u_1^k \u_1^l \nabla_{v^l} \left( v_{\bot}^k M_{\Omega_0} \right) \,.
  \end{align}
Noticing the fact $\nabla_{v^l} M_{\Omega_0} = \div_v (M_{\Omega_0} \nabla_v v^l ) = -M_{\Omega_0} \L_{\mathrm{FP}} (v^l) $ and the self-adjoint property of Fokker-Planck operator $\L_{\mathrm{FP}}$, we have
  \begin{align}
    \skp{\nabla_{v^l} M_{\Omega_0}}{\chi(v) \tfrac{v_{\bot}^i}{\abs{v_{\bot}}}}
    & = - \skp{ \L_{\mathrm{FP}} (v^l) }{ \L_{\mathrm{FP}}^{-1} (v_{\bot}^i) }_M \no\\
    & = - \skp{ v^l }{ v_{\bot}^i }_M
    = - \agl{v^l v_{\bot}^i M_{\Omega_0} }
    = -c_0 P_{\Omega_0^\bot}^{il} \,.
  \end{align}
We can also get, by noticing the odd-order factors $v_{\bot}$ in integrand,
  \begin{align}
    \skp{ \nabla_{v^l} \left( v_{\bot}^k M_{\Omega_0} \right) }{\chi(v) \tfrac{v_{\bot}^i}{\abs{v_{\bot}}}}
    = \delta_{kl} \agl{ \chi(v) \tfrac{v_{\bot}^i}{\abs{v_{\bot}}} M_{\Omega_0} }
      + \skp{ v_{\bot}^k \left( v^l + V'(|v|) \tfrac{v^{l}}{|v|} \right) M_{\Omega_0} }{\chi(v) \tfrac{v_{\bot}^i}{\abs{v_{\bot}}}}
    = 0.
  \end{align}
Then, it follows that,
  \begin{align}\label{calc:limit-j1-divv}
    \skp{ \div_v ( f_1^{\shortparallel} \cdot \widehat{\Omega}_1^{\shortparallel} )}{ \chi(v) \tfrac{v_{\bot}^i}{\abs{v_{\bot}}} }
    = -\rho_0^{-1} \rho_1 \u_1^l P_{\Omega_0^\bot}^{il}
    = -\rho_0^{-1} \rho_1 \u_1^i \,.
  \end{align}

As for the fourth term involving $\div_v ( f_0 \cdot \gamma_{2,n}(j_1^{\shortparallel}) )$ in \cref{eq:skproduct-f1-para}, a similar scheme as above yields
  \begin{align}\label{calc:limit-j1-divv2}
    \skp{ \div_v ( f_0 \cdot \gamma_{2,n} (j_1^\shortparallel))}{ \chi(v) \tfrac{v_{\bot}^i}{\abs{v_{\bot}}} }
    & = -\rho_0 (c_0\rho_0)^{-2} (|\u_1|^2 \Omega_0^l + 2 c_0 \rho_1 \u_1^l) \skp{ \nabla_{v^l} M_{\Omega_0}}{ \chi(v) \tfrac{v_{\bot}^i}{\abs{v_{\bot}}} }
  \no \\[5pt]
    & = c_0^{-2} \rho_0^{-1} ( |\u_1|^2 \Omega_0^l + 2 c_0 \rho_1 \u_1^l) (c_0 P_{\Omega_0^\bot}^{il})
  \no \\[5pt]
    & = 2 \rho_0^{-1} \rho_1 \u_1^i \,.
  \end{align}

On the other hand, the linear terms containing $U_{\text{linear}}(f_1^{\shortparallel})$ in \cref{eq:skproduct-f1-para} can also be calculated. We rewrite it as
  \begin{align}
    U_{\text{linear}}(f_1^{\shortparallel})
    & = \div_v \left( f_1^{\perp} \cdot \widehat{\Omega}_1^{\shortparallel}
              + f_1^{\shortparallel} \cdot \widehat{\Omega}_1^{\perp}
              + f_0 \cdot \gamma_{2,l} (j_1^{\shortparallel}) \right) \\\no
    & = \nabla_v f_1^{\perp} \widehat{\Omega}_1^{\shortparallel}
      + \nabla_v \left[ \rho_1 M_{\Omega_0} + c_0^{-1} \left( v_{\bot} \cdot \u_1 \right) M_{\Omega_0} \right] \widehat{\Omega}_1^{\perp}
      + \nabla_v (\rho_0 M_{\Omega_0}) \gamma_{2,l} (j_1^{\shortparallel}) \\\no
    & = \nabla_{v^l} f_1^{\perp} \left( \widehat{\Omega}_1^{\shortparallel} \right)^l
        + \nabla_{v^l} M_{\Omega_0} \left[ \rho_1 \widehat{\Omega}_1^{\perp} + \rho_0 \gamma_{2,l} (j_1^{\shortparallel}) \right]^l
        + \nabla_{v^l} \left( v_{\bot}^k M_{\Omega_0} \right) c_0^{-1} \u_1^k \left( \widehat{\Omega}_1^{\shortparallel} \right)^l \,.
  \end{align}
Again using the same argument as above, we can obtain from \cref{eq:Omega-1-re,eq:gamma-2-j1} that
  \begin{align}\label{calc:limit-j1-linearU}
    & \skp{ U_{\text{linear}}(f_1^{\shortparallel}) }{ \chi(v) \tfrac{v_{\bot}^i}{\abs{v_{\bot}}} } \\\no
    & = \skp{ \nabla_{v^l} f_1^{\perp} }{ \chi(v) \tfrac{v_{\bot}^i}{\abs{v_{\bot}}} } \left( \widehat{\Omega}_1^{\shortparallel} \right)^l
      + \skp{\nabla_{v^l} M_{\Omega_0}}{\chi(v) \tfrac{v_{\bot}^i}{\abs{v_{\bot}}}} \left[ \rho_1 \widehat{\Omega}_1^{\perp} + \rho_0 \gamma_{2,l} (j_1^{\shortparallel}) \right]^l \\\no & \quad
      + \skp{ \nabla_{v^l} \left( v_{\bot}^k M_{\Omega_0} \right) }{\chi(v) \tfrac{v_{\bot}^i}{\abs{v_{\bot}}}} c_0^{-1} \u_1^k \left( \widehat{\Omega}_1^{\shortparallel} \right)^l
    \\\no
    & = c_{f_1^{\perp}}^{il} (c_0\rho_0)^{-1} \u_1^l - c_0 \rho_1 P_{\Omega_0^\bot}^{il} \widehat{\Omega}_1^{\perp} \\\no & \quad
        + 2 c_0 \rho_0 P_{\Omega_0^\bot}^{il} \left[ (c_0\rho_0)^{-1} (\widehat{\Omega}_1^{\perp} \cdot \u_1) \Omega_0^l
             + (c_0\rho_0)^{-2} (\Omega_0 j_1^{\perp}) \u_1^l
             + \rho_0^{-1} \rho_1 (\widehat{\Omega}_1^{\perp})^l \right]
    \\\no
    & = c_{f_1^{\perp}}^{il} (c_0\rho_0)^{-1} \u_1^l + \rho_0^{-1} \rho_1 (P_{\Omega_0^\bot} j_1^{\perp})^i
        + 2 (c_0\rho_0)^{-1} (\Omega_0 j_1^{\perp}) \u_1^i \,.
  \end{align}
Here the matrix coefficient $c_{f_1^{\perp}}(t,x)$ can be viewed as a known function depending on known $f_1^{\perp}$.

Combining the above calculations \cref{calc:limit-j1-dt,calc:limit-j1-dx,calc:limit-j1-divv,calc:limit-j1-divv2,calc:limit-j1-linearU} together will enable us to rewrite \cref{eq:skproduct-f1-para} as
  \begin{align}\label{eq:limit-j1-TBD}
    & c_0^{-1} P_{\Omega_0^\bot}^{ik} \left[ \tilde{c}_0 \p_t \u_1^k + \tilde{c}_1 (\Omega_0 \cdot \nabla_x) \u_1^k \right]
    + \rho_1 \left\{ \left[ \tilde{c}_0 \p_t \Omega_0^i + \tilde{c}_1 (\Omega_0 \cdot \nabla_x) \Omega_0^i \right] + \tilde{c}_0 P_{\Omega_0^\bot}^{il} \tfrac{\nabla_{x^l} \rho_1}{\rho_1} \right\}
  \no \\[5pt]
    & \qquad + c_0^{-1} (\tilde{c}_2 - \tilde{c}_1) P_{\Omega_0^\bot}^{il} \nabla_{x^l} \Omega_0^k \cdot \u_1^k
      + c_0^{-1} \tilde{c}_2 [ (\u_1 \cdot \nabla_x) \Omega_0^i + \div_{x} \Omega_0 \cdot \u_1^i ]
      + \rho_0^{-1} \rho_1 \u_i^i
    % \no \\[5pt] & \qquad
      % - c_0^{-2} \rho_0^{-1} (|\u_1|^2 \Omega_0^i + c_0 \rho_1 \u_1^i) (\tilde{c}_0 + \tilde{e}_0)
  \no \\[5pt]
    & = -(c_0\rho_0)^{-1} \left[ c_{f_1^{\perp}}^{il} \u_1^l + 2(\Omega_0 j_1^{\perp}) \u_1^i \right] - \rho_0^{-1} \rho_1 (P_{\Omega_0^\bot} j_1^{\perp})^i
        - \skp{ U_{\text{source}}(f_1^{\perp}) }{ \chi(v) \tfrac{v_{\bot}}{\abs{v_{\bot}}} } \,.
  \end{align}

We can get from the limit system \cref{sys:limit-0-order} of leading order functions $(\rho_0, \Omega_0)$ that,
  \begin{align}
    \left[ \tilde{c}_0 \p_t \Omega_0^i + \tilde{c}_1 (\Omega_0 \cdot \nabla_x) \Omega_0^i \right] + \tilde{c}_0 P_{\Omega_0^\bot}^{il} \frac{\nabla_{x^l} \rho_1}{\rho_1}
    = \tilde{c}_0 P_{\Omega_0^\bot}^{il} \tfrac{\nabla_{x^l} \rho_1}{\rho_1} - \sigma \tilde{c}_0 P_{\Omega_0^\bot}^{il} \tfrac{\nabla_{x^l} \rho_0}{\rho_0}
    = \tilde{c}_0 P_{\Omega_0^\bot}^{il} \nabla_{x^l} \ln \tfrac{\rho_1}{\rho_0} \,.
  \end{align}

Notice the simple fact $\u_1 \in {(\RR \Omega_0)^\bot}$ yields $\Omega_0^k \cdot \p_t \u_1^k = - \u_1^k \cdot \p_t \Omega_0^k$, then we can infer from \cref{sys:limit-0-order} that
  \begin{align}
    \Omega_0^i \u_1^k \left[ \tilde{c}_0 \p_t \Omega_0^k + \tilde{c}_1 (\Omega_0 \cdot \nabla_x) \Omega_0^k \right]
    = - \Omega_0^i \u_1^k \left[ \tilde{c}_0 P_{\Omega_0^\bot}^{kl} \nabla_{x^l} \ln \rho_0 \right]
    = - \tilde{c}_0 (\u_1^l \cdot \nabla_{x^l}) \ln \rho_0 \cdot \Omega_0^i \,.
  \end{align}
So we get
  \begin{align}
    P_{\Omega_0^\bot}^{ik} \left[ \tilde{c}_0 \p_t \u_1^k + \tilde{c}_1 (\Omega_0 \cdot \nabla_x) \u_1^k \right]
    & = (\delta_{ik} - \Omega_0^i \Omega_0^k) \left[ \tilde{c}_0 \p_t \u_1^k + \tilde{c}_1 (\Omega_0 \cdot \nabla_x) \u_1^k \right] \no \\
    & = \left[ \tilde{c}_0 \p_t \u_1^i + \tilde{c}_1 (\Omega_0 \cdot \nabla_x) \u_1^i \right]
      + \tilde{c}_0 (\u_1 \cdot \nabla_x) \ln \rho_0 \cdot \Omega_0^i \,.
  \end{align}

Therefore, \cref{eq:limit-j1-TBD} can be reduced to:
  \begin{align} % \label{eq:limit-j1}
    & c_0^{-1} \left[ \tilde{c}_0 \p_t \u_1^i + \tilde{c}_1 (\Omega_0 \cdot \nabla_x) \u_1^i \right]
    + c_0^{-1} \tilde{c}_0 (\u_1 \cdot \nabla_x) \ln \rho_0 \cdot \Omega_0^i
    + \tilde{c}_0 \rho_1 P_{\Omega_0^\bot}^{il} \nabla_{x^l} \ln \tfrac{\rho_1}{\rho_0} \no \\
    & \qquad
      + c_0^{-1} (\tilde{c}_2 - \tilde{c}_1) P_{\Omega_0^\bot}^{il} \nabla_{x^l} \Omega_0^k \cdot \u_1^k
      + c_0^{-1} \tilde{c}_2 [ (\u_1 \cdot \nabla_x) \Omega_0^i + \div_{x} \Omega_0 \cdot \u_1^i ]
      + \rho_0^{-1} \rho_1 \u_i^i
      % \no \\ & \qquad
      % - c_0^{-2} \rho_0^{-1} (\tilde{c}_0 + \tilde{e}_0) (|\u_1|^2 \Omega_0^i + c_0\rho_1 \u_1^i)
  \no \\
    & = c_0^{-1} \left[ \tilde{c}_0 \p_t \u_1^i + \tilde{c}_1 (\Omega_0 \cdot \nabla_x) \u_1^i \right] \no \\
      & \qquad
        + c_0^{-1} \left[ \tilde{c}_0 (\u_1 \cdot \nabla_x) \ln \rho_0 \cdot \Omega_0^i + \tilde{c}_2 (\u_1 \cdot \nabla_x) \Omega_0^i \right]
        + c_0^{-1} (\tilde{c}_2 - \tilde{c}_1) P_{\Omega_0^\bot}^{il} \nabla_{x^l} \Omega_0^k \u_1^k
        + c_0^{-1} \tilde{c}_2 \div_{x} \Omega_0 \cdot \u_1^i \no \\
      & \qquad
        + \tilde{c}_0 \rho_1 P_{\Omega_0^\bot}^{il} \nabla_{x^l} \ln \tfrac{\rho_1}{\rho_0}
        + \rho_0^{-1} \rho_1 \u_i^i
        % - c_0^{-2} \rho_0^{-1} (\tilde{c}_0 + \tilde{e}_0) (|\u_1|^2 \Omega_0^i + c_0\rho_1 \u_1^i)
  \no \\
    & = -(c_0\rho_0)^{-1} \left[ c_{f_1^{\perp}}^{il} \u_1^l + 2(\Omega_0 j_1^{\perp}) \u_1^i \right] - \rho_0^{-1} \rho_1 (P_{\Omega_0^\bot} j_1^{\perp})^i
        - \skp{ U_{\text{source}}(f_1^{\perp}) }{ \chi(v) \tfrac{v_{\bot}}{\abs{v_{\bot}}} } \,,
  \end{align}
which leads us to the desired equation \cref{sys:limit-1st-order-j} by multiplying the constant $c_0 \tilde{c}_0^{-1}$.
\qed

%% ------------------------------------------------------------------------ %%
\section{Proof of the Main Result} % (fold)
\label{sec:proof_of_the_main_result}

% section proof_of_the_main_result (end)

% \section{A-priori Estimate} % (fold)
% \label{sec:a_priori_estimate}

% % section a_priori_estimate (end)

\subsection{Two types of coercivity estimates} % (fold)
\label{sub:two_types_of_coercivity_estimate}

% subsection two_types_of_coercivity_estimate (end)
We provide here two types of coercivity estimates for the linear operator $$\L_{\Omega_0} f = \div_v \left\{ M_{\Omega_0} \nabla_v \left( \frac{f_1}{M_{\Omega_0}} - c_0^{-1} P_{\Omega_0^\bot} v \cdot P_{\Omega_0^\bot} {j_1} \right) \right\},$$ one of which is without any requirements on coefficients $\alpha,\beta$ in the confining potential, and the other one is with such a structure condition. The latter coercivity estimate, referred as type II, is used to deal with the estimates on the remainder equations in \Cref{sub:uniform-esm}, while the former one of type I is for the estimates on the kernel orthogonal parts of lower-orders functions, see \Cref{sub:estimates_on_kernel_orthogonal_parts}.

% Now we go back to \cref{eq:1st-order}. We focus on the kernel orthogonal part $f_1^\perp$, denoted here by only $f_1$ for simplicity.
%   \begin{align}
%     \L_{\Omega_0} f_1 = (\p_t + v \cdot \nabla_x) f_0,
%   \end{align}
% or, by \cref{eq:linearized-op},
%   \begin{align}
%     \div_v \left\{ M_{\Omega_0} \nabla_v \left( \frac{f_1}{M_{\Omega_0}} - c_0^{-1} P_{\Omega_0^\bot} v \cdot P_{\Omega_0^\bot} {j_1} \right) \right\}
%     = (\p_t + v \cdot \nabla_x) f_0.
%   \end{align}
% By the short-hand notations $M_0 = M_{\Omega_0}$, $v_\bot =P_{\Omega_0^\bot} v$, we write
%   \begin{align}
%     \div_v \left\{ M_0 \nabla_v \left( \frac{f_1 - c_0^{-1} (v_\bot j_1)M_0}{M_0} \right) \right\}
%     = (\p_t + v \cdot \nabla_x) f_0.
%   \end{align}

% \begin{remark}
% This above formulation indicates we should consider the whole $f_1 - c_0^{-1} (v_\bot j_1)M_0$ as a test function. (A single $f_1$ was also tried, and failed to close estimate).

% Another point is, to look for a lower bound (or say a coercivity estimate) is very important. This needs a little tricky observation.
% \end{remark}

%% ------------------------------------------------------------------------
% \subsubsection{Coercivity estimate of Type I} % (fold)
% \label{ssub:coercivity_estimate_of_type_i}

% % subsubsection coercivity_estimate_of_type_i (end)

\smallskip\noindent\underline{\textbf{Coercivity estimate of Type I}}.
We focus on the kernel orthogonal part $f_1^\perp$, denoted here by only $f_1$ for simplicity. By the short-hand notations $M_0 = M_{\Omega_0}$, $v_\bot =P_{\Omega_0^\bot} v$, we write
  \begin{align}\label{eq:L-form-var}
    \L_{\Omega_0} f_1 = \div_v \left\{ M_0 \nabla_v \left( \frac{f_1 - c_0^{-1} (v_\bot j_1)M_0}{M_0} \right) \right\}.
    % = (\p_t + v \cdot \nabla_x) f_0.
  \end{align}

\begin{proposition}[Coercivity estimate, I] \label{prop:coercivity-I}
  We have the coercivity estimate for the linear operator $\L_{\Omega_0}$, for $f_1 = f_1^\perp \in \ker^{\bot}(\L_{\Omega_0})$, that
  \begin{align} \label{eq:coercivity-I}
    \skp{-\L_{\Omega_0} f_1}{\frac{f_1 - c_0^{-1} (v_\bot j_1)M_0}{M_0}}
    \ge \lambda_0 \agl{\abs{ \frac{f_1 - c_0^{-1} (v_\bot j_1)M_0}{M_0} }^2 }_M
    \ge \lambda_0 \eta_0 \agl{\frac{|f_1|^2}{M_0}} \,,
  \end{align}
with $\lambda_0$ being the Poincar\'e constant in the Fokker-Planck operator, and $\eta_0$ some small fixed parameter.

\end{proposition}
\proof
This above formulation \cref{eq:L-form-var} indicates we should consider the whole $f_1 - c_0^{-1} (v_\bot j_1)M_0$ as a test function. Observe that, by the Poincar\'e inequality for the Fokker-Planck operator,
  \begin{multline}\label{eq:low-bdd-f-diff}
    \skp{-\L_{\Omega_0} f_1}{\frac{f_1 - c_0^{-1} (v_\bot j_1)M_0}{M_0}} = \agl{M_0 \abs{\nabla_v \left( \frac{f_1 - c_0^{-1} (v_\bot j_1)M_0}{M_0} \right)}^2 } \\
    \ge \lambda_0 \agl{\abs{ \frac{f_1 - c_0^{-1} (v_\bot j_1)M_0}{M_0} }^2 }_M
    = \lambda_0 \agl{ \abs{ \frac{f_1}{M_0} - c_0^{-1} (v_\bot j_1) }^2 }_M ,
  \end{multline}
since we have a zero mean value $\agl{ {\frac{f_1 - c_0^{-1} (v_\bot j_1)M_0}{M_0}} }_M =0$, by $f_1^\bot - c_0^{-1} (v_\bot j_1)M_0 \in \ker^{\bot}(\L_{\mathrm{FP}})$.

We could infer that
  \begin{align}\label{eq:Holder-equlity}
    \abs{(P_{\Omega_0^\bot} {j_1})^k}^2
    = \agl{v_\bot^k \cdot f_1}^2
    \le \agl{|v_\bot^k|^2 M_0} \agl{\frac{|f_1|^2}{M_0}} = c_0 \agl{\frac{|f_1|^2}{M_0}}.
  \end{align}
However, this is not enough for our aim of obtaining a lower bound (only non-negative). We hope to look for in the above inequality an upper bound strictly smaller than 1.

In fact, we notice the equality in \cref{eq:Holder-equlity} holds if and only if $f_1 = c_k v_\bot^k M_0$ for some constant $c_k$, with respect to a fixed $k \in \{ 1,2,\cdots, d \} $. Since a generic $f_1$ can be expressed as $f_1 = g_1 + h_1 = \sum_l c_l v_\bot^l M_0 + h_1$ with $\skpt{g_1}{h_1} =0$, then the above H\"older inequality holds with a strict formulation, i.e.,
  \begin{align}\label{eq:Holder-inequlity}
    \agl{v_\bot^k \cdot f_1}^2
    < \agl{|v_\bot^k|^2 M_0} \agl{\frac{|f_1|^2}{M_0}}.
  \end{align}

% By direct calculations, we have
%   \begin{align}
%     \agl{ \abs{ \frac{f_1}{M_0} - c_0^{-1} (v_\bot j_1) }^2 }_M
%     & = \agl{ \abs{ \frac{f_1}{M_0}}^2 + c_0^{-2} (v_\bot^k j_1^k) (v_\bot^l j_1^l) - 2 c_0^{-1} (v_\bot^k j_1^k) \frac{f_1}{M_0} }_M \\\no
%     & = \agl{\frac{|f_1|^2}{M_0}} + c_0^{-2} \agl{v_\bot^k v_\bot^l M_0} j_1^k j_1^l - 2 c_0^{-1} \agl{v_\bot^k f_1} j_1^k \\\no
%     & = \agl{\frac{|f_1|^2}{M_0}} - c_0^{-1} \abs{(P_{\Omega_0^\bot} {j_{f_1}})}^2
%     \\\no
%     & = \agl{\frac{|h_1|^2}{M_0}} - c_0^{-1} \abs{(P_{\Omega_0^\bot} {j_{h_1}})}^2 \,,
%     % \\\no
%     % & \ge \eta_0 \agl{\frac{|f_1|^2}{M_0}} .
%   \end{align}
% where we have used the expression $f_1 = g_1 + h_1$, and

  % \begin{align}
  %   j_{g_1}^k = \sum_l c_l \int v^k v_\bot^l M_0 \d v
  %   = c_0 \sum_l c_l P_{\Omega_0^\bot}^{kl} \,,
  % \end{align}

Notice the elementary fact $j_{g_1}^k = \sum_l c_l \int v^k v_\bot^l M_0 \d v = c_0 \sum_l c_l P_{\Omega_0^\bot}^{kl}$, then it follows $g_1 - c_0^{-1} (v_\bot j_{g_1}) M_0=0 $ (in fact, this is due to $g_1 \in \ker(\L_{\Omega_0}) $). By direct calculations, we have
  \begin{align}\label{eq:difference-f-h}
    \skp{ \tfrac{f_1}{M_0} - c_0^{-1} (v_\bot j_{f_1})}{f_1}
    & = \skp{ \tfrac{g_1}{M_0} - c_0^{-1} (v_\bot j_{g_1})}{g_1 + h_1}
        + \skp{ \tfrac{h_1}{M_0} - c_0^{-1} (v_\bot j_{h_1})}{g_1 + h_1}
        \no\\
    & = \skp{ \tfrac{h_1}{M_0} - c_0^{-1} (v_\bot j_{h_1})}{h_1}
        \no\\
    & = \agl{\frac{|h_1|^2}{M_0}} - c_0^{-1} \abs{(P_{\Omega_0^\bot} {j_{h_1}})}^2 \,,
  \end{align}
where we have used $\skpt{ \tfrac{h_1}{M_0} - c_0^{-1} (v_\bot j_{h_1})}{g_1} = \skpt{h_1}{ \tfrac{g_1}{M_0} - c_0^{-1} (v_\bot j_{g_1})} = 0$. The orthogonal relation $\skpt{g_1}{h_1} =0$, or equally say, $\skpt{v_\bot^k M_0}{h_1} =0$, enables us to get, for a small enough constant $\eta < 1$, that
  \begin{align}\label{eq:strict-difference}
    \agl{v_\bot^k h_1}^2 = \agl{v_\bot^k \left( 1- \eta M_0 \right) \cdot h_1}^2
    & = \agl{v_\bot^k \left( 1- \eta M_0 \right) \sqrt{M_0} \cdot \tfrac{h_1}{\sqrt{M_0}} }^2 \\\no
    & \le \agl{|v_\bot^k|^2 M_0 \left( 1- \eta M_0 \right)^2} \agl{\frac{|h_1|^2}{M_0}} \\\no
    & \le \agl{|v_\bot|^2 M_0 \left( 1- 2 \eta M_0 + \eta^2 M_0^2 \right)} \agl{\frac{|h_1|^2}{M_0}} \\\no
    & \le c_0 (1- \eta_0) \agl{\frac{|h_1|^2}{M_0}} \,,
  \end{align}
where the basic fact $2 \eta M_0 - \eta^2 M_0^2 \ge \eta M_0 (2- \eta M_0) > \eta M_0$ have been used in the last inequality, and $\eta_1 = \eta c_0^{-1} \agl{|v_\bot^k|^2 M_0^2}$. Combined with \cref{eq:strict-difference}, this inequality yields
  \begin{align}\label{eq:low-bdd-h}
    \skp{ \tfrac{f_1}{M_0} - c_0^{-1} (v_\bot j_{f_1})}{f_1}
    = \agl{\frac{|h_1|^2}{M_0}} - c_0^{-1} \abs{(P_{\Omega_0^\bot} {j_{h_1}})}^2
    \ge \eta_1 \agl{\frac{|h_1|^2}{M_0}} \,.
  \end{align}

On the other hand, by recalling from \cref{sub:derivation_of_macroscopic_system} the fact $f_1 \in \ker^\bot (\L_{\Omega_0})$ is orthogonal to the $\ker (\L_{\Omega_0})$ in the action of $\L_{\mathrm{FP}}^{-1}$, namely, $\skpt{f_1}{\L_{\mathrm{FP}}^{-1} (v_\bot^k)} = \skpt{f_1}{ \tfrac{\chi(v)}{\abs{v_{\bot}}} v_{\bot}^k} = 0$, which implies $\skpt{\tfrac{f_1}{M_0}}{\L_{\mathrm{FP}}^{-1} (\tfrac{g_1}{M_0})}_M = 0$. Combining the decomposition of $f_1 = g_1 + h_1$ and the H\"older inequality ensures that:
  \begin{align}
    \skp{\tfrac{g_1}{M_0}}{\L_{\mathrm{FP}}^{-1} (\tfrac{g_1}{M_0})}_M = - \skp{\tfrac{h_1}{M_0}}{\L_{\mathrm{FP}}^{-1} (\tfrac{g_1}{M_0})}_M
    % & = - \skp{\tfrac{h_1}{M_0}}{\sum_l c_l \tfrac{\chi(v)}{\abs{v_{\bot}}} v_{\bot}^l}_M
    % \\\no &
    \le \agl{\tfrac{|h_1|^2}{M_0}}^\frac{1}{2} \agl{ \abs{\L_{\mathrm{FP}}^{-1} (\tfrac{g_1}{M_0})}^2 M_0}^\frac{1}{2},
  \end{align}
from which, it follows, for some constants $a_1$ and $\tilde{a}_1$, that
  \begin{align}
    \agl{\tfrac{|h_1|^2}{M_0}} \ge \tilde{a}_1 \agl{\tfrac{|g_1|^2}{M_0}} \,,
  \end{align}
and hence, that
  \begin{align}\label{eq:h-larger-f}
    \agl{\tfrac{|h_1|^2}{M_0}} \ge a_1 \agl{\tfrac{|f_1|^2}{M_0}} \,,
  \end{align}
because we can get, by performing straightforward calculations, that
  \begin{align*}
    & \skp{\tfrac{g_1}{M_0}}{\tfrac{g_1}{M_0}}_M
      = \skp{\sum_k c_k v_{\bot}^k}{\sum_l c_l v_{\bot}^l}_M
      = \agl{\sum_{k,l} c_k c_l \tfrac{v_{\bot}^k v_{\bot}^l}{\abs{v_{\bot}}} M_0 }
      = c_0 \sum_k c_k^2 \,, \no\\
    &
      \skp{\tfrac{g_1}{M_0}}{\L_{\mathrm{FP}}^{-1} (\tfrac{g_1}{M_0})}_M
      = \skp{\sum_k c_k v_{\bot}^k}{\sum_l c_l \tfrac{\chi(v)}{\abs{v_{\bot}}} v_{\bot}^l}_M
      = \agl{\sum_{k,l} c_k c_l \chi(v) \tfrac{v_{\bot}^k v_{\bot}^l}{\abs{v_{\bot}}} M_0 }
      = \tilde{c}_0 \sum_k c_k^2 \,, \no\\
    &
      \agl{ \abs{\L_{\mathrm{FP}}^{-1} (\tfrac{g_1}{M_0})}^2 M_0}
      = \skp{\sum_k c_k \tfrac{\chi(v)}{\abs{v_{\bot}}} v_{\bot}^k}{\sum_l c_l \tfrac{\chi(v)}{\abs{v_{\bot}}} v_{\bot}^l}_M
      = \bar{c}_0 \sum_k c_k^2 \,,
  \end{align*}
with $\bar{c}_0$ denoting a fixed constant.

Denoting $\eta_0 = (a_1 \eta_1)^2$, and combining together \cref{eq:low-bdd-f-diff,eq:low-bdd-h,eq:h-larger-f}, the H\"older inequality
  \begin{align*}
    \skp{ \tfrac{f_1}{M_0} - c_0^{-1} (v_\bot j_{f_1})}{f_1} \le \agl{\abs{ \tfrac{f_1 - c_0^{-1} (v_\bot j_1)M_0}{M_0} }^2 }_M^\frac{1}{2} \agl{\tfrac{|f_1|^2}{M_0}}^\frac{1}{2}
  \end{align*}
leads us to the desired coercivity estimate \cref{eq:coercivity-I}, for $f_1 = f_1^\perp$, that:
  \begin{align}
    \skp{-\L_{\Omega_0} f_1}{\frac{f_1 - c_0^{-1} (v_\bot j_1)M_0}{M_0}}
    \ge \lambda_0 \agl{\abs{ \frac{f_1 - c_0^{-1} (v_\bot j_1)M_0}{M_0} }^2 }_M
    \ge \lambda_0 \eta_0 \agl{\frac{|f_1|^2}{M_0}}.
  \end{align}

\endproof
%% ------------------------------------------------------------------------
% \subsubsection{Coercivity estimate of Type II} % (fold)
% \label{ssub:coercivity_estimate_of_type_ii}

% % subsubsection coercivity_estimate_of_type_ii (end)

\smallskip\noindent\underline{\textbf{Coercivity estimate of Type II}}.
We need the coefficients hypothesis $\alpha > 4^3 (c_0^{-1} \lambda_0^{-1} \beta^{1/2} (d+2)^\frac{3}{2} +2)$ in potential function $V(|v|)= \beta \tfrac{|v|^4}{4} - \alpha \tfrac{|v|^2}{2}$, in order to provide more dissipation effect from the $\alpha$-term. An appropriate large $\alpha$ suffices to satisfy the requirement.

\begin{proposition}[Coercivity estimate, II] \label{prop:coercivity}
  Suppose that $\alpha > 4^3 (c_0^{-1} \lambda_0^{-1} \beta^{1/2} (d+2)^\frac{3}{2} +2)$, then we have the coercivity estimate:
  \begin{align} \label{eq:coercivity}
    \skp{-\L_{\Omega_0} f}{\tfrac{f}{M_0}}
    \ge \tfrac{3}{4} \agl{\abs{\nabla_v \left( \tfrac{f}{M_0} \right)}^2}_M + \tfrac{3\alpha}{4c_0} \abs{P_{\Omega_0^\bot} {j_f}}^2 \,.
  \end{align}

\end{proposition}

\proof
By recalling \cref{eq:linearized-op-sep}, we write the linear operator $\L_{\Omega_0} $ as:
  \begin{align}
    \L_{\Omega_0} f
    & = \div_v \left\{ M_0 \nabla_v \left( \tfrac{f}{M_0} \right) - c_0^{-1} P_{\Omega_0^\bot} {j_f} M_0 \right\} \\ \no
    & = \div_v \left\{ M_0 \nabla_v \left( \tfrac{f}{M_0} \right) \right\} - c_0^{-1} P_{\Omega_0^\bot} {j_f} \nabla_v M_0 \\\no
    & = - M_0 \L_{\mathrm{FP}} \left( \tfrac{f}{M_0} \right) - c_0^{-1} P_{\Omega_0^\bot} {j_f} \nabla_v M_0 \,.
  \end{align}
So we have, by the fact $\iint \nabla_v M_0 \rho_f(t,x) \d v \d x=0$,
  \begin{align}
    \skp{-\L_{\Omega_0} f}{\frac{f}{M_0}} % & = \skp{-\L_{\Omega_0} f}{\frac{f}{M_0}-\rho_f} \\\no
    & = \skp{-\div_v \left\{ M_0 \nabla_v \left( \frac{f}{M_0} \right) \right\} }{\frac{f}{M_0}}
        + \skp{\div_v \left( c_0^{-1} P_{\Omega_0^\bot} {j_f} M_0 \right) }{\frac{f}{M_0}}
      \\\no
    & = \skp{M_0 \L_{\mathrm{FP}} \left( \tfrac{f}{M_0} \right) }{\tfrac{f}{M_0}} + \skp{ c_0^{-1} P_{\Omega_0^\bot} {j_f} \nabla_v M_0 }{\tfrac{f}{M_0} - \rho_f(t,x)} \\\no
    & = \agl{M_0 \abs{\nabla_v \left( \tfrac{f}{M_0} \right)}^2} + \skp{c_0^{-1} P_{\Omega_0^\bot} {j_f} \nabla_v M_0}{\tfrac{f}{M_0} - \rho_f(t,x)}.
  % \\\no
  %   & = \agl{M_0 \abs{\nabla_v \left( \tfrac{f}{M_0} \right)}^2} + c_0^{-1} P_{\Omega_0^\bot} {j_f} \cdot \left[ (\alpha-1) (\agl{v_{\bot} f} - \rho_f \agl{v_{\bot} M_0}) - \agl{\beta |v|^2 v_{\bot} (f - \rho_f M_0)} \right] \\\no
  %   & = \agl{M_0 \abs{\nabla_v \left( \tfrac{f}{M_0} \right)}^2} + c_0^{-1} (\alpha-1) \abs{P_{\Omega_0^\bot} {j_f}}^2 - \agl{\beta |v|^2 v_{\bot} (f - \rho_f M_0)}
  \end{align}

The factor $\nabla_v M_0$ can be calculated as, due to the expression $V(|v|)= \beta \tfrac{|v|^4}{4} - \alpha \tfrac{|v|^2}{2}$:
  \begin{align*}
    \nabla_v M_0/M_0 = - (v- \Omega_0) - \nabla_v V'(|v|) \tfrac{v}{|v|}
    = - (v- \Omega_0) - (\beta |v|^2 - \alpha)v
    = (\alpha-1)v - \beta |v|^2 v + \Omega_0,
  \end{align*}
which immediately implies,
  \begin{align}
    & \skp{P_{\Omega_0^\bot} {j_f} \nabla_v M_0}{\frac{f}{M_0} - \rho_f(t,x)} \\\no
    & = P_{\Omega_0^\bot} {j_f} \cdot \left[ (\alpha-1) (\agl{v_{\bot} f} - \rho_f \agl{v_{\bot} M_0}) - \agl{\beta |v|^2 v_{\bot} (f - \rho_f M_0)} \right] \\\no
    & = (\alpha-1) \abs{P_{\Omega_0^\bot} {j_f}}^2 - P_{\Omega_0^\bot} {j_f} \cdot \agl{\beta |v|^2 v_{\bot} (f - \rho_f M_0)},
  \end{align}
where we have used $\agl{v_{\bot} M_0} =0$ for odd factor $v_{\bot}$.

As for the last term in above equality, we firstly write:
  \begin{align*}
    \agl{\beta |v|^2 v_{\bot} (f - \rho_f M_0)} \le \agl{M_0}^{\frac{1}{2}} \cdot \agl{\beta^2 |v|^6 \abs{\tfrac{f}{M_0} - \rho_f}^2 M_0}^{\frac{1}{2}}
    \le \agl{\abs{\nabla_v \Phi}^2 \abs{\tfrac{f}{M_0} - \rho_f}^2 }_M^{\frac{1}{2}}.
  \end{align*}
The vanishing mean value $\agl{\tfrac{f}{M_0} - \rho_f}_M =0$ enables us to get the following Poincar\'e inequality:
  \begin{align}
    \agl{\abs{\nabla_v \Phi}^2 \abs{\tfrac{f}{M_0} - \rho_f}^2 }_M \le C_{\lambda} \agl{\abs{\nabla_v \left( \tfrac{f}{M_0} \right)}^2}_M.
  \end{align}
This Poincar\'e type inequality will be proven in \Cref{lemm:Poincare-ineq} below, in the spirit of \cite{LLZ-07cpam}. We have used $\Delta \Phi \le \tfrac{1}{2}\abs{\nabla_v \Phi}^2 + 2 \beta^{\frac{1}{2}} (d+2)^\frac{3}{2}$. The above Poincar\'e inequlity yields that,
  \begin{align*}
    c_0^{-1} P_{\Omega_0^\bot} {j_f} \cdot \agl{\beta |v|^2 v_{\bot} (f - \rho_f M_0)}
    & \le c_0^{-1} C_{\lambda}^{\frac{1}{2}} \abs{P_{\Omega_0^\bot} {j_f}} \cdot \agl{\abs{\nabla_v \left( \tfrac{f}{M_0} \right)}^2}_M^{\frac{1}{2}} \\\no
    & \le \tfrac{1}{4} \agl{\abs{\nabla_v \left( \tfrac{f}{M_0} \right)}^2}_M + c_0^{-2} C_{\lambda} \abs{P_{\Omega_0^\bot} {j_f}}^2.
  \end{align*}
Notice that the assumption $\alpha> 4^3 (c_0^{-1} \lambda_0^{-1} \beta^{\frac{1}{2}} (d+2)^\frac{3}{2} +2)$ yielding $(\tfrac{\alpha}{4} -1) > 4^2 c_0^{-1} (\lambda_0^{-1} \beta^{\frac{1}{2}} (d+2)^\frac{3}{2} +1) = c_0^{-1} C_{\lambda}$, we get finally
  \begin{align}
    \skp{-\L_{\Omega_0} f}{\tfrac{f}{M_0}}
    & = \agl{M_0 \abs{\nabla_v \left( \tfrac{f}{M_0} \right)}^2} + c_0^{-1} (\alpha-1) \abs{P_{\Omega_0^\bot} {j_f}}^2 - c_0^{-1} P_{\Omega_0^\bot} {j_f} \cdot \agl{\beta |v|^2 v_{\bot} (f - \rho_f M_0)}
  \no \\
    & \ge \tfrac{3}{4} \agl{\abs{\nabla_v \left( \tfrac{f}{M_0} \right)}^2}_M + \tfrac{3\alpha}{4c_0} \abs{P_{\Omega_0^\bot} {j_f}}^2.
  \end{align}
This is the coercivity estimate \cref{eq:coercivity} we are looking for, in which the first-moment $j_f = \agl{vf}$ appears as a dissipation contribution alone.
\endproof

%% ------------------------------------------------------------------------
\subsection{Some useful lemmas} % (fold)
\label{sub:some_useful_lemmas}

% subsection some_useful_lemmas (end)

\subsubsection{The weighted Poincar\'e inequalities} % (fold)
\label{ssub:weighted_poincare_inequality}

% subsubsection weighted_poincare_inequality (end)

We state here the following (weighted) Poincar\'e inequality lemma, helping us deal with the contribution from $g = \tfrac{f_\R}{M_0}$ or its higher-order derivatives $g = \tfrac{\nabla_x^s f_\R}{M_0}$.

% The key part of its proof is to transfer the considered space $L^2(\d q \!\d x)$ into a weighted one $L^2(M\!\d q \!\d x)$ with respect to a Maxwellian weight $M=M(q)$. The proofs are similar as that for Lemmas 3.2--3.3 in \cite{LLZ-07cpam} (or Lemma 1.6 in \cite{JLZ-18sima}), except that an additional mean value term appeared in the present case.

\begin{lemma}[Poincar\'e inequality] \label{lemm:Poincare-ineq}
  We have the weighted Poincar\'e inequalities for any function satisfying $\agl{g}_M = 0$, that
  \begin{align} \label{eq:Poinc-ineq-g}
    \int \abs{g}^2 M_0 \d v & \le \lambda_0 \int \abs{\nabla_v g}^2 M_0 \d v , \\ \label{eq:Poinc-ineq-Phi}
    \int \abs{\nabla_v \Phi}^2 \abs{g}^2 M_0 \d v & \le C_\lambda \int \abs{\nabla_v g}^2 M_0 \d v ,
    % \nm{\nabla_v f}_{\lxv}^2 &
    %   \ls \iint \abs{ \nabla_v \left( \tfrac{f}{\sqrt M} \right) }^2 M \d q \d x + \abs{\rho}_{\lx}^2, \\
    % \nm{q\nabla_q U f}_{\lxv}^2 & % = \iint |\nabla_q U (qg)|^2 M \d q \d x
    %   \ls \iint \abs{ \agl{q} \nabla_q \left( \tfrac{f}{\sqrt M} \right) }^2 M \d q \d x + \abs{\rho}_{\lx}^2,
  \end{align}
with $C_\lambda = 4^2 (\lambda_0^{-1} \beta^{\frac{1}{2}} (d+2)^\frac{3}{2} +1) $.

\end{lemma}

\begin{remark}
  Indeed, we will not use directly the Poincar\'e inequalities to deal with $g = \tfrac{f_\R}{M_0}$, but instead of a variant formulation of $g = \tfrac{f_\R}{M_0} - \agl{\tfrac{f_\R}{M_0}}_M = \tfrac{f_\R}{M_0} - \rho_\R(t,x)$, due to the requirement $\agl{g}_M =0$. However, the quantities involving higher-order derivatives could not be dealed with similarly, since a similar structure is missing.
\end{remark}

\proof
The spirit of the proof is similar as that in \cite[Lemmas 3.2]{LLZ-07cpam} (or, see \cite[Lemma 1.6]{JLZ-18sima}). The first result is exactly the weighted Poincar\'e inequality with constant $\lambda_0$ as the first positive eigenvalue of the Schr\"odinger operator, see \cite[Lemma 2.1]{ABCD-19mbe}. Here we focus on the second result \cref{eq:Poinc-ineq-Phi}.

By performing a straightforward calculus, it follows:
  \begin{align*}
    \Delta_v \Phi & = d (1- \alpha) + \beta (d+2) |v|^2, \\
    \abs{\nabla_v \Phi}^2 & = \left[ (\beta |v|^2 + 1- \alpha) v - \Omega_0 \right]^2
    \ge \tfrac{1}{2} \left[ (\beta |v|^2 + 1- \alpha) v \right]^2 - |\Omega_0|^2
    = \tfrac{1}{2} (\beta |v|^2 + 1- \alpha)^2 |v|^2 -1.
  \end{align*}
Denote $X = \beta |v|^2 + 1- \alpha$, then we get, $X \le \beta |v|^2$ for largely $|v|$, and moreover,
  \begin{align*}
    \tfrac{1}{d+2} \Delta_v \Phi \le X, \quad \tfrac{1}{2} \abs{\nabla_v \Phi}^2 \ge \tfrac{1}{4\beta} X^3.
  \end{align*}
The fundamental Young inequality $(d+2) X \le \tfrac{1}{4\beta} X^3 + 2 \beta^\frac{1}{2} (d+2)^\frac{3}{2} $ leads us to
  \begin{align}
    \Delta_v \Phi \le \tfrac{1}{2} \abs{\nabla_v \Phi}^2 + 2 \beta^\frac{1}{2} (d+2)^\frac{3}{2} .
  \end{align}

By noticing the fact $\nabla_v M_0 = -\nabla_v \Phi M_0$, and by integration by parts, we can infer
  \begin{align}\label{eq:potential-poinc}
    \int \abs{\nabla_v \Phi}^2 \abs{g}^2 M_0 \d v
    & = - \int \nabla_v \Phi g^2 \nabla_v M_0 \d v \\\no
    & = \int (\Delta_v \Phi g^2 + 2 \nabla_v \Phi g \nabla_v g) M_0 \d v \\\no
    & \le \int \left( \tfrac{1}{2} \abs{\nabla_v \Phi}^2 + 2 \beta^\frac{1}{2} (d+2)^\frac{3}{2} \right)g^2 M_0 \d v \\\no & \qquad
      + \tfrac{1}{4} \int \abs{\nabla_v \Phi}^2 \abs{g}^2 M_0 \d v + 4\int \abs{\nabla_v g}^2 M_0 \d v \\\no
    & \le \tfrac{3}{4} \int \abs{\nabla_v \Phi}^2 \abs{g}^2 M_0 \d v + \left[ 2 \lambda_0^{-1} \beta^\frac{1}{2} (d+2)^\frac{3}{2} +4 \right] \int \abs{\nabla_v g}^2 M_0 \d v,
  \end{align}
where we have used in the penultimate line the H\"older inequality, and the first result \cref{eq:Poinc-ineq-g} in the last line. Therefore we get,
  \begin{align}
    \int \abs{\nabla_v \Phi}^2 \abs{g}^2 M_0 \d v \le 4 \left[ 2 \lambda_0^{-1} \beta^\frac{1}{2} (d+2)^\frac{3}{2} +4 \right] \int \abs{\nabla_v g}^2 M_0 \d v,
  \end{align}
which is exactly the desired result \cref{eq:Poinc-ineq-Phi}, with the constant $C_\lambda = 4^2 (\lambda_0^{-1} \beta^{\frac{1}{2}} (d+2)^\frac{3}{2} +1) $. We also refer it as the Poincar\'e inequality in the sequel.
\endproof

%% ------------------------------------------------------------------------ %%
\subsubsection{Estimates on remainders in Taylor's formula} % (fold)
\label{ssub:Taylor-remainder}

% subsubsection Taylor-remainder (end)

For the convenience of readers, we give some common results concerning the product estimate and the composition action of smooth function, both of which will be frequently used in the rest of the paper.
\begin{lemma}[\cite{Majda-84b,Taylor-11b}]\label{lemm:Sobolev-product}
  We have the estimate for products in Sobolev spaces:
    \begin{align}\label{eq:Sobolev-product}
      \nm{f g}_{H^s} \le C \nm{f}_{H^s} \nm{g}_{L^\infty} + C \nm{f}_{L^\infty} \nm{g}_{H^s} \,.
    \end{align}
  Let $F$ a smooth function satisfying $F(0)=0$. Then for $u \in H^s \cap L^\infty$, there exists a constant $C$ depending only on $s$ and $\nm{u}_{L^\infty} $, such that
    \begin{align}\label{eq:smooth-function}
      \nm{F(u)}_{H^s} \le C \nm{u}_{H^s}\,.
    \end{align}
\end{lemma}

In order to deal with the remainders in Taylor's formula, we need the following lemma:
\begin{lemma}\label{lemm:Taylor-remainders}
  We have,
  \begin{align}\label{esm:Taylor-remainders-1R}
    \nm{\tfrac{( \Omega_{f^\eps} - \Omega_0 )}{\eps}}_{H^s_x}
    & \le C_{j_0,j_1,j_2} (1 + \eps \nm{j_\R}_{H^s_x} + \eps^3 \nm{j_\R}_{H^s_x} \nm{j_\R}_{H^2_x}) \\\no
    & \le C_{f_0,f_1,f_2} \left[ 1 + \eps \nm{f_\R}_{\H^s_xL^2_v(M^{-1})} + \eps^3 \nm{f_\R}_{\H^s_xL^2_v(M^{-1})} \nm{f_\R}_{\H^2_x L^2_v(M^{-1})} \right] \,,
    \\
    \label{esm:Taylor-remainders-3R}
    \nm{\widehat{\Omega}_{3,R}}_{\hx{s}}
    & \le C_{f_0,f_1,f_2} (1 + \nm{j_\R}_{\hx{s}} + \eps \nm{j_\R}_{\hx{s}} \nm{j_\R}_{\hx{2}} + \eps^3 \nm{j_\R}_{\hx{s}} \nm{j_\R}^2_{\hx{2}}) \\\no
    & \le C_{f_0,f_1,f_2} \Big[ 1+ \nm{f_\R}_{\H^s_xL^2_v(M^{-1})} + \eps \nm{f_\R}_{\H^s_xL^2_v(M^{-1})} \nm{f_\R}_{\H^2_x L^2_v(M^{-1})} \\\no & \qquad \qquad \quad + \eps^3 \nm{f_\R}_{\H^s_xL^2_v(M^{-1})} \nm{f_\R}^2_{\H^2_x L^2_v(M^{-1})} \Big] \,.
  \end{align}

\end{lemma}

\proof
The quantity $\tfrac{(\Omega_{f^\eps} - \Omega_0)}{\eps}$ is exactly the remainder of order 1 in Taylor's formula, since we have
  \begin{align}
   \Omega_{f^\eps} = \Omega_0 + \eps \tfrac{\d}{\d \eps} \Omega_{f^\eps} \big|_{\eps = \eta}.
  \end{align}
By employing a shorthand notation $j_\eps = j_{f^\eps}$, it follows
  \begin{align}\label{eq:Omega-Remd-1}
    \frac{\d}{\d \eps} \Omega_{f^\eps}^i & = \frac{\d}{\d \eps} \frac{j_{\eps}^i}{\abs{j_{\eps}}}
    = \left( \delta_{ik} - \frac{j_{\eps}^i}{\abs{j_{\eps}}} \frac{j_{\eps}^k}{\abs{j_{\eps}}} \right) \frac{ \tfrac{\d}{\d \eps} j_\eps^k }{\abs{j_\eps}}
    = P_{\Omega_{f^\eps}^\bot}^{ik} \frac{j_1^k + \eps (j_2 + j_R)^k}{\abs{j_\eps}}
    \\\no
    & = P_{\Omega_0^\bot}^{ik} \frac{j_1^k + \eps (j_2 + j_R)^k}{\abs{j_0}} + \left( \frac{j_0^i}{\abs{j_0}} \frac{j_0^k}{\abs{j_0}} \frac{1}{\abs{j_0}} - \frac{j_{\eps}^i}{\abs{j_{\eps}}} \frac{j_{\eps}^k}{\abs{j_{\eps}}} \frac{1}{\abs{j_{\eps}}} \right) \cdot [j_1^k + \eps (j_2 + j_R)^k] \\\no
    & = F_1^i (\eps) + F_2^i (\eps) \cdot [j_1^k + \eps (j_2 + j_R)^k].
  \end{align}
The factor $F_2(\eps)$ in the second term, satisfying $F_2(0)=0$, can be bounded by the composition action of smooth function in \Cref{lemm:Sobolev-product}:
  \begin{align}
    \nm{F_2(\eps)}_{H^s_x} \le C \nm{j_\eps}_{H^s_x} \le C \nm{j_0+ \eps j_1 + \eps^2 j_2}_{H^s_x} + C \eps^2 \nm{j_R}_{H^s_x},
  \end{align}
with the coefficient $C$ depends on the bound of $\nm{j_\eps}_{L^\infty_x}$ from up and especially from below. By using the expansion formula $j_\eps = j_0 + \eps j_1 + \eps^2 j_2 + \eps^2 j_\R$, the upper and lower bound could be established if they hold initially, due to the smallness of parameter $\eps$ and the desired boundness of $\eps \nm{j_\R}_{L^\infty_x} \le \eps \nm{f_\R}_{\H^2_x L^2_v(M^{-1})} \le C $. We give a more precise description in \Cref{prop:nonzero-j} below.

This inequality, combined with the control for the first term $\nm{F_1(\eps)}_{H^s_x} \le C \nm{j_1 + \eps j_2}_{H^s_x} + C \eps \nm{j_R}_{H^s_x}$, yields the first conclusion in \Cref{lemm:Taylor-remainders}:
  \begin{align}
    \nm{\tfrac{( \Omega_{f^\eps} - \Omega_0 )}{\eps}}_{H^s_x}
    & = \nm{ \left( \tfrac{\d}{\d \eps} \Omega_{f^\eps} \right) \big|_{\eps = \eta} }_{H^s_x} \\\no
    & \le \nm{F_1^i (\eta)}_{H^s_x}
          + \nm{ F_2^i (\eta)}_{H^s_x} \nm{j_1 + \eta (j_2 + j_\R)}_{H^2_x}
          + \nm{ F_2^i (\eta)}_{H^2_x} \nm{j_1 + \eta (j_2 + j_\R)}_{H^s_x} \\\no
    & \le C_{j_0,j_1,j_2} (1 + \eps \nm{j_\R}_{H^s_x})
          + C_{j_0,j_1,j_2} (1 + \eps^2 \nm{j_\R}_{H^s_x}) (1 + \eta \nm{j_\R}_{H^2_x}) \\\no & \quad
          + C_{j_0,j_1,j_2} (1 + \eps^2 \nm{j_\R}_{H^2_x}) (1 + \eta \nm{j_\R}_{H^s_x}) \\\no
    & \le C_{j_0,j_1,j_2} (1 + \eps \nm{j_\R}_{H^s_x} + \eps^3 \nm{j_\R}_{H^s_x} \nm{j_\R}_{H^2_x}) \\\no
    & \le C_{f_0,f_1,f_2} \left[ 1 + \eps \nm{f_\R}_{\H^s_xL^2_v(M^{-1})} + \eps^3 \nm{f_\R}_{\H^s_xL^2_v(M^{-1})} \nm{f_\R}_{\H^2_x L^2_v(M^{-1})} \right] \,.
  \end{align}
where we have used the product estimates in Sobolev spaces $H^s_x$ ($s \ge 2$), stated in \Cref{lemm:Sobolev-product}.

For the second result in \Cref{lemm:Taylor-remainders}, notice it is the remainder of order 3 in Taylor's formula. So its proof is in the same spirit with above process, but more complicated. Recalling the expansion formula in \cref{eq:Omega-ansatz}, we have
  \begin{align}\label{eq:RemD3}
    \widehat{\Omega}_{3,R} = \frac{1}{3!}\frac{\d^3}{\d \eps^3} \Omega_{f^\eps} \big|_{\eps = \eta}.
  \end{align}

Performing straightforward calculations implies:
  \begin{align}
    \frac{\d}{\d \eps} \Omega_{f^\eps}^i
    = \frac{\d}{\d \eps} \left( \frac{j_{\eps}^i}{\abs{j_{\eps}}} \right)
    = \left( \delta_{ik} - \frac{j_{\eps}^i}{\abs{j_{\eps}}} \frac{j_{\eps}^k}{\abs{j_{\eps}}} \right)
        \frac{ \tfrac{\d}{\d \eps} j_\eps^k }{\abs{j_\eps}}
    = P_{\Omega_{f^\eps}^\bot}^{ik} \frac{ \tfrac{\d}{\d \eps} j_\eps^k }{\abs{j_\eps}}
    = \frac{ \tfrac{\d}{\d \eps} j_\eps^i }{\abs{j_\eps}}
      - \frac{j_{\eps}^i}{\abs{j_{\eps}}} \cdot \frac{ \tfrac{\d}{\d \eps} (\tfrac{1}{2} \abs{j_\eps}^2) }{\abs{j_\eps}^2} \,,
  \end{align}
and
  \begin{align}
    \frac{\d^2}{\d \eps^2} \Omega_{f^\eps}^i & = \frac{\d}{\d \eps} \left( P_{\Omega_{f^\eps}^\bot}^{ik} \frac{ \tfrac{\d}{\d \eps} j_\eps^k }{\abs{j_\eps}} \right)
    = \frac{\d}{\d \eps} \left( \frac{ \tfrac{\d}{\d \eps} j_\eps^i }{\abs{j_\eps}} \right)
      - \frac{\d}{\d \eps} \left[ \frac{j_{\eps}^i}{\abs{j_{\eps}}} \cdot \frac{ \tfrac{\d}{\d \eps} (\tfrac{1}{2} \abs{j_\eps}^2) }{\abs{j_\eps}^2} \right] \\\no
    & =\left( \frac{ \tfrac{\d^2}{\d \eps^2} j_\eps^i }{\abs{j_\eps}} - \frac{ \tfrac{\d}{\d \eps} j_\eps^i }{\abs{j_\eps}} \cdot \frac{ \tfrac{\d}{\d \eps} (\tfrac{1}{2} \abs{j_\eps}^2) }{\abs{j_\eps}^2} \right)
        - \left( \frac{ \tfrac{\d}{\d \eps} j_\eps^i }{\abs{j_\eps}} - \frac{j_{\eps}^i}{\abs{j_{\eps}}} \frac{ \tfrac{\d}{\d \eps} (\tfrac{1}{2} \abs{j_\eps}^2) }{\abs{j_\eps}^2} \right) \frac{ \tfrac{\d}{\d \eps} (\tfrac{1}{2} \abs{j_\eps}^2) }{\abs{j_\eps}^2}  \\\no & \qquad
        - \frac{j_{\eps}^i}{\abs{j_{\eps}}} \left[ \frac{ \tfrac{\d^2}{\d \eps^2} (\tfrac{1}{2} \abs{j_\eps}^2) }{\abs{j_\eps}^2} - 2 \left( \frac{ \tfrac{\d}{\d \eps} (\tfrac{1}{2} \abs{j_\eps}^2) }{\abs{j_\eps}^2} \right)^2 \right]
  \\\no
    & = \frac{ \tfrac{\d^2}{\d \eps^2} j_\eps^i }{\abs{j_\eps}}
        - \frac{j_{\eps}^i}{\abs{j_{\eps}}} \frac{ \tfrac{\d^2}{\d \eps^2} (\tfrac{1}{2} \abs{j_\eps}^2) }{\abs{j_\eps}^2}
        % - \frac{j_{\eps}^i}{\abs{j_{\eps}}} \cdot \left( \frac{ \tfrac{\d}{\d \eps} j_\eps }{\abs{j_\eps}} \right)^2
        - 2 \frac{ \tfrac{\d}{\d \eps} j_\eps^i }{\abs{j_\eps}} \cdot \frac{ \tfrac{\d}{\d \eps} (\tfrac{1}{2} \abs{j_\eps}^2) }{\abs{j_\eps}^2}
        + 3 \frac{j_{\eps}^i}{\abs{j_{\eps}}} \left( \frac{ \tfrac{\d}{\d \eps} (\tfrac{1}{2} \abs{j_\eps}^2) }{\abs{j_\eps}^2} \right)^2
  \\\no
    & = P_{\Omega_{f^\eps}^\bot}^{ik} \frac{ \tfrac{\d^2}{\d \eps^2} j_\eps^k }{\abs{j_\eps}}
        - \frac{j_{\eps}^i}{\abs{j_{\eps}}} \cdot \left( \frac{ \tfrac{\d}{\d \eps} j_\eps }{\abs{j_\eps}} \right)^2
        - 2 \frac{ \tfrac{\d}{\d \eps} j_\eps^i }{\abs{j_\eps}} \cdot \frac{ \tfrac{\d}{\d \eps} (\tfrac{1}{2} \abs{j_\eps}^2) }{\abs{j_\eps}^2}
        + 3 \frac{j_{\eps}^i}{\abs{j_{\eps}}} \left( \frac{ \tfrac{\d}{\d \eps} (\tfrac{1}{2} \abs{j_\eps}^2) }{\abs{j_\eps}^2} \right)^2
  \\\no
    & = P_{\Omega_{f^\eps}^\bot}^{ik} \frac{ \tfrac{\d^2}{\d \eps^2} j_\eps^k }{\abs{j_\eps}}
        + \Q_2 \left[
              \frac{j_{\eps}^i}{\abs{j_{\eps}}}
                \left( \tfrac{ \tfrac{\d}{\d \eps} j_\eps }{\abs{j_\eps}} \right)^2,\,
              \frac{j_{\eps}^i}{\abs{j_{\eps}}}
                \left( \tfrac{ \tfrac{\d}{\d \eps} (\tfrac{1}{2} \abs{j_\eps}^2) }{\abs{j_\eps}^2} \right)^2,\,
              \frac{ \tfrac{\d}{\d \eps} j_\eps^i }{\abs{j_\eps}}
                \frac{ \tfrac{\d}{\d \eps} (\tfrac{1}{2} \abs{j_\eps}^2) }{\abs{j_\eps}^2}
              \right].
  \end{align}
Here we have used the basic fact in the penultimate line, that
  \begin{align*}
    \tfrac{\d^2}{\d \eps^2} (\tfrac{1}{2} \abs{j_\eps}^2)
    = \tfrac{\d}{\d \eps} \left[ \tfrac{\d}{\d \eps} (\tfrac{1}{2} \abs{j_\eps}^2) \right]
    = \frac{\d}{\d \eps} \left( j_\eps^k \cdot \tfrac{\d}{\d \eps} j_\eps^k \right)
    = j_\eps^k \cdot \tfrac{\d^2}{\d \eps^2} j_\eps^k + \left( \tfrac{\d}{\d \eps} j_\eps^k \right)^2.
  \end{align*}
We also denoted by $\Q_2$ a homogeneous polynomial formulation of order 2 along with the $i$-th fixed direction for $i \in \left\{ 1,2, \cdots, d \right\}$, with respect to its derivative arguments, but except for its highest derivatives on any single argument. For example, the first argument contains a zero-derivative factor $\frac{j_{\eps}^i}{\abs{j_{\eps}}}$ and a quadratic product of one-derivative factor $\tfrac{ \tfrac{\d}{\d \eps} j_\eps }{\abs{j_\eps}}$, and thus is referred here as containing two derivatives.

By virtue of above formulation, and by noticing the fact $\tfrac{\d^3}{\d \eps^3} (\tfrac{1}{2} \abs{j_\eps}^2) = j_\eps^k \cdot \tfrac{\d^3}{\d \eps^3} j_\eps^k + 3 \tfrac{\d}{\d \eps} j_\eps^k \cdot \tfrac{\d^2}{\d \eps^2} j_\eps^k$, we can calculate:
  \begin{align}
    & \frac{\d^3}{\d \eps^3} \Omega_{f^\eps}^i \no\\
    & = \frac{\d}{\d \eps}
      \left\{
        \frac{ \tfrac{\d^2}{\d \eps^2} j_\eps^i }{\abs{j_\eps}}
        - \frac{j_{\eps}^i}{\abs{j_{\eps}}} \frac{ \tfrac{\d^2}{\d \eps^2} (\tfrac{1}{2} \abs{j_\eps}^2) }{\abs{j_\eps}^2}
        - 2 \frac{ \tfrac{\d}{\d \eps} j_\eps^i }{\abs{j_\eps}} \cdot \frac{ \tfrac{\d}{\d \eps} (\tfrac{1}{2} \abs{j_\eps}^2) }{\abs{j_\eps}^2}
        + 3 \frac{j_{\eps}^i}{\abs{j_{\eps}}} \left( \frac{ \tfrac{\d}{\d \eps} (\tfrac{1}{2} \abs{j_\eps}^2) }{\abs{j_\eps}^2} \right)^2
      \right\}
  \no\\
    & = \frac{ \tfrac{\d^3}{\d \eps^3} j_\eps^i }{\abs{j_\eps}}
        - \frac{j_{\eps}^i}{\abs{j_{\eps}}} \frac{ \tfrac{\d^3}{\d \eps^3} (\tfrac{1}{2} \abs{j_\eps}^2) }{\abs{j_\eps}^2}
        - 3 \frac{ \tfrac{\d^2}{\d \eps^2} j_\eps^i }{\abs{j_\eps}} \cdot \frac{ \tfrac{\d}{\d \eps} (\tfrac{1}{2} \abs{j_\eps}^2) }{\abs{j_\eps}^2}
        - 3 \frac{ \tfrac{\d}{\d \eps} j_\eps^i }{\abs{j_\eps}} \cdot \frac{ \tfrac{\d^2}{\d \eps^2} (\tfrac{1}{2} \abs{j_\eps}^2) }{\abs{j_\eps}^2} \no\\ & \qquad
        + 9 \frac{j_{\eps}^i}{\abs{j_{\eps}}} \frac{ \tfrac{\d}{\d \eps} (\tfrac{1}{2} \abs{j_\eps}^2) }{\abs{j_\eps}^2} \frac{ \tfrac{\d^2}{\d \eps^2} (\tfrac{1}{2} \abs{j_\eps}^2) }{\abs{j_\eps}^2}
        + 9 \frac{\tfrac{\d}{\d \eps} j_{\eps}^i}{\abs{j_{\eps}}} \left( \frac{ \tfrac{\d}{\d \eps} (\tfrac{1}{2} \abs{j_\eps}^2) }{\abs{j_\eps}^2} \right)^2
        - 15 \frac{j_{\eps}^i}{\abs{j_{\eps}}} \left( \frac{ \tfrac{\d}{\d \eps} (\tfrac{1}{2} \abs{j_\eps}^2) }{\abs{j_\eps}^2} \right)^3
  \no\\ \label{eq:remainder-3R}
    & = P_{\Omega_{f^\eps}^\bot}^{ik} \frac{ \tfrac{\d^3}{\d \eps^3} j_\eps^k }{\abs{j_\eps}}
        - 3 \frac{j_{\eps}^i}{\abs{j_{\eps}}}
              \left( \frac{\tfrac{\d}{\d \eps} j_\eps^k}{\abs{j_\eps}}
              \cdot \frac{\tfrac{\d^2}{\d \eps^2} j_\eps^k}{\abs{j_\eps}} \right)
        - 3 \frac{ \tfrac{\d^2}{\d \eps^2} j_\eps^i }{\abs{j_\eps}}
              \cdot \frac{ \tfrac{\d}{\d \eps} (\tfrac{1}{2} \abs{j_\eps}^2) }{\abs{j_\eps}^2}
        - 3 \frac{ \tfrac{\d}{\d \eps} j_\eps^i }{\abs{j_\eps}}
              \cdot \frac{ \tfrac{\d^2}{\d \eps^2} (\tfrac{1}{2} \abs{j_\eps}^2) }{\abs{j_\eps}^2}
            \\\no & \qquad
        + 9 \frac{j_{\eps}^i}{\abs{j_{\eps}}}
              \frac{ \tfrac{\d}{\d \eps} (\tfrac{1}{2} \abs{j_\eps}^2) }{\abs{j_\eps}^2}
              \frac{ \tfrac{\d^2}{\d \eps^2} (\tfrac{1}{2} \abs{j_\eps}^2) }{\abs{j_\eps}^2}
        + 9 \frac{\tfrac{\d}{\d \eps} j_{\eps}^i}{\abs{j_{\eps}}}
              \left( \frac{ \tfrac{\d}{\d \eps} (\tfrac{1}{2} \abs{j_\eps}^2) }{\abs{j_\eps}^2} \right)^2
        - 15 \frac{j_{\eps}^i}{\abs{j_{\eps}}}
              \left( \frac{ \tfrac{\d}{\d \eps} (\tfrac{1}{2} \abs{j_\eps}^2) }{\abs{j_\eps}^2} \right)^3
  \\\no
    & = P_{\Omega_{f^\eps}^\bot}^{ik} \frac{ \tfrac{\d^3}{\d \eps^3} j_\eps^k }{\abs{j_\eps}}
        + \Q_3 \left[
              \frac{j_{\eps}^i}{\abs{j_{\eps}}}
                \left( \frac{\tfrac{\d}{\d \eps} j_\eps^k}{\abs{j_\eps}} \cdot \frac{\tfrac{\d^2}{\d \eps^2} j_\eps^k}{\abs{j_\eps}} \right),\,
              \frac{j_{\eps}^i}{\abs{j_{\eps}}}
                \left( \tfrac{ \tfrac{\d}{\d \eps} (\tfrac{1}{2} \abs{j_\eps}^2) }{\abs{j_\eps}^2} \right)^3,\,
              \frac{ \tfrac{\d}{\d \eps} j_\eps^i }{\abs{j_\eps}}
                \frac{ \tfrac{\d^2}{\d \eps^2} (\tfrac{1}{2} \abs{j_\eps}^2) }{\abs{j_\eps}^2},\, \cdots,\,
              \right]
  \\ \label{eq:remainder-3R-homo}
    & = P_{\Omega_{f^\eps}^\bot}^{ik} \frac{ \tfrac{\d^3}{\d \eps^3} j_\eps^k }{\abs{j_\eps}}
        + \Q_3 \left[
                \frac{j_{\eps}^i}{\abs{j_{\eps}}}
                  \left( \frac{\tfrac{\d}{\d \eps} j_\eps^k}{\abs{j_\eps}} \cdot \frac{\tfrac{\d^2}{\d \eps^2} j_\eps^k}{\abs{j_\eps}} \right) ,\,
                \frac{ \tfrac{\d^{\tau_1}}{\d \eps^{\tau_1}} j_\eps^i }{\abs{j_\eps}}
                  \left( \frac{ \tfrac{\d}{\d \eps} (\tfrac{1}{2} \abs{j_\eps}^2) }{\abs{j_\eps}^2} \right)^{\tau_2}
                  \left( \frac{ \tfrac{\d^2}{\d \eps^2} (\tfrac{1}{2} \abs{j_\eps}^2) }{\abs{j_\eps}^2} \right)^{\tau_3} \right] \,,
  \end{align}
where the indices satisfy constraints:
  \begin{align}
  \begin{cases}
    \tau_1 + \tau_2 + 2\tau_3 = 3, \\
    \tau_1 \in \left\{ 0,1,2 \right\}, \quad \tau_2 \in \left\{ 0,1,2,3 \right\},
    \quad \tau_3 \in \left\{ 0,1 \right\}.
   \end{cases}
  \end{align}

The truncated expansion ansatz \cref{eq:expansion-truncate} yields $j_\eps = j_0 + \eps j_1 + \eps^2 (j_2 + j_R)$, then the first term in above expression equals to zero, i.e., $P_{\Omega_{f^\eps}^\bot}^{ik} \frac{ \tfrac{\d^3}{\d \eps^3} j_\eps^k }{\abs{j_\eps}}=0$. The leading order term in $\frac{\d^3}{\d \eps^3} \Omega_{f^\eps}$ is:
  \begin{align}
    % \frac{\d^3}{\d \eps^3} \Omega_{f^\eps}^i =
    -\tfrac{6}{\abs{j_\eps}^2} \Omega_0^i j_1 P_{\Omega_0^\bot} (j_2+ j_R)
    -\tfrac{6}{\abs{j_\eps}^2} (\Omega_0 j_1) [P_{\Omega_0^\bot} (j_2+ j_R)]^i
    -\tfrac{6}{\abs{j_\eps}^2} [\Omega_0 (j_2+ j_R)] (P_{\Omega_0^\bot} j_1)^i,
  \end{align}
which immediately yields, by recalling \cref{eq:RemD3},
  \begin{align}
    \widehat{\Omega}_{3,R}^i = -\tfrac{1}{\abs{j_\eps}^2} \Omega_0^i j_1 P_{\Omega_0^\bot} (j_2+ j_R) - \tfrac{1}{\abs{j_\eps}^2} (\Omega_0 j_1) [P_{\Omega_0^\bot} (j_2+ j_R)]^i - \tfrac{1}{\abs{j_\eps}^2} [\Omega_0 (j_2+ j_R)] (P_{\Omega_0^\bot} j_1)^i
    + \O(\eps) \,.
    % + \O(\eps (j_2+j_\R)^2 ) \,.
  \end{align}
% Therefore, we can get
%   \begin{align}
%     \nm{\widehat{\Omega}_{3,R}}_\lx \le C \nm{j_2 + j_\R}_\lx \le C_{f_0,f_1,f_2} [1 + \nm{j_\R}_\lx + \O (\eps \nm{j_\R^2}_\lx)],
%   \end{align}
% and similarly, get a control on the higher-order derivatives:
%   \begin{align}\label{esm:remainder-3R-j}
%     \nm{\widehat{\Omega}_{3,R}}_{\hx{s}} \le C \nm{j_2 + j_\R}_{\hx{s}} \le C_{f_0,f_1,f_2} [1 + \nm{j_\R}_{\hx{s}} + \O (\eps \nm{j_\R^2}_{\hx{s}})].
%   \end{align}

In order to get an $H^s_x$-estimate, we control the first factor in $\Q_3$, by rewriting
  \begin{align}
    \frac{j_{\eps}^i}{\abs{j_{\eps}}}
        \left( \frac{\tfrac{\d}{\d \eps} j_\eps^k}{\abs{j_\eps}} \cdot \frac{\tfrac{\d^2}{\d \eps^2} j_\eps^k}{\abs{j_\eps}} \right)
    & = \left[ \tfrac{j_0^i}{\abs{j_0}^3} + \left( \tfrac{j_{\eps}^i}{\abs{j_{\eps}}^3} - \tfrac{j_0^i}{\abs{j_0}^3} \right) \right]
      \cdot  \left( {\tfrac{\d}{\d \eps} j_\eps^k} \cdot {\tfrac{\d^2}{\d \eps^2} j_\eps^k} \right) \\\no
    & = 2 \tfrac{j_0^i}{\abs{j_0}^3} [j_1 + 2 \eps (j_2+j_\R)] (j_2+ j_R)
          + \left( \tfrac{j_{\eps}^i}{\abs{j_{\eps}}^3} - \tfrac{j_0^i}{\abs{j_0}^3} \right) [j_1 + 2 \eps (j_2+j_\R)] (j_2+ j_R).
  \end{align}
Notice again the lemma of composition action of smooth function:
  \begin{align}
    \nm{ \tfrac{j_{\eps}^i}{\abs{j_{\eps}}^3} - \tfrac{j_0^i}{\abs{j_0}^3} }_{H^s_x}
    \le C \nm{j_\eps}_{H^s_x}
    \le C_{f_0,f_1,f_2} (1 + \eps^2 \nm{j_2 + j_\R}_{\hx{s}}),
  \end{align}
with the coefficient $C$ depends on the bound of $\nm{j_\eps}_{L^\infty_x}$, then it follows from the Sobolev embedding inequality $H^2 \subset L^\infty$, that
  \begin{align}\label{esm:Q3-factor1}
    \nm{ \frac{j_{\eps}^i}{\abs{j_{\eps}}}
        \left( \frac{\tfrac{\d}{\d \eps} j_\eps^k}{\abs{j_\eps}} \cdot \frac{\tfrac{\d^2}{\d \eps^2} j_\eps^k}{\abs{j_\eps}} \right) }_{H^s_x}
    & \le C_{f_0,f_1} \nm{ [1 + \eps (j_2+j_\R)] (j_2+ j_\R) }_{H^s_x} \\\no
          & \quad + C_{f_0,f_1,f_2} (1 + \eps^2 \nm{j_2 + j_\R}_{\hx{s}}) \nm{ [1 + \eps (j_2+j_\R)] (j_2+ j_\R) }_{H^2_x} \\[5pt] \no
          & \quad + C_{f_0,f_1,f_2} (1 + \eps^2 \nm{j_2 + j_\R}_{\hx{2}}) \nm{ [1 + \eps (j_2+j_\R)] (j_2+ j_\R) }_{H^s_x} \\[5pt] \no
    & \le C_{f_0,f_1,f_2} (1 + \nm{j_\R}_{\hx{s}} + \eps \nm{j_\R^2}_{\hx{s}}) \\\no & \quad
          + C_{f_0,f_1,f_2} (1+ \eps^2 \nm{j_\R}_{\hx{s}}) (1 + \nm{j_\R}_{\hx{2}} + \eps \nm{j_\R^2}_{\hx{2}}) \\\no & \quad
          + C_{f_0,f_1,f_2} (1+ \eps^2 \nm{j_\R}_{\hx{2}}) (1 + \nm{j_\R}_{\hx{s}} + \eps \nm{j_\R^2}_{\hx{s}}) \\\no
    & \le C_{f_0,f_1,f_2} (1 + \nm{j_\R}_{\hx{s}} + \eps \nm{j_\R}_{\hx{s}} \nm{j_\R}_{\hx{2}} + \eps^3 \nm{j_\R}_{\hx{s}} \nm{j_\R}^2_{\hx{2}}) \,,
  \end{align}
where we have used the product estimates \cref{eq:Sobolev-product}, stated in \Cref{lemm:Sobolev-product}.

Next, we rewrite the second factor in $\Q_3$ as, by denoting $\tau = 1+ 2\tau_2 + 2 \tau_3 =3,\,5,\,7$,
  \begin{align}\label{eq:Q3-factor2}
    & \frac{ \tfrac{\d^{\tau_1}}{\d \eps^{\tau_1}} j_\eps^i }{\abs{j_\eps}}
      \cdot \left( \frac{ \tfrac{\d}{\d \eps} (\tfrac{1}{2} \abs{j_\eps}^2) }{\abs{j_\eps}^2} \right)^{\tau_2}
      \cdot \left( \frac{ \tfrac{\d^2}{\d \eps^2} (\tfrac{1}{2}\abs{j_\eps}^2) }{\abs{j_\eps}^2} \right)^{\tau_3}
      \\\no
    & = \left[ \tfrac{1}{\abs{j_0}^\tau} + \left( \tfrac{1}{\abs{j_{\eps}}^\tau} - \tfrac{1}{\abs{j_0}^\tau} \right) \right]
        \tfrac{\d^{\tau_1}}{\d \eps^{\tau_1}} j_\eps^i
          \cdot \left[ \tfrac{\d}{\d \eps} (\tfrac{1}{2} \abs{j_\eps}^2) \right]^{\tau_2}
          \cdot \left[ \tfrac{\d^2}{\d \eps^2} (\tfrac{1}{2}\abs{j_\eps}^2) \right]^{\tau_3}  \\\no
    & = \left[ \tfrac{1}{\abs{j_0}^\tau} + \left( \tfrac{1}{\abs{j_{\eps}}^\tau} - \tfrac{1}{\abs{j_0}^\tau} \right) \right]
      \cdot (1+ \eps^{2-\tau_1} j_\R) \cdot [1+ (\eps j_\R)^{\tau_2}] \cdot (1 + j_\R^{\tau_3}) + \text{Higher Orders}\,.
  \end{align}
The worst terms are those with smaller powers of $\eps$, namely, with bigger $\tau_1,\,\tau_3$ and a smaller $\tau_2$. So we focus on the cases $(\tau_1,\tau_2,\tau_3) =(2,1,0) $ and $(1,0,1)$, both of which result in
  \begin{align}
    & (1+ \eps^{2-\tau_1} j_\R) \cdot [1+ (\eps j_\R)^{\tau_2}] \cdot (1 + j_\R^{\tau_3}) \\\no
    & = 1 + \eps^{2-\tau_1} j_\R + \eps^{\tau_2} j_\R^{\tau_2} + j_\R^{\tau_3}
        + \eps^{2-\tau_1 + \tau_2} j_\R^{1+\tau_2} + \eps^{2-\tau_1} j_\R^{1+\tau_3} + \eps^{\tau_2} j_\R^{\tau_2 + \tau_3}
        + \eps^{2-\tau_1 + \tau_2} j_\R^{1+\tau_2+\tau_3} \\\no
    & \sim 1 + j_\R % + \eps j_\R
             + \eps j_\R^2 \,,
  \end{align}
with the worst term $\eps j_\R^2$ arising from $\eps^{2-\tau_1 + \tau_2} j_\R^{1+\tau_2+\tau_3}$. Then \cref{eq:Q3-factor2} results in, following the same argument as that in obtaining \cref{esm:Q3-factor1},
  \begin{align}\label{esm:Q3-factor2}
    & \nm{ \frac{ \tfrac{\d^{\tau_1}}{\d \eps^{\tau_1}} j_\eps^i }{\abs{j_\eps}}
      \cdot \left( \frac{ \tfrac{\d}{\d \eps} (\tfrac{1}{2} \abs{j_\eps}^2) }{\abs{j_\eps}^2} \right)^{\tau_2}
      \cdot \left( \frac{ \tfrac{\d^2}{\d \eps^2} (\tfrac{1}{2}\abs{j_\eps}^2) }{\abs{j_\eps}^2} \right)^{\tau_3} }_{H^s_x} \\\no
    & \le C_{f_0,f_1,f_2} (1 + \nm{j_\eps}_{H^2_x}) (1 + \nm{j_\R}_{\hx{s}} + \eps \nm{j_\R^2}_{\hx{s}})
          + \nm{j_\eps}_{H^2_x} (1 + \nm{j_\R}_{\hx{2}} + \eps \nm{j_\R^2}_{\hx{2}}) \\ \no
    & \le C_{f_0,f_1,f_2} (1 + \nm{j_\R}_{\hx{s}} + \eps \nm{j_\R}_{\hx{s}} \nm{j_\R}_{\hx{2}} + \eps^3 \nm{j_\R}_{\hx{s}} \nm{j_\R}^2_{\hx{2}}) \,.
  \end{align}

Therefore, combining together \cref{esm:Q3-factor1} and \cref{esm:Q3-factor2} enables us to get a control on the higher-order derivatives:
  \begin{align}\label{esm:remainder-3R-j}
    \nm{\widehat{\Omega}_{3,R}}_{\hx{s}}
      \le C_{f_0,f_1,f_2} (1 + \nm{j_\R}_{\hx{s}} + \eps \nm{j_\R}_{\hx{s}} \nm{j_\R}_{\hx{2}} + \eps^3 \nm{j_\R}_{\hx{s}} \nm{j_\R}^2_{\hx{2}}) \,.
  \end{align}
Hence, the desired result \cref{esm:Taylor-remainders-3R} in \Cref{lemm:Taylor-remainders} follows.

In addition, a similar scheme as above yields:
  \begin{align}\label{esm:remainder-3R-0}
    \nm{\widehat{\Omega}_{3,R}}_\lx
    & \le C_{f_0,f_1,f_2} (1 + \nm{j_\R}_\lx + \eps^2 \nm{j_\R}_{\hx{2}} + \eps \nm{j_\R}_\lx \nm{j_\R}_{\hx{2}} + \eps^3 \nm{j_\R}_\lx \nm{j_\R}^2_{\hx{2}}) \\\no
    & \le C_{f_0,f_1,f_2} \Big[ 1+ \nm{f_\R}_{L^2_{x,v}(M^{-1})} + \eps^2 \nm{f_\R}_{\H^2_x L^2_v(M^{-1})} + \eps \nm{f_\R}_{L^2_{x,v}(M^{-1})} \nm{f_\R}_{\H^2_x L^2_v(M^{-1})} \\\no & \qquad \qquad \quad + \eps^3 \nm{f_\R}_{L^2_{x,v}(M^{-1})} \nm{f_\R}^2_{\H^2_x L^2_v(M^{-1})} \Big] \,.
  \end{align}

As a comment here, note that the above two cases corresponds to the third and fourth terms in the expression \cref{eq:remainder-3R}, which can be recast as
  \begin{align}
    & \frac{ \tfrac{\d^2}{\d \eps^2} j_\eps^i }{\abs{j_\eps}}
      \cdot \frac{ \tfrac{\d}{\d \eps} (\tfrac{1}{2} \abs{j_\eps}^2) }{\abs{j_\eps}^2}
      + \frac{ \tfrac{\d}{\d \eps} j_\eps^i }{\abs{j_\eps}}
          \cdot \frac{ \tfrac{\d^2}{\d \eps^2} (\tfrac{1}{2} \abs{j_\eps}^2) }{\abs{j_\eps}^2} \\\no
    & = \left[ \tfrac{1}{\abs{j_0}^3}
              + \left( \tfrac{1}{\abs{j_{\eps}}^3} - \tfrac{1}{\abs{j_0}^3} \right)
      \right]
      \cdot
        \left[ \tfrac{\d^2}{\d \eps^2} j_\eps^i \cdot \tfrac{\d}{\d \eps} (\tfrac{1}{2} \abs{j_\eps}^2)
                + \tfrac{\d}{\d \eps} j_\eps^i \cdot \tfrac{\d^2}{\d \eps^2} (\tfrac{1}{2} \abs{j_\eps}^2)
        \right] \,.
  \end{align}
This formula also leads us to the same control as above. However, we hope to emphasize that it provides a unified way to estimate the higher-order expansions by analyzing the homogeneous function (like $\Q_3$).
\endproof

\subsection{A priori estimate uniformly-in-$\eps$} % (fold)
\label{sub:uniform-esm}

% subsection uniform-estimate (end)

\subsubsection{Basic estimate} % (fold)
\label{ssub:L2-esm}

% subsubsection main_result (end)
For simplicity, we omit the supscript $\eps$ in the remainder equation \cref{eq:Rem-order} in this section, i.e.:
  \begin{multline}\label{eq:Rem-order-re}
    (\p_t + v \cdot \nabla_x) f_\R - \frac{1}{\eps} \L_{\Omega_0} f_\R
    = - \div_v \left( f_0 \widehat{\Omega}_{3,R} \right)
      - \div_v \left[ (f_1 + \eps f_2) \widehat{\Omega}_\R^\eps \right]
    % \\\no
      - \div_v \left[ f_\R \tfrac{\left( \Omega_{f^\eps} - \Omega_0 \right)}{\eps} \right]
        - \S(f_1,f_2) \,.
  \end{multline}

Multiply \cref{eq:Rem-order-re} by $\tfrac{f_\R}{M_0}$, and take integration over $v$ and $x$, then we get, for the first term,
  \begin{align}
    \sskp{\p_t f_\R}{\tfrac{f_\R}{M_0}} & = \frac{1}{2} \frac{\d}{\d t} \iint \frac{\abs{f_\R}^2}{M_0} \d v \d x - \iint {\tfrac{1}{2} \p_t M_0^{-1} \cdot \abs{f_\R}^2} \d v \d x \\\no
    & = \frac{1}{2} \frac{\d}{\d t} \iint \frac{\abs{f_\R}^2}{M_0} \d v \d x + \iint \tfrac{1}{2} v \cdot \p_t \Omega_0 \frac{\abs{f_\R}^2}{M_0} \d v \d x\,,  % \\\no & =
  \end{align}
and for the second term,
  \begin{align}
    \sskp{v \cdot \nabla_x f_\R}{\tfrac{f_\R}{M_0}} = \agll{v \cdot \tfrac{1}{2} \nabla_x (\tfrac{\abs{f_\R}^2}{M_0}) }
    & = \tfrac{1}{2} \iint \nabla_x \cdot (v \tfrac{\abs{f_\R}^2}{M_0}) \d v \d x - \iint {\tfrac{1}{2} v \cdot \nabla_x M_0^{-1} \abs{f_\R}^2} \d v \d x \\\no
    & = 0 + \iint \tfrac{1}{2} v \otimes v : \nabla_x \Omega_0 \frac{\abs{f_\R}^2}{M_0} \d v \d x\,.
  \end{align}

Recalling the fact $|v|^2 \le |\nabla_v \Phi|^2$, the same argument as that for the Poincar\'e inequality \eqref{eq:Poinc-ineq-Phi} in \Cref{lemm:Poincare-ineq} yields,
  \begin{align}\label{eq:vv-poincare}
    & \abs{\iint v \otimes v : \nabla_x \Omega_0 \frac{\abs{f_\R}^2}{M_0} \d v \d x} \\\no
    & \le \nm{\nabla_x \Omega_0}_{L^\infty_x} {\iint |\nabla_v \Phi|^2 \abs{ \frac{f_\R}{M_0} }^2 M_0 \d v \d x}  \\\no
    & \le C_{\beta} \nm{\nabla_x \Omega_0}_{L^\infty_x} \left[ \iint \abs{\nabla_v \left( \frac{f_\R}{M_0} \right)}^2 M_0 \d v \d x + \iint \abs{\frac{f_\R}{M_0}}^2 M_0 \d v \d x \right] \,,
  \end{align}
where an additional integral involving $\tfrac{f_\R}{M_0}$ itself arose here, due to its nonzero mean value $\agl{\tfrac{f_\R}{M_0}}_\M = \rho_\R \neq 0$, comparing \cref{eq:potential-poinc}.

Similarly, we have
  \begin{align}
    & \abs{\iint v \p_t \Omega_0 \frac{\abs{f_\R}^2}{M_0} \d v \d x} \\\no
    & \le \nm{\p_t \Omega_0}_{L^\infty_x}
            \left\{ \iint \frac{\abs{f_\R}^2}{M_0} \d v \d x \right\}^{\frac{1}{2}}
      \cdot \left\{ \iint |\nabla_v \Phi|^2 \abs{ \frac{f_\R}{M_0} }^2 M_0 \d v \d x \right\}^{\frac{1}{2}} \\\no
    & \le C_{\beta} \nm{\p_t \Omega_0}_{L^\infty_x}
            \left\{ \iint \frac{\abs{f_\R}^2}{M_0} \d v \d x \right\}^{\frac{1}{2}}
      \cdot \left\{ \iint \abs{\nabla_v \left( \frac{f_\R}{M_0} \right)}^2 M_0 \d v \d x + \iint \frac{\abs{f_\R}^2}{M_0} \d v \d x \right\}^{\frac{1}{2}} \,.
  \end{align}

  % \begin{align}
  %   \abs{\iint v \p_t \Omega_0 \frac{\abs{f_\R}^2}{M_0} \d v \d x}
  %   & \le \nm{\p_t \Omega_0}_{L^\infty_x}
  %           \left\{ \iint \frac{\abs{f_\R}^2}{M_0} \d v \d x \right\}^{\frac{1}{2}}
  %     \cdot \left\{ \iint |\nabla_v \Phi|^2 \abs{ \frac{f_\R}{M_0} }^2 M_0 \d v \d x \right\}^{\frac{1}{2}}   \\\no
  %   & \le C_{\beta} \nm{\p_t \Omega_0}_{L^\infty_x}
  %           \left\{ \iint \frac{\abs{f_\R}^2}{M_0} \d v \d x \right\}^{\frac{1}{2}}
  %     \cdot \left\{ \iint \abs{\nabla_v \left( \frac{f_\R}{M_0} \right)}^2 M_0 \d v \d x + \iint \frac{\abs{f_\R}^2}{M_0} \d v \d x \right\}^{\frac{1}{2}} \,.
  % \end{align}

Combining together with the coercivity estimate on linear operator $\L_{\Omega_0}$ in \Cref{prop:coercivity}, the left-hand side of \eqref{eq:Rem-order-re} can be controlled from below, as follows:
  \begin{align}
    & \sskp{(\p_t + v \cdot \nabla_x) f_\R}{\tfrac{f_\R}{M_0}}
      - \frac{1}{\eps} \sskp{ \L_{\Omega_0} f_\R }{\tfrac{f_\R}{M_0}}
  \\\no
    & \ge \frac{1}{2} \frac{\d}{\d t} \iint \frac{\abs{f_\R}^2}{M_0} \d v \d x
        + \frac{3}{4\eps} \iint \abs{\nabla_v (\tfrac{f_\R}{M_0})}^2 M_0 \d v \d x
        + \frac{3\alpha}{4c_0 \eps} \int \abs{P_{\Omega_0^\bot} {j_\R}}^2 \d x\\\no
    & \quad - C_{\beta} \nm{\nabla_x \Omega_0}_{L^\infty_x} \left[ \iint \abs{\nabla_v \left( \frac{f_\R}{M_0} \right)}^2 M_0 \d v \d x + \iint \abs{\frac{f_\R}{M_0}}^2 M_0 \d v \d x \right] \\\no
    & \quad - C_{\beta} \nm{\p_t \Omega_0}_{L^\infty_x}
            \left\{ \iint \frac{\abs{f_\R}^2}{M_0} \d v \d x \right\}^{\frac{1}{2}}
      \cdot \left\{ \iint \abs{\nabla_v (\tfrac{f_\R}{M_0})}^2 M_0 \d v \d x + \iint \frac{\abs{f_\R}^2}{M_0} \d v \d x \right\}^{\frac{1}{2}}
  \\[5pt] \no
    & \ge \frac{1}{2} \frac{\d}{\d t} \nm{f_\R}^2_{\MV}
        + \frac{3}{4\eps} \nm{\nabla_v (\tfrac{f_\R}{M_0})}^2_{\M}
        + \frac{3\alpha}{4c_0 \eps} \nm{P_{\Omega_0^\bot} {j_\R}}_\lx^2 \\\no
    & \qquad - C_{\beta,\Omega_0}
      \left[ \nm{\nabla_v (\tfrac{f_\R}{M_0})}^2_{\M} + \nm{f_\R}_{\MV} \nm{\nabla_v (\tfrac{f_\R}{M_0})}_{\M} + \nm{f_\R}_{\MV}^2 \right] \,,
  \end{align}
where the constant $C_{\beta,\Omega_0}$ depends on $\nm{(\p_t \Omega_0, \nabla_x \Omega_0)}_{L^\infty_x}\ls \nm{(\rho_0,\nabla_x \Omega_0)}_\hx{3}$, the coefficient $\beta$ and other irrelevant constants.

As for the first and second terms on the right-hand side of \eqref{eq:Rem-order-re}, we notice they are linear dependent on the remainder function $f_\R^\eps$. We can infer that:
% the following lemmas:
  \begin{align}
    - & \sskp{\div_v \left( f_0 \widehat{\Omega}_{3,R} \right)}{\tfrac{f_\R}{M_0}} \no\\
    & = \sskp{ \rho_0 M_0 \widehat{\Omega}_{3,R} }{ \nabla_v (\tfrac{f_\R}{M_0}) } \no\\
    & \le C_{\rho_0} \nm{\widehat{\Omega}_{3,R}}_\lx \nm{\nabla_v (\tfrac{f_\R}{M_0})}_{\M} \\\no
    % & \le C_{f_0,f_1,f_2} \nm{\nabla_v (\tfrac{f_\R}{M_0})}_{\M} \Big[ 1+ \nm{f_\R}_{\MV} + \eps^2 \nm{f_\R}_{\H^2_x L^2_v(M^{-1})} \\\no & \hspace*{4.5cm}
    %       + \eps \nm{f_\R}_{\MV} \nm{f_\R}_{\H^2_x L^2_v(M^{-1})} \\\no & \hspace*{4.5cm} + \eps^3 \nm{f_\R}_{\MV} \nm{f_\R}^2_{\H^2_x L^2_v(M^{-1})} \Big] \,,  \\\no
    & \le C_{f_0,f_1,f_2} \nm{\nabla_v (\tfrac{f_\R}{M_0})}_{\M} \Big[ 1+ \nm{f_\R}_{\MV} + \eps^2 \nm{f_\R}_{\H^2_x L^2_v(M^{-1})} \\\no & \hspace*{4cm}
          + \eps \nm{f_\R}_{\MV} \nm{f_\R}_{\H^2_x L^2_v(M^{-1})} + \eps^3 \nm{f_\R}_{\MV} \nm{f_\R}^2_{\H^2_x L^2_v(M^{-1})} \Big] \,,
  \end{align}
and
  \begin{align}
    - & \sskp{ \div_v \left[ (f_1 + \eps f_2) \widehat{\Omega}_\R^\eps \right] }{\tfrac{f_\R}{M_0}} \no\\
    & = \sskp{ (f_1 + \eps f_2) \left[ (c_0\rho_0)^{-1} P_{\Omega_0^\bot} j_\R + \eps \widehat{\Omega}_{3,R} \right] }{ \nabla_v (\tfrac{f_\R}{M_0}) } \\\no
    & \le C_{\rho_0} \nm{f_1 + \eps f_2}_{L^\infty_x L^2_v(M^{-1})} \left[ \nm{P_{\Omega_0^\bot} j_\R}_\lx + \eps \nm{\widehat{\Omega}_{3,R}}_\lx \right] \nm{\nabla_v (\tfrac{f_\R}{M_0})}_{\M} \\\no
    & \le C_{\rho_0,f_1,f_2} \nm{P_{\Omega_0^\bot} j_\R}_\lx \nm{\nabla_v (\tfrac{f_\R}{M_0})}_{\M} + C_{f_0,f_1,f_2} \eps \nm{\widehat{\Omega}_{3,R}}_\lx \nm{\nabla_v (\tfrac{f_\R}{M_0})}_{\M} \,.
  \end{align}
Here the constant $C_{\rho_0,f_1,f_2}$ depends on $\nm{\rho_0}_{L^\infty_x}$ and $\nm{(f_1,f_2)}_{\hxv{2}(M^{-1})}$, both of which are bounded eventually by $\nm{(\rho_0,\Omega_0)}_\hx{12}$ and $\nm{(\rho_1,\Omega_1)}_\hx{4}$.

We turn to consider the nonlinear term $\div_v [ f_\R \tfrac{( \Omega_{f^\eps} - \Omega_0 )}{\eps} ]$, in the third term of the right-hand side of \cref{eq:Rem-order-re}:
  \begin{align}
    - & \sskp{ \div_v \left[ f_\R \cdot \tfrac{( \Omega_{f^\eps} - \Omega_0 )}{\eps} \right] }{\tfrac{f_\R}{M_0}} \\\no
    & = \sskp{ f_\R \tfrac{\left( \Omega_{f^\eps} - \Omega_0 \right)}{\eps} }{ \nabla_v (\tfrac{f_\R}{M_0}) } \\\no
    & \le \abs{\tfrac{( \Omega_{f^\eps} - \Omega_0 )}{\eps}}_{L^\infty_x} \nm{f_\R}_{\MV} \nm{\nabla_v (\tfrac{f_\R}{M_0})}_{\M} \\\no
    & \le \nm{\tfrac{( \Omega_{f^\eps} - \Omega_0 )}{\eps}}_{H^2_x} \nm{f_\R}_{\MV} \nm{\nabla_v (\tfrac{f_\R}{M_0})}_{\M} \\\no
    & \le C_{f_0,f_1,f_2} \left( 1 + \eps \nm{f_\R}_{\H^2_x L^2_v(M^{-1})} + \eps^3 \nm{f_\R}^2_{\H^2_x L^2_v(M^{-1})} \right) \nm{f_\R}_{\MV} \nm{\nabla_v (\tfrac{f_\R}{M_0})}_{\M} \,.
  \end{align}

The last term $\S(f_1,f_2)$ is a known source term, so we have
  \begin{align}
    \sskp{ \S(f_1,f_2) }{\tfrac{f_\R}{M_0}} \le C(\S) \nm{f_\R}_{\MV}\,.
  \end{align}

Therefore, we can combine all above estimate to derive
  \begin{align}
    \frac{1}{2} & \frac{\d}{\d t} \nm{f_\R}^2_{\MV}
        + \frac{3}{4\eps} \nm{\nabla_v (\tfrac{f_\R}{M_0})}^2_{\M}
        + \frac{3\alpha}{4c_0 \eps} \nm{P_{\Omega_0^\bot} {j_\R}}_\lx^2 \\\no
    & \le C_{\beta,\Omega_0}
      \left[ \nm{\nabla_v (\tfrac{f_\R}{M_0})}^2_{\M} + \nm{f_\R}_{\MV} \nm{\nabla_v (\tfrac{f_\R}{M_0})}_{\M} + \nm{f_\R}_{\MV}^2 \right] % \\\no & \quad
        + C_{\rho_0,f_1,f_2} \nm{P_{\Omega_0^\bot} j_\R}_\lx \nm{\nabla_v (\tfrac{f_\R}{M_0})}_{\M} \\\no & \quad
        + C_{f_0,f_1,f_2} (1+ \eps) \nm{\nabla_v (\tfrac{f_\R}{M_0})}_{\M} \Big[ 1+ \nm{f_\R}_{\MV} + \eps^2 \nm{f_\R}_{\H^2_x L^2_v(M^{-1})} + \eps \nm{f_\R}_{\MV} \nm{f_\R}_{\H^2_x L^2_v(M^{-1})} \\\no & \hspace*{5.8cm}
          + \eps^3 \nm{f_\R}_{\MV} \nm{f_\R}^2_{\H^2_x L^2_v(M^{-1})} \Big]
      % \\\no & \quad
      % + C_{\rho_0,f_1,f_2} \left( \nm{P_{\Omega_0^\bot} j_\R}_\lx + \eps \nm{f_\R}_{\MV} \right) \nm{\nabla_v (\tfrac{f_\R}{M_0})}_{\M} \\\no &
      % \quad + (C_{f_0,f_1,f_2} + C \eps^2 \nm{f_\R}_{\hxv{2}(\MV)}) \nm{f_\R}_{\MV} \nm{\nabla_v (\tfrac{f_\R}{M_0})}_{\M} %
      % \\\no & \quad
      % + C_{f_0,f_1,f_2} \left( 1 + \eps \nm{f_\R}_{\H^2_x L^2_v(M^{-1})} + \eps^3 \nm{f_\R}^2_{\H^2_x L^2_v(M^{-1})} \right) \nm{f_\R}_{\MV} \nm{\nabla_v (\tfrac{f_\R}{M_0})}_{\M}
      + C(\S) \nm{f_\R}_{\MV} \\\no
    & \le C_{\beta,\Omega_0} \nm{\nabla_v (\tfrac{f_\R}{M_0})}^2_{\M}
        + C_{\rho_0,f_1,f_2} \nm{P_{\Omega_0^\bot} j_\R}_\lx \nm{\nabla_v (\tfrac{f_\R}{M_0})}_{\M}
        + C_{f_0,f_1,f_2} (\nm{f_\R}_{\MV}^2 + \nm{f_\R}_{\MV})
        \\\no & \quad
        + C_{f_0,f_1,f_2} \nm{\nabla_v (\tfrac{f_\R}{M_0})}_{\M} \Big[ 1+ \nm{f_\R}_{\MV} + \eps^2 \nm{f_\R}_{\H^2_x L^2_v(M^{-1})} + \eps \nm{f_\R}_{\MV} \nm{f_\R}_{\H^2_x L^2_v(M^{-1})} \\\no & \hspace*{4.8cm}
          + \eps^3 \nm{f_\R}_{\MV} \nm{f_\R}^2_{\H^2_x L^2_v(M^{-1})} \Big] \,.
  \end{align}

Corresponding to the definition of energy and dissipation functionals \cref{eq:functionals-ED}, we define, for each index $k \ge 0$:
  \begin{align}\label{eq:functionals-ED-order-eps}
    E_{k,\eps}(t) = \eps^k \nm{\tfrac{\nabla_x^k f_\R^\eps}{M_0}}_{\lxv(\M)}^2, \
    D_{k,\eps}(t) = \eps^{k-1} \left( \nm{\nabla_v (\tfrac{\nabla_x^k f_\R^\eps}{M_0})}_{\lxv(\M)}^2 + \frac{\alpha}{c_0} \nm{P_{\Omega_0^\bot} {\nabla_x^k j_\R^\eps}}_{L^2_x}^2 \right) .
  \end{align}
For brevity, we denote $E_{0,\eps}(t) = E_{0}(t)$ when $k=0$. Then we have
  \begin{multline}
    \tfrac{1}{2} \tfrac{\d}{\d t} E_0 + \tfrac{3}{4} D_{0,\eps}
    \le \eps (C_{\beta,\Omega_0} + C_{\rho_0,f_1,f_2}) D_{0,\eps} + C_{f_0,f_1,f_2} (E_0^{\frac{1}{2}} + E_0) \\
          + \eps^{\frac{1}{2}} C_{f_0,f_1,f_2} D_{0,\eps}^{\frac{1}{2}}
          \left( 1 + E_0^{\frac{1}{2}} + \eps E_{2,\eps}^{\frac{1}{2}} + E_0^{\frac{1}{2}} E_{2,\eps}^{\frac{1}{2}} + \eps E_0^{\frac{1}{2}} E_{2,\eps} \right) \,.
  \end{multline}
It follows from the H\"older inequality that, for sufficiently small $\eps$ (such that, $\eps (C_{\beta,\Omega_0} + C_{\rho_0,f_1,f_2}) + \eps^{\frac{1}{2}} C_{f_0,f_1,f_2} < 1/4 $),
  \begin{align}\label{esm:ZeroOrder}
    \tfrac{\d}{\d t} E_0 + D_{0,\eps}
    & \le C_{f_0,f_1,f_2} (E_0^{\frac{1}{2}} + E_0)
        + \eps^{\frac{1}{2}} C_{f_0,f_1,f_2}
          \left( 1 + E_0 + \eps^2 E_{2,\eps} + E_0 E_{2,\eps} + \eps^2 E_0 E^2_{2,\eps} \right) \\\no
    & \le C_{f_0,f_1,f_2} (1 + E_0) + \eps^{\frac{1}{2}} C_{f_0,f_1,f_2} (E_0 E_{2,\eps} + \eps^2 E_{2,\eps} + \eps^2 E_0 E^2_{2,\eps}) \,.
  \end{align}

%% ------------------------------------------------------------------------ %%
\subsubsection{Higher-order derivatives estimate} % (fold)
\label{ssub:higher-order-esm}

% subsubsection Higher-order-estimates (end)

Applying higher-order spatial derivative operator $\nabla_x^s$ (with $s\ge 3$) on the remainder equation \eqref{eq:Rem-order-re}, we get
  \begin{multline}\label{eq:Remd-HighDeri}
    (\p_t + v \cdot \nabla_x) (\nabla_x^s f_\R) - \frac{1}{\eps} \nabla_x^s \L_{\Omega_0} f_\R
    = - \div_v \left[ \nabla_x^s ( f_0 \widehat{\Omega}_{3,R}) \right]
      - \div_v \left\{ \nabla_x^s \left [ (f_1 + \eps f_2) \widehat{\Omega}_\R^\eps \right] \right\} \\[3pt]
      - \div_v \left\{ \nabla_x^s \left[ f_\R \tfrac{(\Omega_{f^\eps} - \Omega_0)}{\eps} \right] \right\}
        - \nabla_x^s \S(f_1,f_2) \,.
  \end{multline}
We now consider the inner product with $\tfrac{\nabla_x^s f_\R}{M_0}$ term by term in the following contents.

By the expression of linear operator $\L_{\Omega_0}$, a commutator arises here:
  \begin{align}
    [\nabla_x^s,\, \L_{\Omega_0}] f_\R
    & = \nabla_x^s \L_{\Omega_0} f_\R - \L_{\Omega_0} \nabla_x^s f_\R \\\no
    & = \div_v \left\{ - \sum_{\mathclap{\substack{s_1+s_2 \le s\\s_2 \ge 1}}} \nabla_x^{s_1} f_\R \nabla_x^{s_2} \Omega_0
        - c_0^{-1} \sum_{\mathclap{\substack{s_1+s_2 \le s\\s_2 \ge 1}}} \nabla_x^{s_1} j_\R \nabla_x^{s_2} (P_{\Omega_0^\bot} M_0) \right\}.
   \end{align}
Considering the inner product, we have, for the first term,
  \begin{align}
    \sskp{- \div_v (\nabla_x^{s_1} f_\R \nabla_x^{s_2} \Omega_0)}{\tfrac{\nabla_x^s f_\R}{M_0}}
    & \le \nm{\nabla_x^{s_2} \Omega_0}_{L^\infty_x}
        \nm{\nabla_x^{s_1} f_\R}_{\MV} \nm{\nabla_v (\tfrac{\nabla_x^s f_\R}{M_0})}_{\M}  \\\no
    & \le C_{\Omega_0} \nm{\nabla_x^{s_1} f_\R}_{\MV} \nm{\nabla_v (\tfrac{\nabla_x^s f_\R}{M_0})}_{\M},
  \end{align}
and, for the second term,
  \begin{align}
    \sskp{- \div_v [\nabla_x^{s_1} j_\R \nabla_x^{s_2} (P_{\Omega_0^\bot} M_0)]}{\tfrac{\nabla_x^s f_\R}{M_0}}
    & \le \nm{\nabla_x^{s_1} j_\R}_{\lx}
          \nm{\nabla_x^{s_2} (P_{\Omega_0^\bot} M_0)}_{L^\infty_xL^2_v(\MV)}
        \nm{\nabla_v (\tfrac{\nabla_x^s f_\R}{M_0})}_{\M}  \\\no
    & \le C_{\Omega_0} \nm{\nabla_x^{s_1} f_\R}_{\MV} \nm{\nabla_v (\tfrac{\nabla_x^s f_\R}{M_0})}_{\M},
  \end{align}
where we have used the basic fact:
  \begin{multline*}
    {\nabla_x^{s_2} (P_{\Omega_0^\bot} M_0)}
    \\= \sum_{\mathclap{s_{21}+s_{22} = s_2}} \nabla_x^{s_{21}} P_{\Omega_0^\bot} \left[ v_\bot \nabla_x^{s_{22}} \Omega_0 M_0 + (v_\bot \nabla_x^{s_{22}-1} \Omega_0) \otimes (v_\bot \nabla_x \Omega_0) M_0 + \cdots + (v_\bot \nabla_x \Omega_0)^{\otimes s_{22}}  M_0 \right],
  \end{multline*}
and hence, $\nm{\nabla_x^{s_2} (P_{\Omega_0^\bot} M_0)}_{L^\infty_xL^2_v(\MV)} \le C_{\Omega_0} \agl{M_0} = C_{\Omega_0}$.

As a consequence, the commutator can be controlled by
  \begin{align}
    \sskp{[\nabla_x^s,\, \L_{\Omega_0}] f_\R}{\tfrac{\nabla_x^s f_\R}{M_0}}
     \le C_{\Omega_0} \nm{\nabla_x^{s_1} f_\R}_{\MV} \nm{\nabla_v (\tfrac{\nabla_x^s f_\R}{M_0})}_{\M}.
  \end{align}
Combining this inequality and the coercivity estimate for linear operator:
  \begin{align}
    - \frac{1}{\eps} \sskp{ \L_{\Omega_0} \nabla_x^s f_\R }{\tfrac{\nabla_x^s f_\R}{M_0}}
    \ge \frac{3}{4} \frac{1}{\eps} \nm{\nabla_v (\tfrac{\nabla_x^s f_\R}{M_0})}_{\M}^2
        + \frac{3\alpha}{4c_0} \frac{1}{\eps} \nm{P_{\Omega_0^\bot} \nabla_x^{s} j_\R}_{\lx}^2,
  \end{align}
we can infer that
  \begin{align}\label{eq:coercive-HiOd}
    & \frac{3}{4} \frac{1}{\eps} \nm{\nabla_v (\tfrac{\nabla_x^s f_\R}{M_0})}_{\M}^2
        + \frac{3\alpha}{4c_0} \frac{1}{\eps} \nm{P_{\Omega_0^\bot} \nabla_x^{s} j_\R}_{\lx}^2 \no\\
    & \le - \frac{1}{\eps} \sskp{ \nabla_x^s \L_{\Omega_0} f_\R }{\tfrac{\nabla_x^s f_\R}{M_0}}
          + \frac{1}{\eps} \sskp{ [\nabla_x^s,\, \L_{\Omega_0}] f_\R }{\tfrac{\nabla_x^s f_\R}{M_0}} \no\\
    & \le - \frac{1}{\eps} \sskp{ \nabla_x^s \L_{\Omega_0} f_\R }{\tfrac{\nabla_x^s f_\R}{M_0}}
          + \sum_{\mathclap{\substack{s_1+s_2 \le s\\s_2 \ge 1}}} C_{\Omega_0} \frac{1}{\eps} \nm{\nabla_x^{s_1} f_\R}_{\MV} \nm{\nabla_v (\tfrac{\nabla_x^s f_\R}{M_0})}_{\M} \,.
  \end{align}

We employ the notations, for each index $k \ge 0$:
  \begin{align}\label{eq:functionals-ED-order}
    E_{k}(t) = \nm{\tfrac{\nabla_x^k f_\R^\eps}{M_0}}_{\lxv(\M)}^2\,, \quad
    D_{k}(t) = \nm{\nabla_v (\tfrac{\nabla_x^k f_\R^\eps}{M_0})}_{\lxv(\M)}^2 + \frac{\alpha}{c_0} \nm{P_{\Omega_0^\bot} {\nabla_x^k j_\R^\eps}}_{L^2_x}^2 \,.
  \end{align}
Then we can reformulate \cref{eq:coercive-HiOd} as, by employing a small parameter $\nu$,
  \begin{align}\label{eq:coercive-HiOd-nu}
    \nu^s \tfrac{1}{\eps} D_s
    & \le - \nu^s \tfrac{1}{\eps} \sskp{ \nabla_x^s \L_{\Omega_0} f_\R }{\tfrac{\nabla_x^s f_\R}{M_0}}
          + C \frac{1}{\eps} \sum_{\mathclap{s_1 \le s-1}} \nu^s E_{s_1}^{\frac{1}{2}} D_s^{\frac{1}{2}} \\\no
    & \le - \nu^s \tfrac{1}{\eps} \sskp{ \nabla_x^s \L_{\Omega_0} f_\R }{\tfrac{\nabla_x^s f_\R}{M_0}}
          + C \sum_{\mathclap{s_1 \le s-1}} (\tfrac{\nu^{s-s_1}}{\eps})^{\frac{1}{2}} (\nu^{s_1} E_{s_1})^{\frac{1}{2}} \cdot (\nu^s \tfrac{1}{\eps} D_s)^{\frac{1}{2}} \,,
  \end{align}
which yields, by noticing $s-s_1 \ge 1$ and by taking $\nu = \eps$,
  \begin{align}\label{eq:coercive-HiOd-eps}
    D_{s,\eps}
    \le - \eps^{s-1} \sskp{ \nabla_x^s \L_{\Omega_0} f_\R }{\tfrac{\nabla_x^s f_\R}{M_0}}
          + C \sum_{\mathclap{s_1 \le s-1}} E_{s_1,\eps}^{\frac{1}{2}} \cdot D_{s,\eps}^{\frac{1}{2}} \,.
  \end{align}
As a result, the dissipation effect on the right-hand side could be absorbed by the left-hand side, which provides us a possible way to close the energy estimate. Note the appropriate coefficient $\nu^s = \eps^s$ plays an crucial role, which depends on the small parameter $\eps$, and especially on the Sobolev index $s$. This is exactly the reason we get the convergence result in $\hxv{3}(M^{-1})$, see \Cref{thm:limit}.

Perform the inner product on \cref{eq:Remd-HighDeri} with $\tfrac{\nabla_x^s f_\R}{M_0}$, with respect to both variables of $v$ and $x$, then we get, for the first term on the left-hand side,
  \begin{align}
    \sskp{\p_t \nabla_x^s f_\R}{\tfrac{\nabla_x^s f_\R}{M_0}}
    & = \frac{1}{2} \frac{\d}{\d t} \iint \frac{\abs{\nabla_x^s f_\R}^2}{M_0} \d v \d x - \iint {\tfrac{1}{2} \p_t M_0^{-1} \cdot \abs{\nabla_x^s f_\R}^2} \d v \d x \\\no
    & = \frac{1}{2} \frac{\d}{\d t} \iint \frac{\abs{\nabla_x^s f_\R}^2}{M_0} \d v \d x + \iint \tfrac{1}{2} v \cdot \p_t \Omega_0 \frac{\abs{\nabla_x^s f_\R}^2}{M_0} \d v \d x\,,  % \\\no & =
  \end{align}
and for the second term,
  \begin{align}
    \sskp{v \cdot \nabla_x (\nabla_x^s f_\R)}{\tfrac{\nabla_x^s f_\R}{M_0}}
    % = \agll{v \cdot \tfrac{1}{2} \nabla_x (\tfrac{\abs{\nabla_x^s f_\R}^2}{M_0}) }
    & = \tfrac{1}{2} \iint \nabla_x \cdot (v \tfrac{\abs{\nabla_x^s f_\R}^2}{M_0}) \d v \d x - \iint {\tfrac{1}{2} v \cdot \nabla_x M_0^{-1} \abs{\nabla_x^s f_\R}^2} \d v \d x \\\no
    & = \iint \tfrac{1}{2} v \otimes v : \nabla_x \Omega_0 \frac{\abs{\nabla_x^s f_\R}^2}{M_0} \d v \d x\,.
  \end{align}

Arguing as that in deriving the Poincar\'e inequality \eqref{eq:Poinc-ineq-Phi} in \Cref{lemm:Poincare-ineq}, we can get from the fact $|v|^2 \le |\nabla_v \Phi|^2$ that
  \begin{align}
    \Biggl| \iint & v \otimes v : \nabla_x \Omega_0 \frac{\abs{\nabla_x^s f_\R}^2}{M_0} \d v \d x \Biggr| \no\\
    & \le \nm{\nabla_x \Omega_0}_{L^\infty_x} {\iint |\nabla_v \Phi|^2 \abs{ \frac{\nabla_x^s f_\R}{M_0} }^2 M_0 \d v \d x}  \\\no
    & \le C_{\beta,\Omega_0} \left[ \iint \abs{\nabla_v \left( \frac{\nabla_x^s f_\R}{M_0} \right)}^2 M_0 \d v \d x + \iint \abs{\frac{\nabla_x^s f_\R}{M_0}}^2 M_0 \d v \d x \right] \\\no
    & \le C_{\beta,\Omega_0} \nm{\nabla_v (\tfrac{\nabla_x^s f_\R}{M_0})}^2_{\M} +\nm{\nabla_x^s f_\R}^2_{\MV}\,.
  \end{align}
Here we mention that, similar as in \cref{eq:vv-poincare}, an integral involving $\tfrac{\nabla_x^s f_\R}{M_0}$ itself arose due to its nonzero mean value. We also get
  \begin{align}
    \Biggl| \iint & v \p_t \Omega_0 \frac{\abs{\nabla_x^s f_\R}^2}{M_0} \d v \d x \Biggr| \no\\
    & \le \nm{\p_t \Omega_0}_{L^\infty_x}
          \left\{ \iint \frac{\abs{\nabla_x^s f_\R}^2}{M_0} \d v \d x \right\}^{\frac{1}{2}}
          \left\{ \iint |\nabla_v \Phi|^2 \abs{ \frac{\nabla_x^s f_\R}{M_0} }^2 M_0 \d v \d x \right\}^{\frac{1}{2}} \no\\
    & \le C_{\beta,\Omega_0}
            \left\{ \iint \frac{\abs{\nabla_x^s f_\R}^2}{M_0} \d v \d x \right\}^{\frac{1}{2}}
      \cdot \left\{ \iint \abs{\nabla_v \left( \frac{\nabla_x^s f_\R}{M_0} \right)}^2 M_0 \d v \d x + \iint \frac{\abs{\nabla_x^s f_\R}^2}{M_0} \d v \d x \right\}^{\frac{1}{2}} \no\\
    & \le C_{\beta,\Omega_0} \nm{\nabla_v (\tfrac{\nabla_x^s f_\R}{M_0})}^2_{\M} + \nm{\nabla_x^s f_\R}^2_{\MV}\,.
  \end{align}

Together with the higher-order derivatives estimate on linear operator $\L_{\Omega_0}$ in \cref{eq:coercive-HiOd}, we can control the left-hand side of \eqref{eq:Remd-HighDeri} from below, by
  \begin{align}
    & \sskp{(\p_t + v \cdot \nabla_x) \nabla_x^s f_\R}{\tfrac{\nabla_x^s f_\R}{M_0}}
      - \frac{1}{\eps} \sskp{ \nabla_x^s \L_{\Omega_0} f_\R }{\tfrac{\nabla_x^s f_\R}{M_0}}
  \\\no
  %   & \ge \frac{1}{2} \frac{\d}{\d t} \iint \frac{\abs{\nabla_x^s f_\R}^2}{M_0} \d v \d x
  %       + \frac{3}{4\eps} \iint \abs{\nabla_v (\tfrac{\nabla_x^s f_\R}{M_0})}^2 M_0 \d v \d x
  %       + \frac{3\alpha}{4c_0 \eps} \int \abs{P_{\Omega_0^\bot} {\nabla_x^s j_\R}}^2 \d x\\\no
  %   & \quad - C_{\beta} \nm{\nabla_x \Omega_0}_{L^\infty_x} \left[ \iint \abs{\nabla_v \left( \frac{f_\R}{M_0} \right)}^2 M_0 \d v \d x + \iint \abs{\frac{f_\R}{M_0}}^2 M_0 \d v \d x \right] \\\no
  %   & \quad - C_{\beta} \nm{\p_t \Omega_0}_{L^\infty_x}
  %           \left\{ \iint \frac{\abs{f_\R}^2}{M_0} \d v \d x \right\}^{\frac{1}{2}}
  %     \cdot \left\{ \iint \abs{\nabla_v (\tfrac{f_\R}{M_0})}^2 M_0 \d v \d x + \iint \frac{\abs{f_\R}^2}{M_0} \d v \d x \right\}^{\frac{1}{2}}
  % \\[5pt] \no
    & \ge \frac{1}{2} \frac{\d}{\d t} \nm{\nabla_x^s f_\R}^2_{\MV}
        + \frac{3}{4\eps} \nm{\nabla_v (\tfrac{\nabla_x^s f_\R}{M_0})}^2_{\M}
        + \frac{3\alpha}{4c_0 \eps} \nm{P_{\Omega_0^\bot} {\nabla_x^s j_\R}}_\lx^2
        \\\no
    & \qquad - \frac{1}{\eps} C_{\Omega_0} \sum_{\mathclap{s_1 \le s-1}} \nm{\nabla_x^{s_1} f_\R}_{\MV} \nm{\nabla_v (\tfrac{\nabla_x^s f_\R}{M_0})}_{\M}
        - C_{\beta,\Omega_0} \left( \nm{\nabla_v (\tfrac{\nabla_x^s f_\R}{M_0})}^2_{\M} + \nm{\nabla_x^s f_\R}_{\MV}^2 \right) \,.
      % \left[ \nm{\nabla_v (\tfrac{\nabla_x^s f_\R}{M_0})}^2_{\M} + \nm{\nabla_x^s f_\R}_{\MV} \nm{\nabla_v (\tfrac{\nabla_x^s f_\R}{M_0})}_{\M} + \nm{\nabla_x^s f_\R}_{\MV}^2 \right] \,.
  \end{align}
Multiplying by $\eps^s$ on the above inequality implies,
  \begin{align}\label{eq:HighDeri-LHS}
    \eps^s & \sskp{(\p_t + v \cdot \nabla_x) \nabla_x^s f_\R}{\tfrac{\nabla_x^s f_\R}{M_0}}
      - \eps^{s-1} \sskp{ \nabla_x^s \L_{\Omega_0} f_\R }{\tfrac{\nabla_x^s f_\R}{M_0}} \\\no
    & \ge \frac{1}{2} \frac{\d}{\d t} E_{s,\eps} + \frac{3}{4} D_{s,\eps} - C \sum_{\mathclap{s_1 \le s-1}} E_{s_1,\eps}^{\frac{1}{2}} D_{s,\eps}^{\frac{1}{2}} - C_{\beta,\Omega_0} (\eps D_{s,\eps} + E_{s,\eps}) \,.
  \end{align}

Considering the first term on the right-hand side of \cref{eq:Remd-HighDeri}, we can infer
  \begin{align}
    - & \sskp{ \div_v \left[ \nabla_x^s ( f_0 \widehat{\Omega}_{3,R}) \right] }{\tfrac{\nabla_x^s f_\R}{M_0}} \\\no
    & = \sum_{\mathclap{s_1+s_2 \le s}} \sskp{ \nabla_x^{s_2} f_0 \nabla_x^{s_1} \widehat{\Omega}_{3,R} }{\nabla_v (\tfrac{\nabla_x^s f_\R}{M_0})} \\\no &
    \le \sum_{\mathclap{s_1+s_2 \le s}} \nm{\nabla_x^{s_2} f_0}_{L^\infty_x L^2_v(\MV)}
            \nm{\nabla_x^{s_1} \widehat{\Omega}_{3,R}}_{\lx} \nm{ \nabla_v (\tfrac{\nabla_x^s f_\R}{M_0}) }_{\M} \\\no &
    % \le C_{f_0,f_1,f_2} \sum_{\mathclap{s_1 \le s}} \left[ 1+ \nm{\nabla_x^{s_1} f_\R}_{\lxv(\MV)} \right] \nm{\nabla_v (\tfrac{\nabla_x^s f_\R}{M_0})}_{\M} \\\no &
    \le C_{f_0,f_1,f_2} \sum_{\mathclap{s_1 \le s}} \Big[ 1+ \nm{f_\R}_{\H^{s_1}_xL^2_v(M^{-1})} + \eps \nm{f_\R}_{\H^{s_1}_xL^2_v(M^{-1})} \nm{f_\R}_{\H^2_x L^2_v(M^{-1})} \\\no & \hspace*{2.8cm} + \eps^3 \nm{f_\R}_{\H^{s_1}_xL^2_v(M^{-1})} \nm{f_\R}^2_{\H^2_x L^2_v(M^{-1})} \Big] \nm{\nabla_v (\tfrac{\nabla_x^s f_\R}{M_0})}_{\M} \,,
  \end{align}
from which, it follows,
  \begin{align}\label{eq:HighDeri-RHS1}
    - \eps^s \sskp{ \div_v [ \nabla_x^s ( f_0 \widehat{\Omega}_{3,R})] }{\tfrac{\nabla_x^s f_\R}{M_0}}
    \le C \sum_{\mathclap{s_1 \le s}} (1 + E_{s_1,\eps}^{\frac{1}{2}} + E_{s_1,\eps}^{\frac{1}{2}} E_{2,\eps}^{\frac{1}{2}}+ \eps E_{s_1,\eps}^{\frac{1}{2}} E_{2,\eps}) \eps^{\frac{1}{2}} \cdot D_{s,\eps}^{\frac{1}{2}} \,.
  \end{align}

As for the second term on the right-hand side of \cref{eq:Remd-HighDeri}, we get
  \begin{align}
    & - \sskp{ \div_v \left\{ \nabla_x^s [ (f_1 + \eps f_2) \widehat{\Omega}_\R^\eps ] \right\} }{\tfrac{\nabla_x^s f_\R}{M_0}} \\\no
    & \quad = \sskp{ \nabla_x^s \left\{ (f_1 + \eps f_2) \left[ (c_0\rho_0)^{-1} P_{\Omega_0^\bot} j_\R + \eps \widehat{\Omega}_{3,R} \right] \right\} }{ \nabla_v (\tfrac{\nabla_x^s f_\R}{M_0}) } \\\no
    & \quad = \sskp{ (f_1 + \eps f_2) (c_0\rho_0)^{-1} P_{\Omega_0^\bot} (\nabla_x^s j_\R) }{ \nabla_v (\tfrac{\nabla_x^s f_\R}{M_0}) } \\\no & \qquad
        + \sum_{\mathclap{\substack{s_1+s_2 \le s\\s_2 \ge 1}}} \sskp{ \nabla_x^{s_2} \left[ (f_1 + \eps f_2) (c_0\rho_0)^{-1} P_{\Omega_0^\bot} \right] \nabla_x^{s_1} j_\R }{ \nabla_v (\tfrac{\nabla_x^s f_\R}{M_0}) }
        \\\no & \qquad
        % + \eps \sskp{ (f_1 + \eps f_2) \nabla_x^s \widehat{\Omega}_{3,R} }{ \nabla_v (\tfrac{\nabla_x^s f_\R}{M_0}) }
        % + \eps \sskp{ \sum_{\mathclap{\substack{s_1+s_2 \le s\\s_2 \ge 1}}}
        + \eps \sum_{\mathclap{s_1+s_2 \le s}} \sskp{ \nabla_x^{s_2} (f_1 + \eps f_2) \nabla_x^{s_1} \widehat{\Omega}_{3,R} }{ \nabla_v (\tfrac{\nabla_x^s f_\R}{M_0}) }
        \\\no & \quad
      \le \nm{(f_1 + \eps f_2) (c_0\rho_0)^{-1}}_{L^\infty_x L^2_v(\MV)} \nm{P_{\Omega_0^\bot} (\nabla_x^s j_\R)}_{\lx} \nm{\nabla_v (\tfrac{\nabla_x^s f_\R}{M_0})}_{\M} \\\no & \qquad
            + \sum_{\mathclap{\substack{s_1+s_2 \le s\\s_2 \ge 1}}}
                  \nm{ \nabla_x^{s_2} [ (f_1 + \eps f_2) (c_0\rho_0)^{-1} P_{\Omega_0^\bot} ] }_{L^\infty_x L^2_v(\MV)}
                  \nm{\nabla_x^{s_1} j_\R}_{\lx} \nm{\nabla_v(\tfrac{\nabla_x^s f_\R}{M_0})}_{\M}\\\no & \qquad
            % + \nm{(f_1 + \eps f_2) (c_0\rho_0)^{-1}}_{L^\infty_x L^2_v(\MV)} \nm{P_{\Omega_0^\bot} (\nabla_x^s j_\R)}_{\lx} \nm{\nabla_v (\tfrac{\nabla_x^s f_\R}{M_0})}_{\M}
            %   \\\no & \qquad
            + \eps \sum_{\mathclap{s_1+s_2 \le s}} \nm{ \nabla_x^{s_2} (f_1 + \eps f_2)}_{L^\infty_x L^2_v(\MV)}
                    \nm{ \nabla_x^{s_1} \widehat{\Omega}_{3,R}}_{\lx} \nm{ \nabla_v (\tfrac{\nabla_x^s f_\R}{M_0}) }_{\M}
      \\\no & \quad
      \le C_{f_0,f_1,f_2} \nm{P_{\Omega_0^\bot} (\nabla_x^s j_\R)}_{\lx} \nm{\nabla_v (\tfrac{\nabla_x^s f_\R}{M_0})}_{\M}
          + C_{f_0,f_1,f_2} \sum_{\mathclap{\substack{s_1+s_2 \le s\\s_2 \ge 1}}} \nm{\nabla_x^{s_1} f_\R}_{\MV} \nm{\nabla_v (\tfrac{\nabla_x^s f_\R}{M_0})}_{\M} \\\no & \qquad
          + C_{f_0,f_1,f_2} \eps \sum_{\mathclap{s_1 \le s}} \nm{ \nabla_x^{s_1} \widehat{\Omega}_{3,R}}_{\lx} \nm{\nabla_v (\tfrac{\nabla_x^s f_\R}{M_0})}_{\M} \,.
  \end{align}
This implies, by multiplying by the above parameter $\nu^s = \eps^s$,
  \begin{align}\label{eq:HighDeri-RHS2}
    - & \eps^s \sskp{ \div_v \left\{ \nabla_x^s [ (f_1 + \eps f_2) \widehat{\Omega}_\R^\eps ] \right\} }{\tfrac{\nabla_x^s f_\R}{M_0}} \\\no
    & \le C \eps D_{s,\eps} + C \eps^\frac{1}{2} \sum_{\mathclap{s_1 \le s-1}} E_{s_1,\eps}^{\frac{1}{2}} \cdot D_{s,\eps}^{\frac{1}{2}} + C \eps^{\frac{3}{2}} \sum_{\mathclap{s_1 \le s}} (1 + E_{s_1,\eps}^{\frac{1}{2}} + E_{s_1,\eps}^{\frac{1}{2}} E_{2,\eps}^{\frac{1}{2}}+ \eps E_{s_1,\eps}^{\frac{1}{2}} E_{2,\eps}) \cdot D_{s,\eps}^{\frac{1}{2}} \,.
  \end{align}

Concerning the nonlinear term $\div_v \{ \nabla_x^s [ f_\R \tfrac{( \Omega_{f^\eps} - \Omega_0 )}{\eps} ] \}$, in the third term of the right-hand side of \cref{eq:Remd-HighDeri}, we can infer:
  \begin{align}
    - & \sskp{ \div_v \left\{ \nabla_x^s [ f_\R \tfrac{(\Omega_{f^\eps} - \Omega_0)}{\eps} ] \right\} }{\tfrac{\nabla_x^s f_\R}{M_0}} \\\no
    & = \sskp{ \nabla_x^s \left[ f_\R \tfrac{\left( \Omega_{f^\eps} - \Omega_0 \right)}{\eps} \right]  }{ \nabla_v (\tfrac{\nabla_x^s f_\R}{M_0}) } \\\no
    % & = \sum_{\mathclap{s_1+s_2 \le s}} \sskp{ \nabla_x^{s_1} f_\R \nabla_x^{s_2} \left( \tfrac{\Omega_{f^\eps} - \Omega_0}{\eps} \right) }{ \nabla_v (\tfrac{\nabla_x^s f_\R}{M_0}) } \\\no
    & \le C \left[ \nm{ \tfrac{\Omega_{f^\eps} - \Omega_0}{\eps} }_{H^2_x} \nm{f_\R}_{\H^s_xL^2_v(M^{-1})} + \nm{ \tfrac{\Omega_{f^\eps} - \Omega_0}{\eps} }_{H^s_x} \nm{f_\R}_{\H^2_x L^2_v(M^{-1})} \right] \nm{\nabla_v (\tfrac{f_\R}{M_0})}_{\M} \\\no
    & \le C_{f_0,f_1,f_2} \Big[ \nm{f_\R}_{\H^s_xL^2_v(M^{-1})} + \eps \nm{f_\R}_{\H^s_xL^2_v(M^{-1})} \nm{f_\R}_{\H^2_x L^2_v(M^{-1})} \\\no & \hspace*{4cm} + \eps^3 \nm{f_\R}_{\H^s_xL^2_v(M^{-1})} \nm{f_\R}^2_{\H^2_x L^2_v(M^{-1})} \Big] \nm{\nabla_v (\tfrac{f_\R}{M_0})}_{\M} \,,
  \end{align}
which results in,
  \begin{align}\label{eq:HighDeri-RHS3}
    - \eps^s \sskp{ \div_v \left\{ \nabla_x^s [ f_\R \tfrac{(\Omega_{f^\eps} - \Omega_0)}{\eps} ] \right\} }{\tfrac{\nabla_x^s f_\R}{M_0}}
    \le C \left( E_{s,\eps}^{\frac{1}{2}} + E_{s,\eps}^{\frac{1}{2}} E_{2,\eps}^{\frac{1}{2}} + \eps E_{s,\eps}^{\frac{1}{2}} E_{2,\eps} \right) \cdot \eps^{\frac{1}{2}} D_{s,\eps}^{\frac{1}{2}} \,.
  \end{align}

The last source term $\S(f_1,f_2)$ can be bounded by
  \begin{align}
    \eps^s \sskp{ \nabla_x^s \S(f_1,f_2) }{\tfrac{\nabla_x^s f_\R}{M_0}} \le \eps^s C(\S) \nm{\nabla_x^s f_\R}_{\MV} \le C E_{s,\eps}^{\frac{1}{2}} \,.
  \end{align}
This inequality, combined with \cref{eq:HighDeri-RHS1,eq:HighDeri-RHS2,eq:HighDeri-RHS3,eq:HighDeri-LHS}, enables us to obtain
  \begin{align}
    & \frac{1}{2} \frac{\d}{\d t} E_{s,\eps} + \frac{3}{4} D_{s,\eps} \\\no
    & \le C_{\beta,\Omega_0} (\eps D_{s,\eps} + E_{s,\eps})
        + C \sum_{\mathclap{s_1 \le s-1}} E_{s_1,\eps}^{\frac{1}{2}} D_{s,\eps}^{\frac{1}{2}}
        + C \sum_{\mathclap{s_1 \le s}} (1 + E_{s_1,\eps}^{\frac{1}{2}} + E_{s_1,\eps}^{\frac{1}{2}} E_{2,\eps}^{\frac{1}{2}}+ \eps E_{s_1,\eps}^{\frac{1}{2}} E_{2,\eps}) \cdot \eps^{\frac{1}{2}} D_{s,\eps}^{\frac{1}{2}} \\\no & \quad
        + C \eps D_{s,\eps} + C \eps^\frac{1}{2} \sum_{\mathclap{s_1 \le s-1}} E_{s_1,\eps}^{\frac{1}{2}} \cdot D_{s,\eps}^{\frac{1}{2}} + C \eps^{\frac{3}{2}} \sum_{\mathclap{s_1 \le s}} (1 + E_{s_1,\eps}^{\frac{1}{2}} + E_{s_1,\eps}^{\frac{1}{2}} E_{2,\eps}^{\frac{1}{2}}+ \eps E_{s_1,\eps}^{\frac{1}{2}} E_{2,\eps}) \cdot D_{s,\eps}^{\frac{1}{2}} \\\no & \quad
        + C \eps^{\frac{1}{2}} ( E_{s,\eps}^{\frac{1}{2}} + E_{s,\eps}^{\frac{1}{2}} E_{2,\eps}^{\frac{1}{2}} + \eps E_{s,\eps}^{\frac{1}{2}} E_{2,\eps} ) \cdot D_{s,\eps}^{\frac{1}{2}}
        + C E_{s,\eps}^{\frac{1}{2}}
    \\\no
    & \le \eps (C_{\beta,\Omega_0} + C) D_{s,\eps} + C (E_{s,\eps} + E_{s,\eps}^{\frac{1}{2}})
          + C (1 + \eps^{\frac{1}{2}}) \sum_{\mathclap{s_1 \le s-1}} E_{s_1,\eps}^{\frac{1}{2}} D_{s,\eps}^{\frac{1}{2}}
        \\\no & \quad
          + C (1 + \eps) \sum_{\mathclap{s_1 \le s}} (1 + E_{s_1,\eps}^{\frac{1}{2}} + E_{s_1,\eps}^{\frac{1}{2}} E_{2,\eps}^{\frac{1}{2}}+ \eps E_{s_1,\eps}^{\frac{1}{2}} E_{2,\eps}) \cdot \eps^{\frac{1}{2}} D_{s,\eps}^{\frac{1}{2}}
    \\\no
    & \le \eps (C_{\beta,\Omega_0} + C) D_{s,\eps} + C (1 + E_{s,\eps})
          + C \sum_{\mathclap{s_1 \le s-1}} E_{s_1,\eps}^{\frac{1}{2}} D_{s,\eps}^{\frac{1}{2}}
        \\\no & \quad
        + C \eps^{\frac{1}{2}} \sum_{\mathclap{s_1 \le s}} (1 + E_{s_1,\eps}^{\frac{1}{2}} + E_{s_1,\eps}^{\frac{1}{2}} E_{2,\eps}^{\frac{1}{2}} + \eps E_{s_1,\eps}^{\frac{1}{2}} E_{2,\eps}) \cdot D_{s,\eps}^{\frac{1}{2}} \,,
  \end{align}
which reduces to, by virtue of the H\"older inequality and by taking a sufficiently small parameter $\eps$,
  \begin{align}\label{esm:HighOrder}
    \frac{\d}{\d t} E_{s,\eps} + D_{s,\eps}
    & \le C (1 + E_{s,\eps}) + C^2 \sum_{\mathclap{s_1 \le s-1}} E_{s_1,\eps}
          + C \eps^{\frac{1}{2}} \sum_{\mathclap{s_1 \le s}} (1 + E_{s_1,\eps} + E_{s_1,\eps} E_{2,\eps} + \eps E_{s_1,\eps} E_{2,\eps}^2) \no\\
          % + C \sum_{\mathclap{s_1 \le s}} (E_{s_1,\eps} + E_{s_1,\eps} E_{2,\eps} + \eps^2 E_{s_1,\eps} E_{2,\eps}^2) \\\no
    % & \le C (1 + \sum_{\mathclap{s_1 \le s}} E_{s_1,\eps} ) + C \eps^{\frac{1}{2}} \sum_{\mathclap{s_1 \le s}} (1 + E_{s_1,\eps}^{\frac{1}{2}} + E_{s_1,\eps}^{\frac{1}{2}} E_{2,\eps}^{\frac{1}{2}}+ \eps E_{s_1,\eps}^{\frac{1}{2}} E_{2,\eps}) \cdot D_{s,\eps}^{\frac{1}{2}} \,.
    & \le C (1 + \sum_{\mathclap{s_1 \le s}} E_{s_1,\eps} ) + C \eps^{\frac{1}{2}} \sum_{\mathclap{s_1 \le s}} (E_{s_1,\eps} E_{2,\eps} + \eps^2 E_{s_1,\eps} E_{2,\eps}^2 ) \,.
  \end{align}
Noticing that the higher-order derivatives (with orders less than $s$) could be estimated by the same way, and moreover, by the same formulation as that in the highest order derivative estimate.

Finally, by using the definitions for energy and dissipation functionals \cref{eq:functionals-ED}, or speaking equivalently,
$$\E = \sum_{0 \le k \le s} E_{k,\eps}, \quad \text{ and } \D = \sum_{0 \le k \le s} D_{k,\eps} \,,$$
we summarize all estimates for $0 \le k \le s$ (recalling \cref{esm:ZeroOrder} for $k=0$ case), and obtain
  \begin{align}
    \frac{\d}{\d t} \E(t) + \D(t)
    \le C (1 + \E) + C \eps^{\frac{1}{2}} \left( \E^2 + \eps^2 \E^3 \right) \,.
  \end{align}
Note the coefficient $C>0$ depends on the known functions $f_0,\, f_1$ and $f_2$, all of which in turn depend only on the initial data $\|(\rho_0^{\IN}, \Omega_0^{\IN})\|_{H^{s+9}_x}$. The desired \emph{a-priori} estimate uniformly in $\eps$, \cref{esm:apriori-uniform} in \Cref{prop:uniform-esm}, is hence proved. \qed

% %% ------------------------------------------------------------------------ %%
% % \subsubsection{Non-zero property of moment} % (fold)
% % \label{ssub:non_zero_property_of_moment}

% % % subsubsection non_zero_property_of_moment (end)

\subsubsection{Persistence of moment} % (fold)
\label{ssub:persistence_of_moment}

% subsubsection persistence_of_moment (end)

Note that we have claimed in the above process that the bound of $\nm{j_\eps}_{L^\infty_x}$ from up and especially from below. Actually, we can give a more precise estimate on this quantity, see the following \Cref{prop:nonzero-j}.

It should be pointed out, this proposition also ensures the property of moment function $j_{f^\eps}$ away from zero data, provided that it holds initially. In other words, the orientation function $\Omega_{f^\eps} = \tfrac{j_{f^\eps}}{\abs{j_{f^\eps}}}$ have a persistence property.

\begin{proposition}[Persistence] \label{prop:nonzero-j}
  As for the moment function $j_{f^\eps} = \agl{v f^\eps}$, suppose that its initial data is away from zero, then this property can persist in the lifespan of solution.

  More precisely, there exists a time $T>0$, such that, for any $t\in [0,T]$,
  \begin{align}
    \tfrac{1}{2} \inf_{x \in \mathbb{R}^3} \abs{j^\IN_{f^\eps}} & \le \abs{j_{f^\eps}(t)} \le \tfrac{3}{2} \sup_{x \in \mathbb{R}^3} \abs{j^\IN_{f^\eps}} \,.
  \end{align}
\end{proposition}

\proof Concerning the total moment $j_{f^\eps} = \agl{v f^\eps}$ (denoted as $j_\eps$, for brevity), we turn back to the original kinetic Cucker-Smale equation \cref{eq:CS}, which can be recast as:
  \begin{align}\label{eq:total-f-eps}
    (\p_t + v \cdot \nabla_x) f^\eps = \frac{1}{\eps} \div_v \left[ M_0 \nabla_v (\tfrac{f^\eps}{M_0}) \right]
    - \div_v \left[ f^\eps \tfrac{\left( \Omega_{f^\eps} - \Omega_0 \right)}{\eps} \right] \,.
  \end{align}
Here we have used the fact $\div_v \left\{ \nabla_v f + f [(v- \Omega_0) + \nabla_v V] \right\} = \div_v \left[ M_0 \nabla_v (\tfrac{f}{M_0}) \right]$.

Multiply by $v$ on this equation and integrate with respect to $v \in \RR^3$, then it follows that,
  \begin{align}
    \p_t j_\eps + \agl{v \otimes v : \nabla_x f^\eps} + \frac{1}{\eps} \agl{M \nabla_v (\tfrac{f^\eps}{M})}
    = \agl{f^\eps} \tfrac{\left( \Omega_{f^\eps} - \Omega_0 \right)}{\eps} \,,
  \end{align}
where we have used the following facts, by integration by parts,
  \begin{align*}
    - \frac{1}{\eps} \agl{v \div_v \left[ M_0 \nabla_v (\tfrac{f^\eps}{M_0}) \right]} = \frac{1}{\eps} \agl{M \nabla_v (\tfrac{f^\eps}{M})},
  \quad \text{ and }
    - \agl{v \div_v \left[ f^\eps \tfrac{\left( \Omega_{f^\eps} - \Omega_0 \right)}{\eps} \right]}
    = \agl{f^\eps} \tfrac{\left( \Omega_{f^\eps} - \Omega_0 \right)}{\eps} .
  \end{align*}

By multiplying by $j^{n+1}$, we infer that
  \begin{align}
    \frac{1}{2} \frac{\d}{\d t} \abs{j_\eps}^2
    & = - \agl{v \otimes v : \nabla_x f^\eps} \cdot j_\eps
        - \frac{1}{\eps} \agl{M \nabla_v (\tfrac{f^\eps}{M})} \cdot j_\eps
        + \agl{f^\eps} \tfrac{\left( \Omega_{f^\eps} - \Omega_0 \right)}{\eps} \cdot j_\eps
    \\\no
    & = (\sum_{i=1,2,3} I_i) \cdot j_\eps \,,
  \end{align}
which enables us to get
  \begin{align}\label{eq:moment-bound-j-eps}
    - \abs{\sum_{i=1,2,3} I_i} \le \frac{\d}{\d t} \abs{j_\eps} \le \abs{\sum_{i=1,2,3} I_i}.
  \end{align}

We are left to estimate the terms $I_i$'s. The H\"older inequality implies
  \begin{align*}
    \abs{I_1} &
    \le \agl{|v \otimes v|^2 M}^\frac{1}{2} \agl{\frac{\abs{\nabla_x f^\eps}^2}{M}}^\frac{1}{2}
    \le C \nm{\nabla_x f^\eps}_{\lv(\MV)} \,, \\
    \abs{I_2} &
    \le \frac{1}{\eps} \agl{M}^\frac{1}{2} \agl{M \abs{\nabla_v (\tfrac{f^\eps}{M})}^2}^\frac{1}{2}
    \le \nm{\nabla_v (\tfrac{f_1 + \eps f_2}{M})}_{\lv(M)} + \eps \nm{\nabla_v (\tfrac{f_\R^\eps}{M})}_{\lv(M)} \,,
  \end{align*}
where we have used the expansion formula $f^\eps = \rho_0 M_0 + \eps f_1 + \eps^2 (f_2 + f_\R^\eps)$ in the last line.

As for the last term $I_3$, notice again Taylor's formula with a remainder $\tfrac{(\Omega_{f^\eps} - \Omega_0)}{\eps} = \tfrac{\d}{\d \eps} \Omega_{f^\eps} |_{\eps= \eta}$, and the orthogonal relation $\tfrac{\d}{\d \eps} \Omega_{f^\eps} = P_{\Omega_{f^\eps}^\bot} \tfrac{j_1 + \eps (j_2 + j_R)}{\abs{j_\eps}} \perp \Omega_{f^\eps}= \tfrac{j_\eps}{\abs{j_\eps}} $, due to \cref{eq:Omega-Remd-1}, then we can get $I_3 = 0$.
  % \begin{align}
  %   \tfrac{\left( \Omega_{f^\eps} - \Omega_0 \right)}{\eps} =
  %   \frac{\d}{\d \eps} \Omega_{f^\eps}
  %   = P_{\Omega_{f^\eps}^\bot} \frac{j_1 + \eps (j_2 + j_R)}{\abs{j_\eps}}
  % \end{align}

Inserting the above inequalities into \cref{eq:moment-bound-j-eps}, we can infer that
  \begin{align}
    \frac{\d}{\d t} \abs{j_\eps}
    & \ge - C_{f_1,f_2} \left( 1 + \nm{\nabla_x f^\eps}_{\lv(\MV)} + \eps \nm{\nabla_v (\tfrac{f_\R^\eps}{M})}_{\lv(M)} \right) \no\\
    & \ge - C_{f_0,f_1,f_2} \left( 1 + \eps^2 \nm{\nabla_x f_\R^\eps}_{\lv(\MV)} + \eps \nm{\nabla_v (\tfrac{f_\R}{M})}_{\lv(M)} \right),
  \end{align}
which immediately implies
  \begin{align}\label{eq:lower-j-eps}
    \abs{j_\eps(t)}
    & \ge \abs{ j_\eps^{\IN} } - C \int_0^t \left( 1 + \eps^2 \nm{\nabla_x f_\R^\eps}_{\lv(\MV)} + \eps \nm{\nabla_v (\tfrac{f_\R^\eps}{M})}_{\lv(M)} \right) \d t \\\no
    & \ge \inf_{x \in \mathbb{R}^3} \abs{ j_\eps^{\IN} }
          - C t - C \int_0^t \left( \eps^2 \sup_{x \in \mathbb{R}^3} \nm{\nabla_x f_\R^\eps}_{\lv(\MV)} + \eps \sup_{x \in \mathbb{R}^3} \nm{\nabla_v (\tfrac{f_\R^\eps}{M})}_{\lv(M)} \right) \d t \,.
  \end{align}

By straightforward calculations, we get
  \begin{align}\label{eq:Fubini-pre}
    \sup_{x} \nm{g}_{\lv(\MV)}
    \le \nm{ \nm{ \tfrac{g}{\sqrt{M}} }_{L^\infty_x} }_{\lv}
    & \le C \nm{ \nm{ \tfrac{g}{\sqrt{M}} }_{\hx{2}} }_{\lv} \\\no
    & \le C \left( \nm{ \nabla_x^2 (\tfrac{g}{\sqrt{M}}) }_{\lxv} + \nm{ \nabla_x (\tfrac{g}{\sqrt{M}}) }_{\lxv} + \nm{ \tfrac{g}{\sqrt{M}} }_{\lxv} \right)  \\\no
    & = C \sum_{i=1,2} \xi_i + C \nm{ \tfrac{g}{\sqrt{M}} }_{\lxv} \,.
  \end{align}
As for the first term $\xi_1$, we write
  \begin{align}
    \xi_1 & = \nm{ \nabla_x^2 (\tfrac{g}{\sqrt{M}}) }_{\lxv} \\\no
    & \le \nm{ \tfrac{\nabla_x^2 g}{\sqrt{M}} }_{\lxv} + \nm{ \tfrac{\nabla_x g}{\sqrt{M}} (v_\bot \cdot \nabla_x \Omega_0) }_{\lxv} + \nm{ \tfrac{g}{\sqrt{M}} \left[ (v_\bot \nabla_x \Omega_0)\otimes (v_\bot \nabla_x \Omega_0) - (v_\bot \nabla_x^2 \Omega_0) \right] }_{\lxv}.
  \end{align}
In order to control the last term
  \begin{align*}
    \nm{ \tfrac{g}{\sqrt{M}} (v_\bot \nabla_x \Omega_0)\otimes (v_\bot \nabla_x \Omega_0) }_{\lxv}
    = \left\{\iint |v \otimes v|^2 \frac{\abs{g}^2}{M_0} \d v \cdot |\nabla_x \Omega_0 \otimes \nabla_x \Omega_0|^2 \d x \right\}^{\frac{1}{2}},
  \end{align*}
recall the fact $|v|^4 \le |\nabla_v \Phi|^2$, then we can apply the same argument as that for the Poincar\'e inequality \eqref{eq:Poinc-ineq-Phi} in \Cref{lemm:Poincare-ineq} to get
  \begin{align}\label{eq:vvvv-poincare}
    & \abs{\iint |v|^4 \frac{\abs{g}^2}{M_0} \d v |\nabla_x \Omega_0 \otimes \nabla_x \Omega_0|^2\d x} \\\no
    & \le C \nm{\nabla_x \Omega_0}^4_{L^\infty_x} {\iint |\nabla_v \Phi|^2 \abs{ \frac{g}{M_0} }^2 M_0 \d v \d x}  \\\no
    & \le C \nm{\nabla_x \Omega_0}^4_{L^\infty_x} \left[ \iint \abs{\nabla_v \left( \frac{g}{M_0} \right)}^2 M_0 \d v \d x + \iint \abs{\frac{g}{M_0}}^2 M_0 \d v \d x \right] \,.
  \end{align}
We point out that, an additional integral involving $\tfrac{g}{M_0}$ itself arose here, due to its nonzero mean value $\agl{\tfrac{g}{M_0}}_\M = \rho_g \neq 0$, comparing \cref{eq:potential-poinc}. This implies
  \begin{align}
    \nm{ \tfrac{g}{\sqrt{M}} (v_\bot \nabla_x \Omega_0)\otimes (v_\bot \nabla_x \Omega_0) }_{\lxv}
    \le C \left( \nm{\nabla_v (\tfrac{g}{M})}_{\lxv(\M)} + \nm{g}_{\lxv(\MV)} \right).
  \end{align}
Similarly, it follows
  \begin{align*}
    \nm{ \tfrac{g}{\sqrt{M}} (v_\bot \nabla_x^2 \Omega_0) }_{\lxv} & \le C \left( \nm{\nabla_v (\tfrac{g}{M})}_{\lxv(\M)} + \nm{g}_{\lxv(\MV)} \right) \,, \\
    \nm{ \tfrac{\nabla_x g}{\sqrt{M}} (v_\bot \nabla_x \Omega_0) }_{\lxv} & \le C \left( \nm{\nabla_v (\tfrac{\nabla_x g}{M})}_{\lxv(\M)} + \nm{\nabla_x g}_{\lxv(\MV)} \right) \,,
  \end{align*}
so we can infer that
  \begin{align}
    \xi_1 =\nm{ \nabla_x^2 (\tfrac{g}{\sqrt{M}}) }_{\lxv}
    \le C \left( \nm{g}_{\hxv{2}(\MV)} + \nm{\nabla_v (\tfrac{g}{M})}_{\hxv{1}(\M)} \right).
  \end{align}
The term $\xi_2$ can be dealt with in a similar and easier manner, and it follows
  \begin{align}
    \xi_2 = \nm{ \nabla_x (\tfrac{g}{\sqrt{M}}) }_{\lxv}
    & \le \nm{ \tfrac{\nabla_x g}{\sqrt{M}} }_{\lxv} + \nm{ \tfrac{g}{\sqrt{M}} (v_\bot \cdot \nabla_x \Omega_0) }_{\lxv} \\\no
    & \le C \left( \nm{g}_{\hxv{1}(\MV)} + \nm{\nabla_v (\tfrac{g}{M})}_{\lxv(\M)} \right).
  \end{align}
Plugging the two estimates for $\xi_1$ and $\xi_2$ into \cref{eq:Fubini-pre}, we get, for $g = \nabla_x f_\R$:
  \begin{align}\label{eq:Fubini}
    \sup_{x} \nm{g}_{\lv(\MV)} \le C \left( \nm{g}_{\hxv{2}(\MV)} + \nm{\nabla_v (\tfrac{g}{M})}_{\hxv{1}(\M)} \right).
  \end{align}
Perform a similar argument, we can also derive that
  \begin{align}\label{eq:Fubini-2}
    \sup_{x} \nm{\nabla_v (\tfrac{f_\R}{M})}_{\lv(M)}
    & \le \nm{ \nm{ \nabla_v (\tfrac{f_\R}{M}) \sqrt{M} }_{L^\infty_x} }_{\lv}
      \le C \nm{ \nm{ \nabla_v (\tfrac{f_\R}{M}) \sqrt{M} }_{\hx{2}} }_{\lv} \\\no
    & \le C \sum_{k=0,1,2} \nm{ \nabla_v \nabla_x^k (\tfrac{f_\R}{\sqrt{M}}) \sqrt{M} }_{\lxv}
            % \left( \nm{ \nabla_v \nabla_x^2 (\tfrac{f_\R}{\sqrt{M}}) \sqrt{M} }_{\lxv}
            % + \nm{ \nabla_v \nabla_x (\tfrac{f_\R}{\sqrt{M}}) \sqrt{M}}_{\lxv} + \nm{ \nabla_v (\tfrac{f_\R}{\sqrt{M}}) \sqrt{M} }_{\lxv} \right) \\\no
            + C \nm{ \nabla_v \nabla_x (\tfrac{f_\R}{\sqrt{M}}) \nabla_x \sqrt{M} }_{\lxv} \\\no & \quad \
            + C \nm{ \nabla_v (\tfrac{f_\R}{\sqrt{M}}) \nabla_x^2 \sqrt{M} }_{\lxv}
            + C \nm{ \nabla_v (\tfrac{f_\R}{\sqrt{M}}) \nabla_x \sqrt{M} }_{\lxv} \\\no
    & \le C \left( \nm{\nabla_v (\tfrac{f_\R}{M})}_{\hxv{2}(\M)} + \nm{f_\R}_{\hxv{1}(\MV)} \right) \,.
  \end{align}
Combining \cref{eq:Fubini} for $g = \nabla_x f_\R$, and \cref{eq:Fubini-2} together imply
  \begin{align}\label{eq:lower-useful}
    & \eps^2 \sup_{x \in \mathbb{R}^3} \nm{\nabla_x f_\R}_{\lv(\MV)} + \eps \sup_{x \in \mathbb{R}^3} \nm{\nabla_v (\tfrac{f_\R}{M})}_{\lv(M)} \\\no
    & \le C \eps^2 \left( \nm{\nabla_x f_\R}_{\hxv{2}(\MV)} + \nm{\nabla_v (\tfrac{\nabla_x f_\R}{M})}_{\hxv{1}(\M)} \right)
        + C \eps \left( \nm{\nabla_v (\tfrac{f_\R}{M})}_{\hxv{2}(\M)} + \nm{f_\R}_{\hxv{1}(\MV)} \right) \\\no
    & \le C \eps^\frac{1}{2} \sum_{k=0,1,2,3} \eps^\frac{k}{2} \nm{\nabla_x^k f_\R}_{\lxv(\MV)} + C \eps^\frac{1}{2} \sum_{k=0,1,2} \eps^\frac{k-1}{2} \nm{\nabla_v (\tfrac{\nabla_x^k f_\R}{M})}_{\lxv(\M)} \,.
    % \\\no
    % & \le C \eps^\frac{1}{2} \nm{f_\R^{\eps,\IN}}_{\hxv{3}(\MV)} \,.
  \end{align}

On the other hand, note that we have used the uniform-in-$\eps$ estimate in \Cref{prop:uniform-esm} yields that the energy bound (taking $s=3$):
  \begin{align}
    \E(t) + \D(t) \le C \E(0),
  \end{align}
i.e.,
  \begin{align}
    \sum_{k=0,1,2,3} \eps^k \nm{\nabla_x^k f_\R}^2_{\lxv(\MV)} + \int_0^t \sum_{k=0,1,2,3} \eps^{k-1} \nm{\nabla_v (\tfrac{\nabla_x^k f_\R}{M})}^2_{\lxv(\M)} \d t
    \le C \nm{f_\R^{\eps,\IN}}^2_{\hxv{3}(\MV)} \,.
  \end{align}
By inserting \cref{eq:lower-useful} into \cref{eq:lower-j-eps}, we can infer
  \begin{align}
    \abs{j_\eps(t)}
    & \ge \abs{ j_\eps^{\IN} } - Ct - C \eps^\frac{1}{2} \int_0^t \sum_{k=0,1,2,3} \eps^\frac{k}{2} \nm{\nabla_x^k f_\R}_{\lxv(\MV)} \d t \no\\
    & \quad - C \eps^\frac{1}{2} \int_0^t \sum_{k=0,1,2} \eps^\frac{k-1}{2} \nm{\nabla_v (\tfrac{\nabla_x^k f_\R}{M})}_{\lxv(\M)} \d t
    \no\\
    & \ge \inf_{x \in \mathbb{R}^3} \abs{ j_\eps^{\IN} }
          - C t
          - C t \eps^\frac{1}{2} \sum_{k=0,1,2,3} \eps^\frac{k}{2} \nm{\nabla_x^k f_\R}_{L^\infty_t\lxv(\MV)} \no\\ & \qquad
          - C t^\frac{1}{2} \eps^\frac{1}{2} \left[ \int_0^t \sum_{k=0,1,2} \eps^{k-1} \nm{\nabla_v (\tfrac{\nabla_x^k f_\R}{M})}^2_{\lxv(\M)} \d t \right]^\frac{1}{2} \no\\
    & \ge \inf_{x \in \mathbb{R}^3} \abs{ j_\eps^{\IN} }
          - C t
          - C (t + t^\frac{1}{2}) \eps^\frac{1}{2} \nm{f_\R^{\eps,\IN}}_{\hxv{3}(\MV)} \,.
  \end{align}

Thus, there exists a positive time $T$:
  \begin{align}
   T = \min \left\{ \tfrac{\abs{ j_\eps^{\IN} }}{4C (1 + \eps^\frac{1}{2} \nm{f_\R^{\eps,\IN}}_{\hxv{3}(\MV)})},\, \left[ \tfrac{\abs{ j_\eps^{\IN} }}{ 4C \eps^\frac{1}{2} \nm{f_\R^{\eps,\IN}}_{\hxv{3}(\MV)} } \right] ^2 \right\} >0 \,,
   \end{align}
which eventually depends only on the initial data $f^{\eps,\IN}$, such that for any $t \in [0,T]$,
  \begin{align}
    \abs{j_\eps(t)} \ge \inf_{x \in \mathbb{R}^3} \abs{ j_\eps^{\IN} } - C t \left( 1 + \eps^\frac{1}{2} \nm{f_\R^{\eps,\IN}}_{\hxv{3}(\MV)} \right)  - C t^\frac{1}{2} \eps^\frac{1}{2} \nm{f_\R^{\eps,\IN}}_{\hxv{3}(\MV)}
    \ge \tfrac{1}{2} \inf_{x \in \mathbb{R}^3} \abs{j^{n+1,\IN}}.
  \end{align}
Following the same line as above, we can get the upper bound for $\abs{j_\eps(t)}$, and hence we conclude that
  \begin{align}
    \tfrac{1}{2} \inf_{x \in \mathbb{R}^3} \abs{ j_\eps^{\IN} } \le \abs{j_\eps(t)} \le \tfrac{3}{2} \sup_{x \in \mathbb{R}^3} \abs{ j_\eps^{\IN} } \,.
  \end{align}
This completes the whole proof of \Cref{prop:nonzero-j}.

\endproof

%% ------------------------------------------------------------------------ %%
\subsection{Completion of proof of main result} % (fold)
\label{sub:completion_of_main_result}

% subsection completion_of_main_result (end)

% \smallskip\noindent\emph{Proof of \Cref{thm:limit}}.
We give here the proof of \Cref{thm:limit} by virtue of \Cref{prop:uniform-esm}. Given initial data $f_R^{\eps,\IN}$, it follows from the local existence result of the remainder equation \eqref{eq:Rem-order} (see Lemma \ref{lemm:Rem-local} below) that, there exists a unique local solution $f_R^\eps(t,x,\omega)$ on $[0,\ T_\eps]$ with $T_\eps > 0$ be the maximal lifespan.

Denote $\widetilde{E}_1 = 2 e^{C T_1} (\E(0) + 2 C T_1)$ with $T_1 \in (0, T_0]$ being the lifespan of the solution to the next order macroscopic system \cref{sys:limit-1st-order-rho,sys:limit-1st-order-j}, (recall here $T_0$ is the lifespan of the hydrodynamic system \cref{sys:limit-0-order}). Moreover, let $\widetilde{T}_1 = \sup_{t \in [0,\ T_\eps]} \left\{ \E(t) \le \widetilde{E}_1 \right\}$.

If $T_\eps < T_1$, the Gr\"onwall inequality, together with the \emph{a-priori} estimate \Cref{prop:uniform-esm}, enable us to get, for any $t \in [0,\ \widetilde{T}_1]$, that
  \begin{align}
    \E(t) \le e^{C t} (\E(0) + 2 C t) \le \tfrac{1}{2} \widetilde{E}_1,
  \end{align}
by noticing that $\eps^{\frac{1}{2}} (\E^2 + \eps^2 \E^3) \le \eps^{\frac{1}{2}} (\widetilde{E}_1^2 + \eps^2 \widetilde{E}_1^3) \le 1$, for sufficiently small $\eps$. This implies that the solution can be continued beyond the time $T_\eps$, which is a contradiction with the maximal property of $T_\eps$. So, the maximal existing time can be continued beyond $T_1$, i.e., $T_\eps \ge T_1$, and hence the whole proof of Theorem \ref{thm:limit} is completed. \qed

%% ------------------------------------------------------------------------ %%
\subsection{Local existence of remainder equation} % (fold)
\label{sub:local_existence_of_remainder_equation}

% subsection local_existence_of_remainder_equation (end)

In the sequel we study the local existence of remainder equation \eqref{eq:Rem-order} (or equivalently, \eqref{eq:Rem-order-re}), as follows:
\begin{lemma}[Local Existence] \label{lemm:Rem-local}
  Let $s\ge 3$. There exists some $T_* > 0$ (depending on $\eps>0$), such that the Cauchy problem of the remainder equation \eqref{eq:Rem-order} admits a unique solution $f_R^\eps$, satisfying
  \begin{align}
    f_R^\eps (t,x,v) \in L^\infty(0, T_*;\ \H^s_x L^2_v(M^{-1}) ) \,. % \cap L^2(0, T_*;\ H^s_x H^1_v(M^{-1}) ) \,,
  \end{align}

\end{lemma}

% Indeed, define the energy and dissipation functionals as follows,
%   \begin{align}\label{eq:functionals-local}
%     \mathbb{E}(t) = \sum_{0\le k \le s} \nm{\tfrac{\nabla_x^k f_\R^\eps}{M_0}}_\M^2\,, \quad
%     \mathbb{D}(t) = \sum_{0\le k\le s} \nm{\nabla_v (\tfrac{\nabla_x^k f_\R^\eps}{M_0})}_\M^2 + \frac{\alpha}{c_0} \nm{P_{\Omega_0^\bot} {\nabla_x^k j_\R^\eps}}_{L^2_x}^2 \,.
%   \end{align}
% then we have
%   \begin{align}\label{esm:apriori-local}
%     \tfrac{\d}{\d t} \mathbb{E}(t) + \mathbb{D}(t) \le C \left( 1 + \mathbb{E}(t) + \mathbb{E}^2(t) + \mathbb{E}^3(t) \right) \,,
%   \end{align}
% which implies the local existence result in \Cref{lemm:Rem-local}.

\smallskip\noindent\emph{Sketch Proof of \Cref{lemm:Rem-local}}. The proof of \Cref{lemm:Rem-local} proceeds by an iteration scheme. More precisely, the $(n+1)$-th step iteration scheme can be constructed as follows,
  \begin{multline}\label{eq:Rem-order-iteration}
    (\p_t + v \cdot \nabla_x) f_\R^{n+1} - \frac{1}{\eps} \L_{\Omega_0} f_\R^{n+1}
    = - \div_v \left[ f_\R^{n+1} \tfrac{\left( \Omega_{f}^{\eps,n} - \Omega_0 \right)}{\eps} \right] \\
      - \div_v \left( f_0 \widehat{\Omega}_{3,R}^{n} \right)
      - \div_v \left[ (f_1 + \eps f_2) \widehat{\Omega}_\R^{\eps,n} \right]
    % \\\no
        - \S(f_1,f_2) \,.
  \end{multline}
where we have still dropped the superscript $\eps$ in $f_\R^{\eps}$ for brevity. The iteration begins with a given $f_R^0$, and the initial data is given by $f_\R^{n+1} \big|_{t=0} = f_R^{n+1,\,\IN}$.

% Note that the first equation is linear with respect to the distribution function $f$, so given a $v_n$, we can determine a $f_n$ by using the first equation; Next combining the given $v_n$ and the determined $f_n$, we can determine $v_{n+1}$ by solving the equation of velocity field.

% \begin{align} \label{eq:iteration} \tag*{(Re)$_{n+1}$}
% \begin{cases}
% \p_t f_n + u_0 \nabla_x f_n + \nabla_{\omega}\cdot(\mathcal{F}_0 f_n)
% = \frac{1}{\eps}\mathcal{L}_{\Omega_0} f_n - \frac{1}{\sqrt\eps} h_0 - h_1(v_n, f_n),
% \\[3pt]
% \p_t v_{n+1} + v_0 \cdot \nabla_x v_{n+1} + v_{n+1} \cdot \nabla_x v_0 + \sqrt\eps v_n \cdot \nabla_x v_{n+1} \\[3pt]
% \hspace*{5cm} + \nabla_x\cdot Q_{f_n} = -\nabla_x p_{n+1} + \Delta_x v_{n+1},
% \\[3pt]
% \nabla_x \cdot v_{n+1} =0,
% \end{cases}
% \end{align}
% where $h_1(v_n, f_n) = v_n \cdot \nabla_x (f_0 + \sqrt\eps f_n) + \nabla_\omega \cdot [\mathcal{P}_{\omega^\perp}B(v_n) \omega (f_0 + \sqrt\eps f_n)]$.

Our task is to prove that the sequence $\{f_\R^{n}\}_{n=0}^{\infty}$ is a Cauchy sequence, then the limit point $f_R$ will be the solution we are looking for. For that, we define the solution set $\mathcal{S}_T(A)$ by
\begin{align}\label{eq:local-space}
\mathcal{S}_T(A)
= \left\{
    f_\R (t,x,v) \in L^\infty(0, T_*;\ \H^s_x L^2_v(M^{-1}) ) \,
  \left|\, % \sup_{t \in [0, T)} \E(t) = \le A
    \begin{aligned}
      \sup_{t \in [0, T)} \E(t) & + \int_0^T \D(t) \d t \le A\,;
    \\[3pt]
      \tfrac{1}{2} \inf_{x \in \mathbb{R}^3} \abs{j^\IN_{f^\eps}} & \le \abs{j_{f^\eps}(t)} \le \tfrac{3}{2} \sup_{x \in \mathbb{R}^3} \abs{j^\IN_{f^\eps}} \,.
    \end{aligned}
  \right.
  \right\} \,,
\end{align}
where the constant $A$ could be determined. Note that, we will use, in the $(n+1)$-th step, the upper and lower bound of the total moment $j_\eps^n=j_{f^{\eps,n}}$ coming from the $n$-th step, in order to control the quantities $(\Omega_{f}^{\eps,n} - \Omega_0)/\eps$, $\widehat{\Omega}_{3,\R}^{n}$, and $\widehat{\Omega}_\R^{\eps,n}$ on the right-hand side of \cref{eq:Rem-order-iteration}. So we have employed this condition in the above solution space.

Assuming we have got $f_\R^n \in \mathcal{S}_T(A)$, the energy bound of the $(n+1)$-th step in \cref{eq:local-space} follows a similar process as that in deriving the \emph{a-priori} estimate \Cref{prop:uniform-esm}, so it holds $\sup_{t \in [0, T)} \E(f_\R^{n+1}) + \int_0^T \D(f_\R^{n+1}) \d t \le A$.

We next prove $j_{f^\eps}^{n+1}$ is bounded from up and below. Firstly, as a correspondence to the approximate remainder equation \cref{eq:Rem-order-iteration}, the function $f^{\eps,n+1} $ satisfies:
  \begin{align}
    (\p_t + v \cdot \nabla_x) f^{\eps,n+1} = \frac{1}{\eps} \div_v \left\{ \sigma \nabla_v f^{\eps,n+1} + f^{\eps,n+1} \left[ (v- \Omega_{f^{\eps,n}}) + \nabla_v V \right] \right\} \,,
  \end{align}
which could be recast as (the superscript $\eps$ is also omitted):
  \begin{align}
    (\p_t + v \cdot \nabla_x) f^{n+1} = \frac{1}{\eps} \div_v \left[ M_0 \nabla_v (\tfrac{f^{n+1}}{M_0}) \right]
    - \div_v \left[ f^{n+1} \tfrac{\left( \Omega_{f}^{\eps,n} - \Omega_0 \right)}{\eps} \right] \,.
  \end{align}
This is similar as \cref{eq:total-f-eps} in \Cref{prop:nonzero-j}, so the same argument can be also applied here. It follows that,
  \begin{align}
    \p_t j^{n+1} = - \agl{v \otimes v : \nabla_x f^{n+1}} - \frac{1}{\eps} \agl{M \nabla_v (\tfrac{f^{n+1}}{M})}
    + \agl{f^{n+1}} \tfrac{\left( \Omega_{f}^{\eps,n} - \Omega_0 \right)}{\eps}
    = \sum_{i=1,2,3} I_i\,,
  \end{align}
which enables us to get
  \begin{align}\label{eq:moment-bound}
    - \abs{\sum_{i=1,2,3} I_i} \le \frac{\d}{\d t} \abs{j^{n+1}} \le \abs{\sum_{i=1,2,3} I_i}.
  \end{align}
The estimates on $I_i$'s are also similar as before. The H\"older inequality implies
  \begin{align*}
    \abs{I_1} &
    \le \agl{|v \otimes v|^2 M}^\frac{1}{2} \agl{\frac{\abs{\nabla_x f^{n+1}}^2}{M}}^\frac{1}{2}
    \le C \nm{\nabla_x f^{n+1}}_{\lv(\MV)} \,, \\
    \abs{I_2} &
    \le \frac{1}{\eps} \agl{M}^\frac{1}{2} \agl{M \abs{\nabla_v (\tfrac{f^{n+1}}{M})}^2}^\frac{1}{2}
    \le \nm{\nabla_v (\tfrac{f_1 + \eps f_2}{M})}_{\lv(M)} + \eps \nm{\nabla_v (\tfrac{f_\R^{n+1}}{M})}_{\lv(M)} \,,
  \end{align*}
and for $I_3$,
  \begin{align}
    \abs{I_3}
    \le \agl{M}^\frac{1}{2} \agl{\frac{\abs{f^{n+1}}^2}{M}}^\frac{1}{2} \cdot \abs{\tfrac{\left( \Omega_{f}^{\eps,n} - \Omega_0 \right)}{\eps}}
    \le C_{j^n} \nm{f^{n+1}}_{\lv(\MV)} \,.
  \end{align}
Here we have used the upper and lower bound of $j^n$ due to the induction assumption $f_\R^n \in \mathcal{S}_T(A)$.

Combining the above inequalities and \cref{eq:moment-bound}, we can infer that
  \begin{align}
    \frac{\d}{\d t} \abs{j^{n+1}}
    & \ge - C \left( 1 + \nm{f^{n+1}}_{\lv(\MV)} + \nm{\nabla_x f^{n+1}}_{\lv(\MV)} + \eps \nm{\nabla_v (\tfrac{f_\R^{n+1}}{M})}_{\lv(M)} \right) \no\\
    & \ge - C \left( 1 + \eps^2 \nm{f_\R^{n+1}}_{\lv(\MV)} + \eps^2 \nm{\nabla_x f_\R^{n+1}}_{\lv(\MV)} + \eps \nm{\nabla_v (\tfrac{f_\R^{n+1}}{M})}_{\lv(M)} \right),
  \end{align}
with the constant $C$ depending on the known functions $f_0,f_1,f_2$ and the upper and lower bound of $j^n$.
This immediately yields
  \begin{align}
    \abs{j^{n+1}(t)}
    & \ge \abs{ j^{n+1,\IN} } - Ct - \int_0^t \left(\eps^2 \nm{f_\R^{n+1}}_{\hx{1}\lv(\MV)} + \eps \nm{\nabla_v (\tfrac{f_\R^{n+1}}{M})}_{\lv(M)} \right) \d t \\\no
    & \ge \abs{ j^{n+1,\IN} } - Ct
          - C \eps^\frac{1}{2} \int_0^t \sum_{k=0,1,2,3} \eps^\frac{k}{2} \nm{\nabla_x^k f_\R^{n+1}}_{\lxv(\MV)} \d t \\\no & \qquad
          - C \eps^\frac{1}{2} \int_0^t \sum_{k=0,1,2} \eps^\frac{k-1}{2} \nm{\nabla_v (\tfrac{\nabla_x^k f_\R^{n+1}}{M})}_{\lxv(\M)} \d t
    \\\no
    & \ge \inf_{x \in \mathbb{R}^3} \abs{ j_\eps^{\IN} }
          - C t
          - C t \eps^\frac{1}{2} \sum_{k=0,1,2,3} \eps^\frac{k}{2} \nm{\nabla_x^k f_\R^{n+1}}_{L^\infty_t\lxv(\MV)} \\\no & \qquad
          - C t^\frac{1}{2} \eps^\frac{1}{2} \left[ \int_0^t \sum_{k=0,1,2} \eps^{k-1} \nm{\nabla_v (\tfrac{\nabla_x^k f_\R^{n+1}}{M})}^2_{\lxv(\M)} \d t \right]^\frac{1}{2} \,.
          % \no\\
    % & \ge \inf_{x \in \mathbb{R}^3} \abs{ j_\eps^{\IN} }
    %       - C t
    %       - C (t + t^\frac{1}{2}) \eps^\frac{1}{2} \nm{f_\R^{\eps,\IN}}_{\hxv{3}(\MV)} \,.
  \end{align}
The established energy bound for $(n+1)$-th step enables that, there exists a short time $T>0$, depending on the bound constant $A$ and the initial data $\abs{j^{n+1,\IN}} = \agl{v f^{n+1,\IN}}$, such that for any $t \in [0,T]$,
  \begin{align}
    \abs{j^{n+1}(t)} \ge \inf_{x \in \mathbb{R}^3} \abs{j^{n+1,\IN}} - C t(1+ \eps^\frac{1}{2} A^\frac{1}{2}) - C t^\frac{1}{2} \eps^\frac{1}{2} A^\frac{1}{2}
    \ge \tfrac{1}{2} \inf_{x \in \mathbb{R}^3} \abs{j^{n+1,\IN}}.
  \end{align}
On the other hand, the upper bound for $\abs{j^{n+1}(t)}$ follows the same lines as above, and hence we conclude that
  \begin{align}
    \tfrac{1}{2} \inf_{x \in \mathbb{R}^3} \abs{j^{n+1,\IN}} \le \abs{j^{n+1}(t)} \le \tfrac{3}{2} \sup_{x \in \mathbb{R}^3} \abs{j^{n+1,\IN}} \,.
  \end{align}
This ensures us to complete the iteration process, and leads to the existence result of \Cref{lemm:Rem-local} by employing a standard compactness argument.
\qed

%% ------------------------------------------------------------------------ %%
\section{Estimates on lower-order functions} % (fold)
\label{sec:estimates_on_lower_order_functions}

% section estimates_on_lower_order_functions (end)

In this subsection, we will establish the estimates on the first and second order expansion functions, $f_1$ and $f_2$. Recalling the decomposition in \cref{eq:decomps-f1} $f_1 = f_1^{\shortparallel} + f_1^{\perp}$, and our truncation assumption $f_2 = f_2^{\perp}$, it suffices to deal with the kernel orthogonal functions $f_i^{\perp}(t,x,v) \in \ker^{\bot}(\L_{\Omega_0})$ with $i=1,2$, and the kernel function $f_1^{\shortparallel} = \rho_1 M_0 + c_0^{-1} (v_{\bot} \cdot {\bm u}_1) M_0 \in \ker(\L_{\Omega_0})$.

%% ------------------------------------------------------------------------ %%
\subsection{Estimates on kernel orthogonal parts} % (fold)
\label{sub:estimates_on_kernel_orthogonal_parts}

% subsection estimates_on_kernel_orthogonal_parts (end)
Note $f_i^{\perp}(t,x,v) \in \ker^{\bot}(\L_{\Omega_0})$ ($i=1,2$) are respectively the unique solution to \cref{eq:1st-order,eq:2nd-order}, due to their solvabilities. Moreover, we can establish the following estimates:
\begin{lemma}\label{lemm:f1-f2-perp-esm}
For the kernel orthogonal part of the lower order functions $f_i^{\perp}(t,x,v)$'s, $i=1,2$, we have
  \begin{align}
    f_1^{\perp} \in C([0,T_0];\, \H^{9}_xL^2_v(M^{-1})), \quad \text{ and }
    f_2^{\perp} \in C([0,T_1];\, \H^{6}_xL^2_v(M^{-1}))\,,
  \end{align}
where $T_0$ and $T_1 \le T_0$ are the maximal existing time of the solutions to the hydrodynamic system \cref{sys:limit-0-order}, and the next order macroscopic system \cref{sys:limit-1st-order-rho,sys:limit-1st-order-j}, respectively.
\end{lemma}

\proof The proof is divided into the following four steps.

\smallskip\noindent\textbf{Step I.}
Firstly, we investigate the estimate for $f_1^{\perp}$, by studying \cref{eq:1st-order}, i.e.,
  \begin{align}
    \L_{\Omega_0} f_1^{\perp} = (\p_t + v \cdot \nabla_x) f_0,
  \end{align}
Taking inner product with $\tfrac{f_1^{\perp} - c_0^{-1} (v_\bot j_1^{\perp})M_0}{M_0}$ with respect to both variables $v$ and $x$, it is straightforward to control the right-hand side by
  \begin{align}
    \sskp{(\p_t + v \cdot \nabla_x) f_0}{\tfrac{f_1^{\perp} - c_0^{-1} (v_\bot j_1^{\perp})M_0}{M_0}}
    \le \nm{(\p_t + v \cdot \nabla_x) f_0}_{\MV} \nm{f_1^{\perp} - c_0^{-1} (v_\bot j_1^{\perp})M_0}_{\MV},
  \end{align}
with the bound
  \begin{align*}
    \nm{f_1^{\perp} - c_0^{-1} (v_\bot j_1^{\perp})M_0}^2_{\MV}
    & \le 2 \nm{f_1^{\perp}}^2_{\MV} + 2 c_0^{-1} \iint |v_\bot|^2 M_0 \d v |j_1^{\perp}|^2 \d x \\
    & \le 2 \nm{f_1^{\perp}}^2_{\MV} + 2 \int \abs{\agl{v f_1^{\perp}}}^2 \d x
    \le C\nm{f_1^{\perp}}^2_{\MV}.
  \end{align*}

We now turn to consider the term involving known function $f_0$. Due to $\p_t f_0 = \p_t (\rho_0 M_0) = \p_t \rho_0 M_0 +  \rho_0 M_0 v_\bot \cdot \p_t \Omega_0$, and $\agl{v_\bot M_0} =0$, we get
  \begin{align}
    \nm{\p_t f_0}^2_{\MV} & = \iint \tfrac{\abs{\p_t f_0}^2}{M_0} \d v \d x \\\no
    & = \iint \abs{\p_t \rho_0}^2 M_0 \d v \d x + \iint \abs{\rho_0 \p_t \Omega_0}^2 \cdot |v_\bot|^2 M_0 \d v \d x
        \\\no & \quad
        + \iint 2 \p_t \rho_0 \p_t \Omega_0 \rho_0 \cdot v_\bot M_0 \d v \d x \\\no
    & \le \nm{\p_t \rho_0}^2_{\lx} + c_0 \nm{\rho_0}^2_{L^\infty_x} \nm{\p_t \Omega_0}^2_{\lx} \,.
    % \le C({\nm{\p_t \rho_0}_{\lx}},\, \nm{\rho_0}_{L^\infty_x},\, \nm{\p_t \Omega_0}_{\lx}) \,.
  \end{align}
In the same way, by $v \cdot \nabla_x f_0 = v \cdot \nabla_x \rho_0 M_0 + \rho_0 M_0 (v \otimes v_{\bot}: \nabla_x \Omega_0)$ we can get
  \begin{align}
    \nm{v \cdot \nabla_x f_0}^2_{\MV}
    & = \iint \abs{\nabla_x \rho_0}^2 |v|^2 M_0 \d v \d x + \iint \abs{\rho_0 \nabla_x \Omega_0}^2 |v\otimes v_{\bot}|^2 M_0 \d v \d x
      \\\no & \quad + \iint 2 \rho_0 \nabla_x \Omega_0 (v \cdot \nabla_x \rho_0) (v \otimes v_{\bot}) M_0 \d v \d x
  \\\no
    & \le C \nm{\nabla_x \rho_0}^2_{\lx} + C \nm{\rho_0}^2_{L^\infty_x} \nm{\nabla_x \Omega_0}^2_{\lx} + C \nm{\rho_0}_{L^\infty_x} \nm{\nabla_x \rho_0}_{\lx} \nm{\nabla_x \Omega_0}_{\lx} \,.
  \end{align}
Combining together the above two inequalities gives
  \begin{align}
    \nm{(\p_t + v \cdot \nabla_x) f_0}_{\MV} \le C_{f_0},
  \end{align}
where the known constant $C_{f_0}$ depends on the quantities $\nm{(\p_t \rho_0,\, \p_t \Omega_0,\, \nabla_x \rho_0,\, \nabla_x \Omega_0)}_{\lx}$ and $\nm{\rho_0}_{L^\infty_x}$, and in turn, depends on $\nm{\rho_0}_{H^1_t(H^1_x \cap L^\infty_x)}$ and $\nm{\Omega_0}_{H^1_t H^1_x}$.

Recalling the coercive estimate \Cref{prop:coercivity-I},
  \begin{align}
    \skp{-\L_{\Omega_0} f_1}{\frac{f_1^{\perp} - c_0^{-1} (v_\bot j_1^{\perp})M_0}{M_0}}
    = \agl{M_0 \abs{\nabla_v \left( \frac{f_1 - c_0^{-1} (v_\bot j_1)M_0}{M_0} \right)}^2 }
    \ge \lambda_0 \eta_0 \agl{\frac{|f_1^{\perp}|^2}{M_0}} \,,
  \end{align}
we can obtain
  \begin{align}
    \lambda_0 \eta_0 \nm{f_1^{\perp}}^2_{\MV}
    \le \nm{ \nabla_v \left( \frac{f_1^{\perp} - c_0^{-1} (v_\bot j_1^{\perp})M_0}{M_0} \right) }^2_{\M}
    \le C_{f_0} \nm{f_1^{\perp}}_{\MV} \,,
  \end{align}
which yields the bounds of
  \begin{align}
    \nm{f_1^{\perp}}_{\MV} \le C_{f_0}, \quad \text{ and }
    \nm{ \nabla_v \left( \frac{f_1^{\perp} - c_0^{-1} (v_\bot j_1^{\perp})M_0}{M_0} \right) }^2_{\M} \le C_{f_0} \,.
  \end{align}
Noticing the basic fact $(a-b)^2 \ge \tfrac{1}{2} a^2 - b^2$, and
  \begin{align*}
    \nabla_v \left( \tfrac{f_1^{\perp} - c_0^{-1} (v_\bot j_1^{\perp})M_0}{M_0} \right) = \nabla_v \left( \tfrac{f_1^{\perp}}{M_0} \right) - c_0^{-1} \left( P_{\Omega_0^\bot} j_1^{\perp} \right),
  \end{align*}
it follows that
  \begin{align}
    \nm{ \nabla_v \left( \frac{f_1^{\perp} - c_0^{-1} (v_\bot j_1^{\perp})M_0}{M_0} \right) }^2_{\M}
    & \ge \tfrac{1}{2} \nm{\nabla_v \left( \frac{f_1^{\perp}}{M_0} \right)}^2_{\M} - \nm{c_0^{-1} \left( P_{\Omega_0^\bot} j_1^{\perp} \right)}^2 \agl{M_0} \\\no
    & \ge \tfrac{1}{2} \nm{\nabla_v \left( \frac{f_1^{\perp}}{M_0} \right)}^2_{\M} - \nm{f_1^{\perp}}^2_{\MV}.
  \end{align}
As a result, we get
  \begin{align}
    \nm{\nabla_v \left( \frac{f_1^{\perp}}{M_0} \right)}_{\M} \le C_{f_0} \,.
  \end{align}
Therefore, we have got $\nm{f_1^{\perp}}_{\MV} + \nm{\nabla_v \left( \frac{f_1^{\perp}}{M_0} \right)}_{\M} \le C_{f_0}$.

\smallskip\noindent\textbf{Step II.}
Next, we consider the higher-order derivatives estimates, by applying higher-order spatial derivative operator $\nabla_x^s$ (with $s\ge 2$) on \cref{eq:1st-order}, i.e.,
  \begin{align}\label{eq:1stOd-HighDeri}
    \nabla_x^s \L_{\Omega_0} f_1
    = (\p_t + v \cdot \nabla_x) \nabla_x^s f_0 \,.
  \end{align}

By using again the expression of $\nabla_x^s \L_{\Omega_0}$ with a commutator,
  \begin{align}
    \nabla_x^s \L_{\Omega_0} g & = \L_{\Omega_0} \nabla_x^s g + [\nabla_x^s,\, \L_{\Omega_0}] g  \\\no
    & = \L_{\Omega_0} \nabla_x^s g - \div_v \left\{ \sum_{{\substack{s_1+s_2 \le s\\s_2 \ge 1}}} \nabla_x^{s_1} g \nabla_x^{s_2} \Omega_0
        + c_0^{-1} \sum_{\mathclap{\substack{s_1+s_2 \le s\\s_2 \ge 1}}} \nabla_x^{s_1} j_g \nabla_x^{s_2} (P_{\Omega_0^\bot} M_0) \right\} \,,
   \end{align}
it follows from taking inner product with $\tfrac{\nabla_x^s f_1^{\perp} - c_0^{-1} (v_\bot \nabla_x^s j_1^{\perp})M_0}{M_0}$, that
  \begin{align}
    \sskp{\L_{\Omega_0} \nabla_x^s f_1^{\perp}}{\tfrac{\nabla_x^s f_1^{\perp} - c_0^{-1} (v_\bot \nabla_x^s j_1^{\perp})M_0}{M_0}} = \nm{ \nabla_v \left( \tfrac{\nabla_x^s f_1^{\perp} - c_0^{-1} (v_\bot \nabla_x^s j_1^{\perp})M_0}{M_0} \right) }^2_{\M} \,,
  \end{align}
and for the commutator, similar as in \cref{ssub:higher-order-esm},
  \begin{align}
    & \sskp{[\nabla_x^s,\, \L_{\Omega_0}] f_1^{\perp}}{\tfrac{\nabla_x^s f_1^{\perp} - c_0^{-1} (v_\bot \nabla_x^s j_1^{\perp})M_0}{M_0}} \\\no
    & \le
    C_{\Omega_0} \sum_{\mathclap{s_1 \le s-1}} \nm{\nabla_x^{s_1} f_1^{\perp}}_{\MV} \nm{ \nabla_v \left( \tfrac{\nabla_x^s f_1^{\perp} - c_0^{-1} (v_\bot \nabla_x^s j_1^{\perp})M_0}{M_0} \right) }_{\M}.
    % \\\no
    % & \le C_{\Omega_0} \nm{\nabla_x^{s_1} f_1^{\perp}}^2_{\MV} + \tfrac{1}{2} \nm{ \nabla_v \left( \tfrac{\nabla_x^s f_1^{\perp} - c_0^{-1} (v_\bot \nabla_x^s j_1^{\perp})M_0}{M_0} \right) }^2_{\M}.
  \end{align}
Here the coefficient $C_{\Omega_0} $ depends on the bound of $\nm{\nabla_x^{s_2} \Omega_0}_{L^\infty_x}$ ($s_2 \le s$), as justified before.

Combining the above two inequalities and the H\"older inequality gives
  \begin{align}\label{eq:1stOd-lhs}
    & \sskp{-\nabla_x^s \L_{\Omega_0} f_1^{\perp}}{\tfrac{\nabla_x^s f_1^{\perp} - c_0^{-1} (v_\bot \nabla_x^s j_1^{\perp})M_0}{M_0}} \\\no
    &\ \ge \nm{ \nabla_v \left( \tfrac{\nabla_x^s f_1^{\perp} - c_0^{-1} (v_\bot \nabla_x^s j_1^{\perp})M_0}{M_0} \right) }^2_{\M}
        - C_{\Omega_0} \sum_{\mathclap{s_1 \le s-1}} \nm{\nabla_x^{s_1} f_1^{\perp}}_{\MV} \nm{ \nabla_v \left( \tfrac{\nabla_x^s f_1^{\perp} - c_0^{-1} (v_\bot \nabla_x^s j_1^{\perp})M_0}{M_0} \right) }_{\M} \\\no
    &\ \ge \tfrac{3}{4} \nm{ \nabla_v \left( \tfrac{\nabla_x^s f_1^{\perp} - c_0^{-1} (v_\bot \nabla_x^s j_1^{\perp})M_0}{M_0} \right) }^2_{\M}
            - C_{\Omega_0} \sum_{\mathclap{s_1 \le s-1}} \nm{\nabla_x^{s_1} f_1^{\perp}}^2_{\MV} \,.
  \end{align}

On the other hand, the right-hand side of \cref{eq:1stOd-HighDeri} can be controlled by
  \begin{align}\label{eq:1stOd-rhs}
    & \sskp{(\p_t + v \cdot \nabla_x) \nabla_x^s f_0}{\tfrac{\nabla_x^s f_1^{\perp} - c_0^{-1} (v_\bot \nabla_x^s j_1^{\perp})M_0}{M_0}} \\\no
    & \le \nm{(\p_t + v \cdot \nabla_x) \nabla_x^s f_0}_{\MV} \nm{\nabla_x^s f_1^{\perp} - c_0^{-1} (v_\bot \nabla_x^s j_1^{\perp})M_0}_{\MV} \\\no
    & \le \nm{(\p_t + v \cdot \nabla_x) \nabla_x^s f_0}_{\MV} \nm{\nabla_x^s f_1^{\perp}}_{\MV} \,,
  \end{align}
where we have used $\abs{\nabla_x^s j_1^{\perp}} = \abs{\agl{v \nabla_x^s f_1^{\perp}}} \le \nm{\nabla_x^s f_1^{\perp}}_{\MV}$, by virtue of the H\"older inequality.

Arguing as above, we can get from the fact $f_0 = \rho_0 M_0$, that
  \begin{align}
    \nm{(\p_t + v \cdot \nabla_x) \nabla_x^s f_0}_{\MV} \le C_{s,f_0},
  \end{align}
with the known constant $C_{s,f_0}$ depending on the quantities $\nm{(\p_t \rho_0,\, \p_t \Omega_0,\, \nabla_x \rho_0,\, \nabla_x \Omega_0)}_{\hx{s}}$. This fact, combined together with \cref{eq:1stOd-lhs,eq:1stOd-rhs}, enables us to get
  \begin{align}\label{eq:cycle-bdd}
    \lambda_0 \eta_0 \nm{\nabla_x^s f_1^{\perp}}^2_{\MV}
    & \le \nm{ \nabla_v \left( \tfrac{\nabla_x^s f_1^{\perp} - c_0^{-1} (v_\bot \nabla_x^s j_1^{\perp})M_0}{M_0} \right) }^2_{\M} \\\no
    & \le C_{s,f_0} \nm{\nabla_x^s f_1^{\perp}}_{\MV} + C_{\Omega_0} \sum_{\mathclap{s_1 \le s-1}} \nm{\nabla_x^{s_1} f_1^{\perp}}^2_{\MV} \,,
    % \\\no
    % & \le \tfrac{1}{4} \lambda_0 \eta_0 \nm{\nabla_x^s f_1^{\perp}}_{\MV} + C_{s,f_0} + C_{\Omega_0} \sum_{\mathclap{s_1 \le s-1}} \nm{\nabla_x^{s_1} f_1^{\perp}}^2_{\MV} \,,
  \end{align}
where we have used again the coercive estimate \Cref{prop:coercivity-I}.

The H\"older inequality, $C_{s,f_0} \nm{\nabla_x^s f_1^{\perp}}_{\MV} \le \tfrac{1}{4} \lambda_0 \eta_0 \nm{\nabla_x^s f_1^{\perp}}^2_{\MV} + C_{s,f_0}$, yields immediately that
  \begin{align}
    \nm{\nabla_x^s f_1^{\perp}}^2_{\MV}
    \le C_{s,f_0} + C_{\Omega_0} \sum_{\mathclap{s_1 \le s-1}} \nm{\nabla_x^{s_1} f_1^{\perp}}^2_{\MV} \,.
  \end{align}
A simple induction principle will lead us to the bound
  \begin{align}
    \nm{\nabla_x^s f_1^{\perp}}^2_{\MV}
    \le C_{s,f_0,\Omega_0}.
  \end{align}
We mention the constant $C_{s,f_0,\Omega_0}$ depends on the bounds of $\nm{(\p_t \rho_0,\, \p_t \Omega_0,\, \nabla_x \rho_0,\, \nabla_x \Omega_0)}_{\hx{s}}$ and $\nm{\nabla_x^s \Omega_0}_{L^\infty_x}$. Then by the dynamical system obeyed by $(\rho_0,\, \Omega_0)$ and the Sobolev embedding inequality, $C_{s,f_0,\Omega_0}$ finally depends on $\nm{(\rho_0,\, \Omega_0)}_{\hx{s+2}}$.

Inserting the above bound for $\nabla_x^s f_1^{\perp}$ into \cref{eq:cycle-bdd}, we also infer
  \begin{align}
    \nm{ \nabla_v \left( \tfrac{\nabla_x^s f_1^{\perp} - c_0^{-1} (v_\bot \nabla_x^s j_1^{\perp})M_0}{M_0} \right) }^2_{\M} \le C_{s,f_0,\Omega_0}.
  \end{align}
Therefore, we get higher-order derivatives estimates
  \begin{align*}
    \sum_{k \le s} \left[ \nm{\nabla_x^k f_1^{\perp}}_{\MV} + \nm{\nabla_v \left( \tfrac{\nabla_x^k f_1^{\perp}}{M_0} \right)}_{\M} \right] \le C_{s,f_0,\Omega_0} \,,
  \end{align*}
i.e.,
  \begin{align}
    \nm{f_1^{\perp}}_{\hxv{s}(\MV)} + \nm{\nabla_v \left( \tfrac{f_1^{\perp}}{M_0} \right)}_{\hxv{s}(\M)} \le C_{s,f_0,\Omega_0} \,,
  \end{align}

\smallskip\noindent\textbf{Step III.}
In the third step, we need to estimate the quantity $\p_t f_1^{\perp}$ and its higher derivatives $\nabla_x^s \p_t f_1^{\perp}$. Noticing the formula that
  \begin{align}
    \p_t \L_{\Omega_0} g & = \L_{\Omega_0} \p_t g + [\p_t,\, \L_{\Omega_0}] g  \\\no
    & = \L_{\Omega_0} \p_t g - \div_v \left\{ g \p_t \Omega_0 + c_0^{-1} j_g \p_t (P_{\Omega_0^\bot} M_0) \right\} \,,
   \end{align}
and
  \begin{align}
    \nabla_x^s \p_t \L_{\Omega_0} g & = \L_{\Omega_0} \nabla_x^s \p_t g + [\nabla_x^s \p_t,\, \L_{\Omega_0}] g  \\\no
    & = \L_{\Omega_0} \nabla_x^s \p_t g
        - \div_v \left\{ \sum_{{\substack{s_1+s_2 \le s\\s_2 \ge 1}}} \nabla_x^{s_1} \p_t g \nabla_x^{s_2} \Omega_0
                          + c_0^{-1} \sum_{\mathclap{\substack{s_1+s_2 \le s\\s_2 \ge 1}}} \nabla_x^{s_1} \p_t j_g \nabla_x^{s_2} (P_{\Omega_0^\bot} M_0) \right\}
        \\\no & \hspace*{2.2cm}
        - \div_v \left\{ \sum_{{\substack{s_1+s_2 \le s}}} \nabla_x^{s_1} g \nabla_x^{s_2} \p_t \Omega_0
                          + c_0^{-1} \sum_{\mathclap{\substack{s_1+s_2 \le s}}} \nabla_x^{s_1} j_g \nabla_x^{s_2} \p_t (P_{\Omega_0^\bot} M_0) \right\} \,,
  \end{align}
we can perform the similar scheme as above to obtain the estimates:
  \begin{align}
    % \sum_{k \le s} \left[ \nm{\nabla_x^k \p_t f_1^{\perp}}_{\MV} + \nm{\nabla_v \left( \tfrac{\nabla_x^k \p_t f_1^{\perp}}{M_0} \right)}_{\M} \right] \le C_{t,s,f_0,\Omega_0},
    \nm{\p_t f_1^{\perp}}_{\hxv{s}(\MV)} + \nm{\nabla_v \left( \tfrac{\p_t f_1^{\perp}}{M_0} \right)}_{\hxv{s}(\M)} \le C_{t,s,f_0,\Omega_0},
  \end{align}
with the bound $C_{t,s,f_0,\Omega_0}$ depending upon $\nm{(\p_t + v \cdot \nabla_x) \nabla_x^s \p_t f_0}_{\MV}$ and $\nm{\nabla_x^s \p_t \Omega_0}_{L^\infty_x}$, those of which depend finally on $\nm{(\rho_0,\, \Omega_0)}_{\hx{s+3}}$.

Recalling the local solution $(\rho_0,\, \Omega_0)$ to the hydrodynamic system of \cref{sys:limit-0-order} is proven to be well-posed in $C([0,T_0];\, H^{12}_x(\mathbb{R}^3))$, see the main \Cref{thm:limit}, we can obtain
  \begin{align}
    f_1^{\perp} \in H^1([0,T_0];\, \H^{9}_xL^2_v(M^{-1})).
  \end{align}

\smallskip\noindent\textbf{Step IV.}
The estimate on $f_2^{\perp}$ follows along the same lines of that on $f_1^{\perp}$, by considering \cref{eq:2nd-order}. We will obtain
  \begin{align}
    % \sum_{k \le s} \left[ \nm{\nabla_x^k \p_t f_1^{\perp}}_{\MV} + \nm{\nabla_v \left( \tfrac{\nabla_x^k \p_t f_1^{\perp}}{M_0} \right)}_{\M} \right] \le C_{t,s,f_0,\Omega_0},
    \nm{(f_2^{\perp},\, \p_t f_2^{\perp})}_{\hxv{s}(\MV)} + \nm{\nabla_v \left( \tfrac{(f_2^{\perp},\, \p_t f_2^{\perp})}{M_0} \right)}_{\hxv{s}(\M)} \le C_{s,f_1},
  \end{align}
with the constant $C_{s,f_1}$ depending upon $\nm{(\p_t + v \cdot \nabla_x) \nabla_x^s \p_t f_1}_{\MV}$ and $\nm{\nabla_x^s \p_t j_1}_{L^\infty_x}$. This in turn yields a dependence of $\nm{f_1}_{H^1_t\H^{s+3}_xL^2_v(M^{-1})}$, and finally of the bounds $\nm{(\rho_0,\, \Omega_0)}_{\hx{s+6}}$. Therefore, for known $(\rho_0,\, \Omega_0) \in C([0,T_0];\, H^{12}_x(\mathbb{R}^3))$, we can justify the desired result $f_2^{\perp} \in H^1([0,T_1];\, \H^{6}_xL^2_v(M^{-1}))$. This completes the whole proof of \Cref{lemm:f1-f2-perp-esm}.
\endproof

%% ------------------------------------------------------------------------ %%
\subsection{Estimates on next-order macroscopic equation} % (fold)
\label{sub:estimates_on_next_order_macroscopic_equation}

% subsection estimates_on_next_order_macroscopic_equation (end)

We mention that, when considering the second order function $f_2^{\perp}$, all of information for $f_1$ should be viewed as known terms, including both kernel orthogonal part $f_1^{\perp}$ and kernel part $f_1^{\shortparallel}$. As for the latter, since $f_1^{\shortparallel}(t,x,v) = \rho_1 M_0 + c_0^{-1} (v_{\bot} \cdot {\bm u}_1) M_0 \in \ker(\L_{\Omega_0})$, we can obtain the following lemma.
\begin{lemma}\label{lemm:system-f1-esm}
  Let $m \ge 2$. Given initial data $(\rho_1^{\IN}, \Omega_1^{\IN}) \in H^m_x(\mathbb{R}^3)$, then there exists a local time  $T_1 \in (0, T_0]$, such that $(\rho_1, \Omega_1) \in C([0,T_1];$ $H^m_x(\mathbb{R}^3))$ is a local solution to the Cauchy problem of the first order macroscopic system \cref{sys:limit-1st-order-rho,sys:limit-1st-order-j}.
\end{lemma}

\smallskip\noindent\emph{Sketch Proof of \Cref{lemm:system-f1-esm}}.
Noticing a local existence result is sufficient here, the proof for \Cref{lemm:system-f1-esm} is straightforward. The key ingredients are to treat carefully the possible gained derivatives, and our treatment is to transfer one derivative on the studied quantities $\rho_1,\, \u_1$ to some other known functions $\rho_0,\, \Omega_0$, by taking advantage of integration by parts, or by the orthogonal facts $\Omega_0 \cdot \u_1 = 0$. Here we only show the closed energy estimate to give the lifespan of existence. For simplicity, we consider the case $m=2$.

Applying the operator $\nabla_x^2$ on \cref{sys:limit-1st-order-rho,sys:limit-1st-order-j}, we get
  \begin{align}\label{sys:limit-1st-order-rho-HD}
    \p_t \nabla_x^2 \rho_1 + \div_x \nabla_x^2 (c_0 \Omega_0 \rho_1) + \div_x \nabla_x^2 \u_1 = \nabla_x^2 S_{\rho}(U_{\text{source}}(f_1^{\perp})) \,,
  \end{align}
and
    \begin{align}\label{sys:limit-1st-order-j-HD}
      \p_t \nabla_x^2 \u_1^i & + c_1 \nabla_x^2 [(\Omega_0 \cdot \nabla_x) \u_1^i] \\[5pt] \no
        & + \nabla_x^2 [ (\u_1 \cdot \nabla_x) \ln \rho_0 \cdot \Omega_0^i + c_2 (\u_1 \cdot \nabla_x) \Omega_0^i
      % \\[5pt] \no & \quad
      + (c_2 - c_1) P_{\Omega_0^\bot}^{il} \nabla_{x^l} \Omega_0^k \u_1^k
        + c_2 \div_{x} \Omega_0 \cdot \u_1^i]
    \\[5pt] \no
        & + c_0 \nabla_x^2 [P_{\Omega_0^\bot}^{il} \nabla_{x^l} \rho_1 - \rho_1 P_{\Omega_0^\bot}^{il} \nabla_{x^l} \ln \rho_0]
        + c_0 \tilde{c}_0^{-1} \nabla_x^2 (\rho_0^{-1} \rho_1 \u_1^i)
        \\[5pt] \no
      = & - \nabla_x^2 \left\{ (\tilde{c}_0 \rho_0)^{-1} \left[ c_{f_1^{\perp}}^{il} \u_1^l + 2(\Omega_0 j_1^{\perp}) \u_1^i \right] \right\}
        - c_0 \tilde{c}_0^{-1} \nabla_x^2 \left[ \rho_0^{-1} \rho_1 (P_{\Omega_0^\bot} j_1^{\perp})^i \right]
        - \nabla_x^2 S_{\u}(U_{\text{source}}(f_1^{\perp})) \,.
    \end{align}

In the following we consider the inner product in $L^2_x$ of the two equations with $c_0 \nabla_x^2 \rho_1$ and $\nabla_x^2 \u_1^i$, respectively. Firstly we have
  \begin{align}
    \skp{\p_t \nabla_x^2 \rho_1}{c_0 \nabla_x^2 \rho_1}_{L^2_x}
    & = \frac{1}{2} \frac{\d }{\d t} c_0 \nm{\nabla_x^2 \rho_1}^2_{\lx}, \\
    \skp{\p_t \nabla_x^2 \u_1^i}{\nabla_x^2 \u_1^i}_{L^2_x}
    & = \frac{1}{2} \frac{\d }{\d t} \nm{\nabla_x^2 \u_1}^2_{\lx}.
  \end{align}

As for the second term in \cref{sys:limit-1st-order-rho-HD}, it follows
  \begin{align}
    \skp{\div_x \nabla_x^2 (c_0 \Omega_0 \rho_1)}{c_0 \nabla_x^2 \rho_1}_{L^2_x}
    = c_0^2 \skp{\nabla_x^2 (\Omega_0 \nabla_x \rho_1)}{\nabla_x^2 \rho_1}_{L^2_x}
      + c_0^2 \skp{\nabla_x^2 (\div_x \Omega_0 \rho_1)}{\nabla_x^2 \rho_1}_{L^2_x},
  \end{align}
in which the first term can be controlled, by an integration by parts and the Sobolev embedding inequality, by
  \begin{align*}
    & \skp{\nabla_x^2 (\Omega_0 \nabla_x \rho_1)}{\nabla_x^2 \rho_1}_{L^2_x} \\\no
    & = \skp{\Omega_0 \nabla_x (\nabla_x^2 \rho_1)}{\nabla_x^2 \rho_1}_{L^2_x}
      + \skp{\nabla_x \Omega_0 \nabla_x^2 \rho_1}{\nabla_x^2 \rho_1}_{L^2_x}
      + \skp{\nabla_x^2 \Omega_0 \nabla_x^2 \rho_1}{\nabla_x^2 \rho_1}_{L^2_x} \\\no
    & = - \int \div_x \Omega_0 \abs{\nabla_x^2 \rho_1}^2 \d x
        + \skp{\nabla_x \Omega_0 \nabla_x^2 \rho_1}{\nabla_x^2 \rho_1}_{L^2_x}
        + \skp{\nabla_x^2 \Omega_0 \nabla_x \rho_1}{\nabla_x^2 \rho_1}_{L^2_x} \\\no
    & \ls \nm{\div_x \Omega_0}_{L^\infty_x} \nm{\nabla_x^2 \rho_1}^2_{\lx}
          + \nm{\nabla_x \Omega_0}_{L^\infty_x} \nm{\nabla_x^2 \rho_1}^2_{\lx}
          + \nm{\nabla_x^2 \Omega_0}_{L^4_x} \nm{\nabla_x \rho_1}_{L^4_x} \nm{\nabla_x^2 \rho_1}_{\lx}
    \\\no
    & \ls \nm{\nabla_x \Omega_0}_{\hx{2}} \nm{\rho_1}_{\hx{2}} \nm{\nabla_x^2 \rho_1}_{\lx} \,,
  \end{align*}
and the second term can be treated similarly,
  \begin{align*}
    & \skp{\nabla_x^2 (\div_x \Omega_0 \rho_1)}{\nabla_x^2 \rho_1}_{L^2_x} \\\no
    & = \skp{\div_x \Omega_0 \nabla_x^2 \rho_1}{\nabla_x^2 \rho_1}_{L^2_x}
      + \skp{\nabla_x \div_x \Omega_0 \nabla_x \rho_1}{\nabla_x^2 \rho_1}_{L^2_x}
      + \skp{\nabla_x^2 \div_x \Omega_0 \rho_1}{\nabla_x^2 \rho_1}_{L^2_x} \\\no
    & \ls \nm{\div_x \Omega_0}_{L^\infty_x} \nm{\nabla_x^2 \rho_1}^2_{\lx}
          + \nm{\nabla_x \div_x \Omega_0}_{L^4_x} \nm{\nabla_x \rho_1}_{L^4_x} \nm{\nabla_x^2 \rho_1}_{\lx}
          \\\no & \quad
          + \nm{\nabla_x^2 \div_x \Omega_0}_{L^2_x} \nm{\nabla_x \rho_1}_{L^\infty_x} \nm{\nabla_x^2 \rho_1}_{\lx}
    \\\no
    & \ls \nm{\nabla_x \Omega_0}_{\hx{2}} \nm{\rho_1}_{\hx{2}} \nm{\nabla_x^2 \rho_1}_{\lx} \,.
  \end{align*}
Combining the two inequalities yields
  \begin{align}
    \skp{\div_x \nabla_x^2 (c_0 \Omega_0 \rho_1)}{c_0 \nabla_x^2 \rho_1}_{L^2_x}
    \ls c_0^2 \nm{\nabla_x \Omega_0}_{\hx{2}} \nm{\rho_1}_{\hx{2}} \nm{\nabla_x^2 \rho_1}_{\lx} \,.
  \end{align}

Concerning the second term in \cref{sys:limit-1st-order-j-HD}, by using again the Leibniz formula and the integration by parts, we get, (by omitting the subscripts $L^2_x$ for brevity)
  \begin{align}
    & \skp{\nabla_x^2 [(\Omega_0 \nabla_x) \u_1^i]}{\nabla_x^2 \u_1^i} \\\no
    & = \skp{\Omega_0 \nabla_x (\nabla_x^2 \u_1)}{\nabla_x^2 \u_1^i}
      + \skp{\nabla_x \Omega_0 \nabla_x^2 \u_1}{\nabla_x^2 \u_1^i}
      + \skp{\nabla_x^2 \Omega_0 \nabla_x^2 \u_1}{\nabla_x^2 \u_1^i} \\\no
    % & = - \int \div_x \Omega_0 \abs{\nabla_x^2 \u_1}^2 \d x
    %     + \skp{\nabla_x \Omega_0 \nabla_x^2 \u_1}{\nabla_x^2 \u_1^i}
    %     + \skp{\nabla_x^2 \Omega_0 \nabla_x \u_1}{\nabla_x^2 \u_1^i} \\\no
    % & \ls \nm{\div_x \Omega_0}_{L^\infty_x} \nm{\nabla_x^2 \rho_1}^2_{\lx}
    %       + \nm{\nabla_x \Omega_0}_{L^\infty_x} \nm{\nabla_x^2 \rho_1}^2_{\lx}
    %       + \nm{\nabla_x^2 \Omega_0}_{L^4_x} \nm{\nabla_x \rho_1}_{L^4_x} \nm{\nabla_x^2 \rho_1}_{\lx}
    % \\\no
    & \ls \nm{\nabla_x \Omega_0}_{\hx{2}} \nm{\u_1}_{\hx{2}} \nm{\nabla_x^2 \u_1}_{\lx} \,.
  \end{align}

The third term on the left-hand side of \cref{sys:limit-1st-order-rho-HD} should be considered together with the pressure term in the third line of \cref{sys:limit-1st-order-j-HD}, due to a cancellation relation between them. Precisely speaking, we have,
  \begin{align}
    & \skp{\div_x \nabla_x^2 \u_1}{c_0 \nabla_x^2 \rho_1} + \skp{c_0 \nabla_x^2 P_{\Omega_0^\bot}^{il} \nabla_{x^l} \rho_1}{\nabla_x^2 \u_1^i} \\\no
    & = \skp{\div_x \nabla_x^2 \u_1}{c_0 \nabla_x^2 \rho_1} + \skp{c_0 \nabla_x^2 \nabla_{x^i} \rho_1}{\nabla_x^2 \u_1^i} - \skp{c_0 \nabla_x^2 (\Omega_0^i \Omega_0^l \nabla_{x^l} \rho_1)}{\nabla_x^2 \u_1^i} \\\no
    & = - \skp{c_0 \nabla_x^2 (\Omega_0^i \Omega_0^l \nabla_{x^l} \rho_1)}{\nabla_x^2 \u_1^i} \\\no
    & = - c_0 \skp{\nabla_x^2 (\Omega_0^l \nabla_{x^l} \rho_1)}{\Omega_0^i \nabla_x^2 \u_1^i}
        - c_0 \skp{\nabla_x (\Omega_0^l \nabla_{x^l} \rho_1) \nabla_x \Omega_0^i }{\nabla_x^2 \u_1^i} \\\no
        & \quad
        - c_0 \skp{\Omega_0^l \nabla_{x^l} \rho_1 \nabla_x^2 \Omega_0^i }{\nabla_x^2 \u_1^i} \,.
  \end{align}
Note a factor with one more derivative $\nabla_x^3 \rho_1$ is possible to arise, which needs to be treated carefully. Fortunately, by observing the fact $\Omega_0 \cdot \u_1 =0$, and hence it follows $ \nabla_x^2 (\Omega_0^i \u_1^i ) = \Omega_0^i \nabla_x^2 \u_1^i + \nabla_x \u_1^i \nabla_x \Omega_0^i + \nabla_x \Omega_0^i \nabla_x \u_1^i + \nabla_x^2 \Omega_0^i \u_1^i =0$. Thus, an integration by parts enables us to get
  \begin{align}
    & \skp{\nabla_x^2 (\Omega_0^l \nabla_{x^l} \rho_1)}{\Omega_0^i \nabla_x^2 \u_1^i} \\\no
    & = - \skp{\nabla_x^2 (\Omega_0^l \nabla_{x^l} \rho_1)}
              { \nabla_x \u_1^i \nabla_x \Omega_0^i + \nabla_x \Omega_0^i \nabla_x \u_1^i + \nabla_x^2 \Omega_0^i \u_1^i } \\\no
    & = \skp{\nabla_x (\Omega_0^l \nabla_{x^l} \rho_1)}
            {\div_x \left( \nabla_x \u_1^i \nabla_x \Omega_0^i + \nabla_x \Omega_0^i \nabla_x \u_1^i + \nabla_x^2 \Omega_0^i \u_1^i \right) } \\\no
    & = \skp{\Omega_0^l \nabla_{x^l} \nabla_x \rho_1}{\Delta_x \u_1^i \nabla_x \Omega_0^i + \nabla_x \Omega_0^i \nabla_x^2 \u_1^i} % \\\no & \quad
        + \skp{\nabla_x \Omega_0^l \nabla_{x^l} \rho_1}{\Delta_x \u_1^i \nabla_x \Omega_0^i + \nabla_x \Omega_0^i \nabla_x^2 \u_1^i} \\\no & \quad
        + \skp{\Omega_0^l \nabla_{x^l} \nabla_x \rho_1 + \nabla_x \Omega_0^l \nabla_{x^l} \rho_1}
              {\nabla_x \u_1^i \nabla_x^2 \Omega_0^i + \Delta_x \Omega_0^i \nabla_x \u_1^i + \nabla_x \Delta_x \Omega_0^i \u_1^i + \nabla_x^2 \Omega_0^i \nabla_x \u_1^i} \\\no
    & \le \nm{\Omega_0}_{L^\infty_x} \nm{\nabla_x \Omega_0}_{L^\infty_x} \nm{\nabla_x^2 \rho_1}_{\lx}
              \left( \nm{\Delta_x \u_1}_{\lx} + \nm{\nabla_x^2 \u_1}_{\lx} \right) \\\no & \quad
          + \nm{\nabla_x \Omega_0}^2_{L^\infty_x} \nm{\nabla_x \rho_1}_{\lx}
                \left( \nm{\Delta_x \u_1}_{\lx} + \nm{\nabla_x^2 \u_1}_{\lx} \right)
    \\\no & \quad
          + \left( \nm{\Omega_0}_{L^\infty_x} \nm{\nabla_x^2 \rho_1}_{\lx} + \nm{\nabla_x \Omega_0}_{L^\infty_x} \nm{\nabla_x \rho_1}_{\lx} \right)
            \cdot \Big( \nm{\nabla_x \u_1}_{L^4_x} \nm{\nabla_x^2 \Omega_0}_{L^4_x} \\\no & \hspace*{2cm}
                        + \nm{\nabla_x \u_1}_{L^4_x} \nm{\Delta_x \Omega_0}_{L^4_x}
                        + \nm{\nabla_x \Delta_x \Omega_0}_{L^2_x} \nm{\u_1}_{L^\infty_x}
                        + \nm{\nabla_x^2 \Omega_0}_{L^4_x} \nm{\nabla_x \u_1}_{L^4_x} \Big)
    \\\no
    & \le \nm{\nabla_x \Omega_0}^2_{\hx{2}} \nm{\rho_1}_{\hx{2}} \nm{\u_1}_{\hx{2}} \,.
  \end{align}
Combined together with the facts
  \begin{align*}
    \skp{\Omega_0^l \nabla_{x^l} \rho_1 \nabla_x^2 \Omega_0^i }{\nabla_x^2 \u_1^i}
    \le \nm{\nabla_x^2 \Omega_0}_{L^4_x} \nm{\nabla_x \rho_1}_{L^4_x} \nm{\nabla_x^2 \u_1}_{\lx}
    \ls \nm{\nabla_x \Omega_0}_{\hx{2}} \nm{\rho_1}_{\hx{2}} \nm{\nabla_x^2 \u_1}_{\lx} \,,
  \end{align*}
and
  \begin{align*}
    \skp{\nabla_x (\Omega_0^l \nabla_{x^l} \rho_1) \nabla_x \Omega_0^i }{\nabla_x^2 \u_1^i}
    & = \skp{(\nabla_x \Omega_0^l \nabla_{x^l} \rho_1 + \Omega_0^l \nabla_{x^l} \nabla_x \rho_1) \nabla_x \Omega_0^i }{\nabla_x^2 \u_1^i} \\\no
    & \ls \nm{\nabla_x \Omega_0}^2_{L^\infty_x} \left( \nm{\nabla_x \rho_1}_{\lx} + \nm{\nabla_x^2 \rho_1}_{\lx} \right) \nm{\nabla_x^2 \u_1}_{\lx} \\\no
    & \ls \nm{\nabla_x \Omega_0}^2_{\hx{2}} \nm{\rho_1}_{\hx{2}} \nm{\nabla_x^2 \u_1}_{\lx} \,,
  \end{align*}
we can get
  \begin{align}
    \skp{\div_x \nabla_x^2 \u_1}{c_0 \nabla_x^2 \rho_1} + \skp{c_0 \nabla_x^2 P_{\Omega_0^\bot}^{il} \nabla_{x^l} \rho_1}{\nabla_x^2 \u_1^i}
    & = - \skp{c_0 \nabla_x^2 (\Omega_0^i \Omega_0^l \nabla_{x^l} \rho_1)}{\nabla_x^2 \u_1^i} \\\no
    & \ls c_0 \nm{\nabla_x \Omega_0}^2_{\hx{2}} \nm{\rho_1}_{\hx{2}} \nm{\u_1}_{\hx{2}} \,.
  \end{align}

Next, we deal with the linear terms in the second line of \cref{sys:limit-1st-order-j-HD}, with a linear dependence on $\u_1$. Straightforward calculations ensure
  \begin{align}
    & \skp{\nabla_x^2 [ (\u_1 \cdot \nabla_x) \ln \rho_0 \cdot \Omega_0^i
                      + c_2 (\u_1 \cdot \nabla_x) \Omega_0^i
                      + (c_2 - c_1) P_{\Omega_0^\bot}^{il} \nabla_{x^l} \Omega_0^k \u_1^k
                      + c_2 \div_{x} \Omega_0 \cdot \u_1^i]}
        {\nabla_x^2 \u_1^i}
    \no \\[5pt] \no
    & \ls \nm{\nabla_x \rho_0}_{\hx{2}} \nm{\nabla_x \Omega_0}_{\hx{2}} \nm{\u_1}_{\hx{2}} \nm{\nabla_x^2 \u_1}_{\lx} \,.
  \end{align}
Similarly, it holds
  \begin{align}
    \skp{\nabla_x^2 [ \rho_1 P_{\Omega_0^\bot}^{il} \nabla_{x^l} \ln \rho_0 ]}{\nabla_x^2 \u_1^i}
    % \no \\[5pt] \no &
    \ls \nm{\rho_0}_{\hx{3}} \nm{\Omega_0}_{\hx{2}} \nm{\rho_1}_{\hx{2}} \nm{\nabla_x^2 \u_1}_{\lx} \,.
  \end{align}

Now we are left to consider the last term in the left-hand side with a nonlinearity of $\rho_1 \u_1$, and the linear terms on the right-hand side with their coefficients depending on the known part $f_1^{\perp}$.
As for the nonlinear product term, by the fact the Sobolev space $H^s_x$ ($s \ge 2$) is an algebra, we can infer that
  \begin{align}
    \skp{\nabla_x^2 (\rho_0^{-1} \rho_1 \u_1^i)}{\nabla_x^2 \u_1^i}
    \ls \nm{\rho_0}_{\hx{2}} \nm{\rho_1}_{\hx{2}} \nm{\u_1}_{\hx{2}} \nm{\nabla_x^2 \u_1}_{\lx} \,.
  \end{align}
Recalling the notation $ c_{f_1^{\perp}}^{il}(t,x) = \skpt{ \nabla_{v^l} f_1^{\perp} }{ \chi(v) \tfrac{v_{\bot}^i}{\abs{v_{\bot}}} }$, the right-hand side can be controlled by
  \begin{align}
    & \skp{\nabla_x^2 \left\{ (\tilde{c}_0 \rho_0)^{-1} \left[ c_{f_1^{\perp}}^{il} \u_1^l + 2(\Omega_0 j_1^{\perp}) \u_1^i \right] \right\} }{\nabla_x^2 \u_1^i}_{L^2_x} \\\no
    & \ls \nm{\rho_0}_{\hx{2}} \left( \nm{c_{f_1^{\perp}}}_{\hx{2}} + \nm{\Omega_0}_{\hx{2}} \nm{j_1^{\perp}}_{\hx{2}} \right) \nm{\u_1}_{\hx{2}} \nm{\nabla_x^2 \u_1}_{\lx}
    \\\no
    & \ls \nm{\rho_0}_{\hx{2}} \nm{\Omega_0}_{\hx{2}} \nm{ f_1^{\perp} }_{\hxv{2}(\MV)} \nm{\u_1}_{\hx{2}} \nm{\nabla_x^2 \u_1}_{\lx} \,.
  \end{align}
and similarly,
  \begin{align}
    \skp{\nabla_x^2 \left[ \rho_0^{-1} \rho_1 (P_{\Omega_0^\bot} j_1^{\perp})^i \right] }{\nabla_x^2 \u_1^i}_{L^2_x}
    \ls \nm{\rho_0}_{\hx{2}} \nm{\Omega_0}_{\hx{2}} \nm{ f_1^{\perp} }_{\hxv{2}(\MV)} \nm{\rho_1}_{\hx{2}} \nm{\nabla_x^2 \u_1}_{\lx} \,.
  \end{align}

Finally, taking summation of all the above estimates, we can obtain the energy estimates on $(\rho_1,\, \u_1)$,
  \begin{align}
    & \frac{1}{2} \frac{\d}{\d t} \left( \nm{\rho_1}^2_{\hx{2}} + \nm{\u_1}^2_{\hx{2}} \right) \\\no
    & \ls \nm{\nabla_x \Omega_0}_{\hx{2}} \left( \nm{\rho_1}^2_{\hx{2}} + \nm{\u_1}^2_{\hx{2}} \right)
        + \nm{\nabla_x \Omega_0}^2_{\hx{2}} \nm{\rho_1}_{\hx{2}} \nm{\u_1}_{\hx{2}}
        \\\no & \quad + \nm{\nabla_x \rho_0}_{\hx{2}} \nm{\nabla_x \Omega_0}_{\hx{2}} \nm{\u_1}^2_{\hx{2}}
        + \nm{\rho_0}_{\hx{3}} \nm{\Omega_0}_{\hx{2}} \nm{\rho_1}_{\hx{2}} \nm{\u_1}_{\hx{2}}
        \\\no & \quad + \nm{\rho_0}_{\hx{2}} \nm{\rho_1}_{\hx{2}} \nm{\u_1}^2_{\hx{2}}
        + \nm{\rho_0}_{\hx{2}} \nm{\Omega_0}_{\hx{2}} \nm{ f_1^{\perp} }_{\hxv{2}(\MV)} \left( \nm{\u_1}^2_{\hx{2}} + \nm{\rho_1}_{\hx{2}} \nm{\u_1}_{\hx{2}} \right) \,,
  \end{align}
By denoting $\EE = \nm{\rho_1}^2_{\hx{2}} + \nm{\u_1}^2_{\hx{2}} $, this inequality can be reformulated as
  \begin{align}
    \frac{1}{2} \frac{\d}{\d t} \EE (t)
    \le C_{f_0,f_1^{\perp}} (\EE + \EE^\frac{3}{2}).
  \end{align}
Applying the Gr\"onwall lemma ensures
  \begin{align}
    \EE (t) \le \frac{\EE^\frac{1}{2}(0)}{1- t C_{f_0,f_1^{\perp}} \EE^\frac{1}{2}(0)} \exp \left( { t C_{f_0,f_1^{\perp}} } \right) ,
  \end{align}
which, together with a standard compactness argument, leads to the existence of a local time $T_1 < \min \left\{ T_0,\, C_{f_0,f_1^{\perp}}^{-1} \EE^{-\frac{1}{2}}(0) \right\} $, such that $(\rho_1, \Omega_1) \in C([0,T_1];$ $H^2_x(\mathbb{R}^3))$ is a local solution to the Cauchy problem of the first order macroscopic system \cref{sys:limit-1st-order-rho,sys:limit-1st-order-j}.

Note that, the coefficient $C_{f_0,f_1^{\perp}}$ depends upon the bounds of hydrodynamic quantities $\nm{\rho_0,\,\Omega_0}_{\hx{3}}$, and the function $\nm{f_1^{\perp}}_{\hxv{2}(\MV)} $ (which has been determined by \cref{eq:1st-order}, see the proof in \Cref{lemm:f1-f2-perp-esm}). When we consider the Cauchy problem in a general case of $m \ge 2$, $C_{f_0,f_1^{\perp}}$ will depend on $\nm{\rho_0,\,\Omega_0}_{\hx{m+1}}$, and the function $\nm{f_1^{\perp}}_{\hxv{m}(\MV)} $, and in turn, depends on $\nm{\rho_0,\,\Omega_0}_{\hx{m+3}}$. This completes the whole proof of \Cref{lemm:system-f1-esm}.
\qed

%% ------------------------------------------------------------------------ %%
\section*{Acknowledgments}
The author N. Jiang is supported by the grants from the National Natural Foundation of China under contract Nos. 11971360 and 11731008, and also supported by the Strategic Priority Research Program of Chinese Academy of Sciences, Grant No. XDA25010404. The author Y.-L. Luo is supported by grants from the National Natural Science Foundation of China under contract No. 12201220, the Guang Dong Basic and Applied Basic Research Foundation under contract No. 2021A1515110210, and the Science and Technology Program of Guangzhou, China under the contract No. 202201010497. The author T.-F. Zhang is supported by the grants from the National Natural Foundation of China under contract No. 11871203.

%% ------------------------------------------------------------------------ %%
% \bigskip
% \phantomsection
% \addcontentsline{toc}{section}{\refname}

% \bibliographystyle{aomvar}
  % aomplain, abbrv, siamplain, amsalpha
  % elsarticle-num, elsart-num-sort, elsart-num-names, elsart-harv
  % plain, abbrv, alpha, unsrt, unsrtnat
  % amsplain, amsalpha, aomplain, aomalpha
% \nocite{*} % list all references
% \bibliography{duality-ref}

\begin{thebibliography}{10}

\bibitem{ABCD-19mbe}
P.~Aceves-S\'{a}nchez, M.~Bostan, J.-A. Carrillo, and P.~Degond.
  Hydrodynamic limits for kinetic flocking models of {C}ucker-{S}male type.
  \emph{Math. Biosci. Eng.} \textbf{16} (2019), no.~6, 7883--7910.

\bibitem{ABF8-19m3as}
G.~Albi, N.~Bellomo, L.~Fermo, S.-Y. Ha, J.~Kim, L.~Pareschi, D.~Poyato, and J.~Soler.
  Vehicular traffic, crowds, and swarms: from kinetic theory and multiscale methods to
  applications and research perspectives.
  \emph{Math. Models Methods Appl. Sci.} \textbf{29} (2019), no.~10, 1901--2005.

\bibitem{ASR-19b}
D.~Ars\'{e}nio and L.~Saint-Raymond.
  \emph{From the {V}lasov-{M}axwell-{B}oltzmann system to incompressible
  viscous electro-magneto-hydrodynamics. {V}ol. 1}.
  \emph{EMS Monographs in Mathematics}, European Mathematical Society (EMS), Z\"{u}rich, 2019.

\bibitem{BGL-93cpam}
C.~Bardos, F.~Golse, and C.~D. Levermore.
  Fluid dynamic limits of kinetic equations. {II}. {C}onvergence proofs for the {B}oltzmann equation.
  \emph{Comm. Pure Appl. Math.} \textbf{46} (1993), no.~5, 667--753.

\bibitem{BGL-91jsp}
C.~Bardos, F.~Golse, and D.~Levermore.
  Fluid dynamic limits of kinetic equations. {I}. {F}ormal derivations.
  \emph{J. Statist. Phys.} \textbf{63} (1991),
  no.~1-2, 323--344.

\bibitem{BF-17m3as}
N.~Bellomo and F.~Brezzi.
  Mathematical models of self-propelled particles.
  \emph{Math. Models Methods Appl. Sci.}
  \textbf{27} (2017), no.~6, 997--1004.

\bibitem{BBGO-17b}
N.~Bellomo, A.~Bellouquid, L.~Gibelli, and N.~Outada.
  \emph{A quest towards a mathematical theory of living systems}.
  \emph{Modeling and Simulation in Science, Engineering and Technology}, Birkh\"{a}user/Springer,
  Cham, 2017.

\bibitem{BDT-17b}
N.~Bellomo, P.~Degond, and E.~Tadmor (eds.).
  \emph{Active particles. {V}ol. 1. {A}dvances in theory, models, and applications}.
  \emph{Modeling and Simulation in Science, Engineering and Technology}, Birkh\"{a}user/Springer,
  Cham, 2017.

\bibitem{BC-17m3as}
M.~Bostan and J.~A. Carrillo.
  Reduced fluid models for self-propelled particles interacting through alignment,
  \emph{Math. Models Methods Appl. Sci.} \textbf{27} (2017), no.~7, 1255--1299.

\bibitem{BC-20m3as}
M.~Bostan and J.~A. Carrillo.
  Fluid models with phase transition for kinetic equations in swarming.
  \emph{Math. Models Methods Appl. Sci.} \textbf{30} (2020), no.~10, 2023--2065.

\bibitem{BDMa-22sima}
M.~Briant, A.~Diez, and S.~Merino-Aceituno.
  Cauchy theory for general kinetic {V}icsek models in collective dynamics and mean-field limit approximations.
  \emph{SIAM J. Math. Anal.} \textbf{54} (2022), no.~1, 1131--1168.

\bibitem{Caf-80cpam}
R.~E. Caflisch.
  The fluid dynamic limit of the nonlinear {B}oltzmann equation.
  \emph{Comm. Pure Appl. Math.} \textbf{33} (1980),
  no.~5, 651--666.

\bibitem{CFR-10sima}
J.~A. Carrillo, M.~Fornasier, J.~Rosado, and G.~Toscani.
  Asymptotic flocking dynamics for the kinetic {C}ucker-{S}male model.
  \emph{SIAM J. Math. Anal.} \textbf{42} (2010), no.~1, 218--236.

\bibitem{CFTV-10b}
J.~A. Carrillo, M.~Fornasier, G.~Toscani, and F.~Vecil.
  Particle, kinetic, and hydrodynamic models of swarming,
  in \emph{Mathematical modeling of collective behavior in socio-economic and life sciences}.
  \emph{Model. Simul. Sci. Eng. Technol.}, Birkh\"{a}user Boston, Boston, MA, 2010,
  pp.~297--336.

\bibitem{CHL-17chap}
Y.-P. Choi, S.-Y. Ha, and Z.~Li.
  Emergent dynamics of the {C}ucker-{S}male flocking model and its variants,
  in \emph{Active particles. {V}ol. 1. {A}dvances in theory, models, and applications}.
  \emph{Model. Simul. Sci. Eng. Technol.},
  Birkh\"{a}user/Springer, Cham, 2017, pp.~299--331.

\bibitem{CS-07ieee}
F.~Cucker and S.~Smale.
  Emergent behavior in flocks.
  \emph{IEEE Trans. Automat. Control} \textbf{52} (2007), no.~5,
  852--862.

\bibitem{CS-07jjm}
F.~Cucker and S.~Smale.
  On the mathematics of emergence.
  \emph{Jpn. J. Math.} \textbf{2} (2007), no.~1, 197--227.

\bibitem{MEL-89cpam}
A.~De~Masi, R.~Esposito, and J.~L. Lebowitz.
  Incompressible {N}avier-{S}tokes and {E}uler limits of the
  {B}oltzmann equation.
  \emph{Comm. Pure Appl. Math.} \textbf{42} (1989),
  no.~8, 1189--1214.

\bibitem{DFL-15arma}
P.~Degond, A.~Frouvelle, and J.-G. Liu.
  Phase transitions, hysteresis, and hyperbolicity for
  self-organized alignment dynamics.
  \emph{Arch. Ration. Mech. Anal.}
  \textbf{216} (2015), no.~1, 63--115.

\bibitem{DFL-22krm}
P.~Degond, A.~Frouvelle, and J.-G. Liu.
  From kinetic to fluid models of liquid crystals by the moment method.
  \emph{Kinet. Relat. Models} \textbf{15} (2022), no.~3, 417--465.

\bibitem{DFLMN-14swz}
P.~Degond, A.~Frouvelle, J.-G. Liu, S.~Motsch, and L.~Navoret.
  Macroscopic models of collective motion and self-organization,
  in \emph{S\'{e}minaire {L}aurent {S}chwartz---\'{E}quations aux d\'{e}riv\'{e}es
  partielles et applications. {A}nn\'{e}e 2012--2013},
  \emph{S\'{e}min.
  \'{E}qu. D\'{e}riv. Partielles}, \'{E}cole Polytech., Palaiseau, 2014,
  pp.~Exp. No. I, 27.

\bibitem{DFA-17m3as}
P.~Degond, A.~Frouvelle, and S.~Merino-Aceituno.
  A new flocking model through body attitude coordination.
  \emph{Math. Models Methods Appl. Sci.} \textbf{27}
  (2017), no.~6, 1005--1049.

\bibitem{DFAT-18mms}
P.~Degond, A.~Frouvelle, S.~Merino-Aceituno, and A.~Trescases.
  Quaternions in collective dynamics. \emph{Multiscale Model. Simul.}
  \textbf{16} (2018), no.~1, 28--77.

\bibitem{DLMP-13maa}
P.~Degond, J.-G. Liu, S.~Motsch, and V.~Panferov.
  Hydrodynamic models of self-organized dynamics: derivation and existence theory,
  \emph{Methods Appl. Anal.} \textbf{20} (2013), no.~2, 89--114.

\bibitem{DA-20m3as}
P.~Degond and S.~Merino-Aceituno.
  Nematic alignment of self-propelled particles: from particle to macroscopic dynamics,
  \emph{Math. Models Methods Appl. Sci.} \textbf{30} (2020), no.~10,
  1935--1986.

\bibitem{DMVY-19jmfm}
P.~Degond, S.~Merino-Aceituno, F.~Vergnet, and H.~Yu.
  Coupled self-organized hydrodynamics and {S}tokes models for suspensions of active
  particles.
  \emph{J. Math. Fluid Mech.} \textbf{21} (2019), no.~1, Paper No.
  6, 36.

\bibitem{DM-08m3as}
P.~Degond and S.~Motsch.
  Continuum limit of self-driven particles with orientation interaction.
  \emph{Math. Models Methods Appl. Sci.} \textbf{18} (2008), no.~suppl., 1193--1215.

\bibitem{FK-19apde}
A.~Figalli and M.-J. Kang.
  A rigorous derivation from the kinetic {C}ucker-{S}male model to the pressureless
  {E}uler system with nonlocal alignment.
  \emph{Anal. PDE} \textbf{12} (2019),
  no.~3, 843--866.

\bibitem{GSR-04inv}
F.~Golse and L.~Saint-Raymond.
  The {N}avier-{S}tokes limit of the {B}oltzmann equation for bounded collision
  kernels.
  \emph{Invent. Math.} \textbf{155} (2004), no.~1, 81--161.

\bibitem{Guo-06cpam}
Y.~Guo.
  Boltzmann diffusive limit beyond the {N}avier-{S}tokes approximation.
  \emph{Comm. Pure Appl. Math.} \textbf{59}
  (2006), no.~5, 626--687.

\bibitem{HL-09cms}
S.-Y. Ha and J.-G. Liu.
  A simple proof of the {C}ucker-{S}male flocking dynamics and mean-field limit.
  \emph{Commun. Math. Sci.} \textbf{7} (2009), no.~2, 297--325.

\bibitem{HT-08krm}
S.-Y. Ha and E.~Tadmor.
  From particle to kinetic and hydrodynamic descriptions of flocking.
  \emph{Kinet. Relat. Models} \textbf{1} (2008), no.~3, 415--435.

\bibitem{JLZ-18sima}
N.~Jiang, Y.~Liu, and T.-F. Zhang.
  Global classical solutions to a compressible model for micro-macro polymeric fluids near equilibrium.
  \emph{SIAM J. Math. Anal.}
  \textbf{50} (2018), no.~4, 4149--4179.

\bibitem{JL-22apde}
N.~Jiang and Y.-L. Luo.
  From {V}lasov-{M}axwell-{B}oltzmann system to two-fluid incompressible
  {N}avier-{S}tokes-{F}ourier-{M}axwell system with {O}hm's law: convergence
  for classical solutions.
  \emph{Ann. PDE} \textbf{8} (2022), no.~1, Paper No.
  4, 126.

\bibitem{JLT-22axv}
N.~Jiang, Y.-L. Luo, and S. Tang.
  Grad-{C}aflisch type decay estimates of pseudo-inverse of
  linearized {B}oltzmann operator and application to {H}ilbert expansion of
  compressible {E}uler scaling. preprint,
  arXiv:\href{https://arxiv.org/abs/2206.02677}{2206.02677} [math.AP], 2022.

\bibitem{JLZ-20arma}
N.~Jiang, Y.-L. Luo, and T.-F. Zhang.
  Coupled self-organized hydrodynamics and {N}avier-{S}tokes
  models: local well-posedness and the limit from the self-organized
  kinetic-fluid models.
  \emph{Arch. Ration. Mech. Anal.} \textbf{236} (2020),
  no.~1, 329--387.

\bibitem{JLZ-21axv}
N.~Jiang, Y.-L. Luo, and T.-F. Zhang.
  Hydrodynamic limit of the incompressible
  {N}avier-{S}tokes-{F}ourier-{M}axwell system with {O}hm's law from the
  {V}lasov-{M}axwell-{B}oltzmann system: {H}ilbert expansion approach.
  preprint, arXiv:\href{https://arxiv.org/abs/2007.02286v2}{2007.02286v2}
  [math.AP], 2021.

\bibitem{JXZ-16sima}
N.~Jiang, L.~Xiong, and T.-F. Zhang.
  Hydrodynamic limits of the kinetic self-organized models,
  \emph{SIAM J. Math. Anal.} \textbf{48} (2016), no.~5, 3383--3411.

\bibitem{KCBFL-13physD}
T.~Kolokolnikov, J.~A. Carrillo, A.~Bertozzi, R.~Fetecau, and M.~Lewis.
  Emergent behaviour in multi-particle systems with
  non-local interactions [{E}ditorial].
  \emph{Phys. D} \textbf{260} (2013),,
  1--4.

\bibitem{LLZ-07cpam}
F.-H. Lin, C.~Liu, and P.~Zhang.
  On a micro-macro model for polymeric fluids near equilibrium.
  \emph{Comm. Pure Appl. Math.} \textbf{60} (2007), no.~6,
  838--866.

\bibitem{LW-18jfa}
Y.~Liu and W.~Wang.
  The small {D}eborah number limit of the {D}oi-{O}nsager equation without hydrodynamics.
  \emph{J. Funct. Anal.} \textbf{275} (2018), no.~10, 2740--2793.

\bibitem{Majda-84b}
A.~Majda. \emph{Compressible fluid flow and systems of conservation laws in several space variables}.
  \emph{Applied Mathematical Sciences} \textbf{53}, Springer-Verlag, New York, 1984.

\bibitem{MT-11jsp}
S.~Motsch and E.~Tadmor.
  A new model for self-organized dynamics and its flocking behavior.
  \emph{J. Stat. Phys.}
  \textbf{144} (2011), no.~5, 923--947.

\bibitem{Tad-21nAMS}
E.~Tadmor.
  On the mathematics of swarming: emergent behavior in alignment dynamics.
  \emph{Notices Amer. Math. Soc.} \textbf{68} (2021),
  no.~4, 493--503.

\bibitem{Taylor-11b}
M.~E. Taylor.
  \emph{Partial differential equations {III}. {N}onlinear equations}, second ed.
  \emph{Applied Mathematical Sciences}
  \textbf{117}, Springer, New York, 2011.

\bibitem{VCBCS-95prl}
T.~Vicsek, A.~Czir\'{o}k, E.~Ben-Jacob, I.~Cohen, and O.~Shochet.
  Novel type of phase transition in a system of self-driven particles.
  \emph{Phys. Rev. Lett.} \textbf{75} (1995), no.~6,
  1226--1229.

\bibitem{WZZ-15cpam}
W.~Wang, P.~Zhang, and Z.~Zhang.
  The small {D}eborah number limit of the {D}oi-{O}nsager equation to the {E}ricksen-{L}eslie equation.
  \emph{Comm. Pure Appl. Math.} \textbf{68} (2015), no.~8, 1326--1398.

\end{thebibliography}

\end{document}